\documentclass[11pt]{article}

\usepackage{amssymb}
\usepackage{amstext}
\usepackage{amsmath}
\usepackage{latexsym}
\usepackage[all]{xy}

\newfont{\gothique}{eufm10 scaled 1100}  

\newcommand{\PP}{{\mathbb{P}}}
\newcommand{\CC}{{\bf{C}}}
  
\newcommand{\HH}{H^1(\Theta_{X})}
\newcommand{\TET}{\Theta_{X}}
\newcommand{\KX}{K_{X}}
\newcommand{\HKX}{H^0(\OO_X (K_{X}))}
\newcommand{\TT}{{\cal{T}}}
\newcommand{\GG}{{\cal{G}}}

\newcommand{\OM}{\Omega_{X}}

\newcommand{\OO}{{\cal{O}}}

\newcommand{\FF}{{\cal{F}}}

\newtheorem{thm}{Theorem}[section]
\newtheorem{lem}[thm]{Lemma}
\newtheorem{pro}[thm]{Proposition}
\newtheorem{cor}[thm]{Corollary}
\newtheorem{rem}[thm]{Remark}
\newtheorem{cl}[thm]{Claim}

\newlength{\myskip}
\newenvironment{pf}{
     \addvspace{\myskip}  

     \noindent {\it Proof.$\, $}}
     {$\Box$

     \addvspace{\myskip}
     }

\makeatletter
\renewcommand{\@seccntformat}[1]{\S \/ {\csname the#1\endcsname}\hspace{0.5em}}
\makeatother

\title{ON THE INFINITESIMAL TORELLI THEOREM \\
          FOR REGULAR SURFACES WITH VERY AMPLE CANONICAL DIVISOR }
\author{Igor Reider}

\begin{document}

\bibliographystyle{amsplain}
\maketitle

\setcounter{section}{-1}
\numberwithin{equation}{section}
\begin{abstract}
Let $X$ be a smooth compact complex surface subject to the following conditions:
\begin{enumerate}
\item[(i)]
  the canonical
 line bundle $\OO_X(K_X ) $ is very ample,
\item[(ii)]
the irregularity $q(X): = h^1(\OO_X) =0$,
\item[(iii)]
 $X$ contains no rational normal curves of degree
$\leq (p_g-1)$,
\item[(iv)]
the multiplication map $m_2: Sym^2(\HKX) \longrightarrow H^0 (\OO_X (2K_X))$
is surjective.
\end{enumerate}
It is shown that the Infinitesimal Torelli holds for such $X$.

Our proof is based on the study of the cup-product
$$
H^1 (\Theta_X ) \longrightarrow \HKX^{\ast} \otimes H^1 (\OM)
$$
where $\Theta_X$ (resp. $\OM$) is the holomorphic tangent (resp. cotangent) bundle of $X$. Conceptually, the approach consists of lifting the data of the cohomological cup-product above to the category 
 of complexes of coherent sheaves of $X$. This establishes connections between the geometry of the canonical map and the above cup-product by exhibiting geometrically meaningful objects in the category of (short) exact complexes of coherent sheaves on $X$.

\end{abstract}  
\section{Introduction}  
The classical Torelli theorem says that a smooth projective curve of genus 
$\geq2$ is determined, up to an isomorphism, by its Jacobian and its theta-divisor
(see, e.g. \cite{GH}). The works of Griffiths on the variation of Hodge structure
 allow
one to formulate the Torelli question in arbitrary dimension. Namely, the question
asks if 
a smooth projective manifold of complex dimension $n$ can be recovered, 
up to an isomorphism,
from its Hodge structure of weight $n$, i.e. from the Hodge decomposition
$\displaystyle{H^n(X,{\bf C})= \bigoplus_{p+q=n} H^{p,q} (X)}$
(we refer to \cite{G} and references therein for more details).
The infinitesimal version of the Torelli question arises naturally as follows.
Let $\pi: {\cal X} \longrightarrow B$ be a proper surjective smooth morphism of 
smooth algebraic varieties with $\pi^{-1}(b)= X_b$ being projective manifold of 
complex dimension $n$. Fixing a reference point $o\in B$ we can view the morphism
$\pi$ as a deformation of complex structure on $X_o = \pi^{-1}(o)$. 
This induces the variation
of Hodge structure on $H^n(X_o,{\bf C})$. Following Griffiths one 
obtains the period map 
$
P_n: B \longrightarrow D_n / \Gamma
$,
where $D_n$ is the Griffiths period domain for Hodge structures of weight $n$
on $H^n(X_o,{\bf C})$ and $\Gamma$ is the monodromy group of the family
$\{X_b \}_{b\in B}$, i.e. the image of the representation of the fundamental
group $\pi_1 (B,o)$ on $H^n(X_o,{\bf C})$ (see \cite{G}, \cite{GS}).
The Infinitesimal Torelli question asks if the derivative of the period map
$P_n$ is injective. More precisely, one knows (\cite{GS}) that the period map
is locally liftable to $D_n$, i.e. for any point $b \in B$ there exists a
 neighborhood $U_b$ and a morphism $\tilde{P}_n: U_b  \longrightarrow D_n$
 such that the diagram
 $$
 \xymatrix@R=12pt@C=12pt{
 & D_n \ar[2,0] \\
 U_b \ar[ur]^{\tilde{P}_n} \ar[dr]_{P_n} & \\
 &  D_n / \Gamma }
 $$
 commutes. Then the Infinitesimal Torelli question is about the injectivity
 of the differential of $\tilde{P}_n$
 $$
  (d \tilde{P}_n)_b : T_{B,b} \longrightarrow T_{D_n, \tilde{P}_n (b)}
  $$
  where $T_{B,b}$ (resp. $T_{D_n, \tilde{P}_n (b)}$) is the holomorphic
  tangent space  of $B$ (resp. $D_n$) at $b$ (resp. $\tilde{P}_n (b)$).
  From the work of Griffiths one knows that the image of 
  $(d \tilde{P}_n)_b$ is contained in the subspace
  $\displaystyle{\bigoplus_{p+q=n} Hom(H^{p,q}(X_b), H^{p-1,q+1}(X_b))}$ of 
  $T_{D_n, \tilde{P}_n (b)}$ (Griffiths' transversality of the period map)
  while Kodaira-Spencer theory of deformation of complex structure gives
  the linear map $T_{B,b} \longrightarrow H^1( \Theta_{X_b})$,
  where $\Theta_{X_b}$ is the holomorphic tangent bundle of $X_b$.
  Furthermore, Griffiths shows that the following diagram commutes
  $$
 \xymatrix@R=12pt@C=12pt{
  T_{B,b} \ar[d] \ar[r]^(0.22){(d \tilde{P}_n)_b}& 
  {\displaystyle\bigoplus_{p+q=n} Hom(H^{p,q}(X_b), H^{p-1,q+1}(X_b))} \ar@{=}[d] \\
  H^1( \Theta_{X_b}) \ar[r]^(0.25){p_n}& 
 {\displaystyle \bigoplus_{p+q=n} Hom(H^{p,q}(X_b), H^{p-1,q+1}(X_b))} }
  $$
  where the homomorphism at the bottom is given by the cohomology cup-product
  $$
  H^1( \Theta_{X_b}) \otimes H^q( \Omega^p_{X_b})\longrightarrow
   H^{q+1}( \Omega^{p-1}_{X_b})
   $$
   where $\Omega_{X_b}$ is the holomorphic cotangent bundle of $X_b$ and
   $\Omega^p_{X_b}$ is its $p$-th exterior power and where
   the identification $H^{p,q}(X_b) = H^q( \Omega^p_{X_b})$ (via
   Dolbeault isomorphism) has been used.
   This cohomological interpretation allows one to reformulate
   the Infinitesimal Torelli question as the question about the injectivity of 
   the homomorphism
   \begin{equation}\label{pn}
   p_n:  H^1( \Theta_{X_b}) \longrightarrow 
 \bigoplus_{p+q=n} Hom( H^q( \Omega^p_{X_b}), H^{q+1}( \Omega^{p-1}_{X_b}))
 \end{equation}
 This cohomological interpretation turned out to be tractable in certain cases.
 In particular, it has given rise to the Infinitesimal Torelli theorem
 for hypersurfaces of high degree in an arbitrary projective variety, \cite{Gr}
 (see also \cite{Fl}). For smooth projective surfaces, i.e. $n=2$,
  R. Pardini proved the 
 Infinitesimal Torelli theorem for smooth abelian covers with a "general" building
 data for the abelian cover, \cite{Pa} (see also \cite{Pe}). Our previous works
 settled the question for a large class of irregular surfaces of general type,
 \cite{R1}, and for bicanonical double coverings, i.e.
double coverings branched along a smooth divisor in $|2\KX |$, \cite{R2}.

 One of the results of this paper is the Infinitesimal Torelli theorem
  for 
   any smooth complex projective
 surface $X$ subject to the following conditions:
 \begin{eqnarray}\label{cond}
 (i)& &{\mbox{ the canonical line bundle $\OO_X(K_X)$ of $ X$
   is very ample, }} \nonumber\\
 (ii)& &{\mbox{ the irregularity $q(X)= h^1(\OO_X)= h^0 (\OM)=0$}}, \\
(iii)& &
{\mbox{ $X$ contains no rational normal curves of degree
$\leq (p_g-1)$,}} \nonumber\\
(iv)& &{\mbox{the canonical line bundle is quadratically normal, i.e. the
multiplication map}} \nonumber\\
& &
{\mbox{$m_2:Sym^2 H^0(\OO_X (\KX) \longrightarrow H^0 (\OO_X (2\KX))$
is surjective.}}
\nonumber
\end{eqnarray}

\begin{thm}\label{th}
Let $X$ be subject to (\ref{cond}).
Then the Infinitesimal Torelli theorem holds for $X$.
 More precisely, for $X$ subject to  (\ref{cond})
  the cup-product
\begin{equation}\label{cp}
p_X: H^1(\Theta_{X}) \longrightarrow \left(H^0(\OO_X(K_{X}))\right)^{\ast}
\otimes H^1 ( \Omega_{X})
\end{equation}
is injective.
\end{thm}

The above result could be viewed as a precise generalization of the analogous result for curves, since for a curve $C$ the differential of the period map is
encoded in the cup-product
$$
H^1(\Theta_{C}) \longrightarrow \left(H^0(\OO_X(K_{C}))\right)^{\ast}\otimes H^1 ( \OO_{C})
$$
and its injectivity is equivalent to the quadratic normality of the canonical embedding, since $H^1(\Theta_{C})=H^1 (\OO_C (-K_C))$ and the above map is essentially the dual of the multiplication
$m_2:Sym^2 (H^0(\OO_X(K_{C}))\longrightarrow H^0(\OO_X(2K_{C}))$.  The quadratic normality of the canonical embedding of curves is a classical result of Max Noether, see e.g. \cite{GH}. For canonical surfaces the analogue of
Noether's theorem is not known. Hence the reason for a rather unpleasant assumption $(iv)$ in \eqref{cond}. That assumption can be replaced by a weaker one of a cohomological nature. This way 
the above theorem is a corollary of a more general result. To state it we will
need to recall another cohomological cup-product
\begin{equation}\label{q-form}
\delta_X :Sym^2 (H^1(\Theta_X)) \longrightarrow \mbox{$H^2 (\bigwedge^2 \Theta_X)$} =H^2(\OO_X (-\KX)).
\end{equation}
We view this as a vector-valued quadratic form on $H^1(\Theta_X)$. In particular, we will be interested in the quadratic cone
\begin{equation}\label{q-cone}
\widehat{Q}_X =\{\xi\in H^1(\Theta_X) |\,\,\delta_X(\xi^2)=0 \}
\end{equation}
of isotropic (with respect to $\delta_X$) vectors in $H^1(\Theta_X)$. Our main result can be stated as follows.
\begin{thm}\label{th-tech} 
Let $X$ be a smooth compact complex surface subject to the conditions $(i)-(iii)$ in \eqref{cond}. Then the intersection  of the kernel \,$\ker(p_X)$ of the cup-product \eqref{cp}
and the quadratic cone $\widehat{Q}_X $ is the zero vector: 
$$
\ker(p_X) \bigcap \widehat{Q}_X =\{0 \}.
$$
In particular, the kernel $\ker(p_X)$ of the cup-product \eqref{cp} is a proper subspace of
$H^1 (\Theta_X)$, provided the cone $\widehat{Q}_X$ in \eqref{q-cone} is nonzero.
\end{thm}


Before going further into explaining the main ideas of the paper, let us spell out the reason for Theorem \ref{th} to be a consequence of Theorem \ref{th-tech}. The point is that under the assumption
of quadratic normality $(iv)$ in \eqref{cond}, the kernel $ker(p_X)$ becomes a totally isotropic subspace of $H^1(\Theta_X)$, i.e.
$ker(p_X) \subset \widehat{Q}_X$. This is due to the fact that the quadratic form $\delta_X$ in \eqref{q-form} and the dual of the multiplication map
$$
m^{\ast}_2 : H^2 (\OO_X (-K_X)) \cong ( H^0(\OO_X (2K_X)))^{\ast}  \longrightarrow (Sym^2 H^0(\OO_X (K_X))^{\ast}
$$
 are parts of the following commutative triangle
\begin{equation}\label{diag-comtriang}
\xymatrix@R=12pt@C=8pt{
Sym^2(H^1 (\TET)) \ar[rr]^{\delta_X}  \ar[rd]_(0.4){p^{(2)}}& &
H^2(\OO_X (-K_X))\ar[dl]^(0.45){m^{\ast}_2}\\
& ( Sym^2(\HKX))^{\ast} &
}
\end{equation}
where $p^{(2)}$ is the (second) iterate of the cup-product in \eqref{cp}. We recall that the dual of the above triangle appears in one of the Griffiths' invariants of the infinitesimal variation of Hodge structure (IVHS), see \cite{G}. Namely, the kernel of $(p^{(2)})^{\ast}$ is a system of quadrics in $\PP(H^0(\OO_X (K_X))^{\ast})$ naturally attached to the IVHS of a canonical surface and it contains the space of quadrics containing the canonical image of $X$. Griffiths' insight here is that $\ker((p^{(2)})^{\ast})$ should have a geometric meaning as well. Our considerations suggest that the kernel of the quadratic form $\delta_X$ and, in particular, isotropic Kodaira-Spencer classes are related to the geometry of the canonical map.  


We will now give a detailed outline of the main ideas and steps involved in the proof of Theorem \ref{th-tech}.

\vspace{0.2cm}  
{\bf Outline of the proof of Theorem \ref{th-tech}.}
Our approach toward the study of the cup-product in (\ref{cp}) is based on 
interpreting the cohomology classes of $H^1(\Theta_{X})$ as objects in the category of (short) exact complexes of coherent sheaves on $X$. Namely, we make the following identification
\begin{equation}\label{ident}
 \HH=Ext^1 (\OO_X(\KX), \OM),
\end{equation}
according to which a cohomology class $\xi \in \HH$ can be thought of as
 the corresponding extension, i.e., an exact sequence of sheaves on $X$
\begin{equation}\label{ext}
 \xymatrix{
0\ar[r]& \OM \ar[r]^{i}&\TT_{\xi} \ar[r]^(0.35){p}& \OO_X(\KX)\ar[r]&0 }.
\end{equation}

We fix a nonzero $\xi$ lying in the kernel of  (\ref{cp}), then $p$ induces a surjective
homomorphism
\begin{equation}\label{p}
 \xymatrix{
H^0 (\TT_{\xi} ) \ar[r]^(0.4){p^0}& H^0 (\OO_X (\KX)) }
\end{equation}
due to the fact that the coboundary map
$ H^0 (\OO(\KX)) \longrightarrow H^1 (\OM)$ in the long exact sequence of cohomology groups of (\ref{ext}) is the cup-product with $\xi$, which we assume to be identically zero. This together with the assumption $(ii)$ of (\ref{cond})
imply that (\ref{p}) is an isomorphism.
 Thus the fact that $\xi$ lies in the kernel of the cup-product (\ref{cp})
means that the sheaf $\TT_{\xi}$ in (\ref{ext}) has global sections parametrized by
$H^0 (\OO(\KX))$. This will be recorded by introducing the isomorphism
\begin{equation}\label{alpha}
 \xymatrix{
\alpha: H^0 (\OO_X (\KX)) \ar[r]& H^0 (\TT_{\xi} ),  }
\end{equation}
which is the inverse of $p^0$ in (\ref{p}). The main guiding line of the subsequent arguments is to use the abundance of global sections of $\TT_{\xi}$ to construct a sort of `destabilizing' subobjects of the extension sequence  \eqref{ext}.
This involves the study of the sheaf $\TT_{\xi}$.

\vspace{0.2cm}
{\it Step 1.}
As a first step toward understanding the properties of $\TT_{\xi}$ we show
that it is generated by its global sections almost everywhere.

\begin{lem}\label{lem-glgen}
$\TT_{\xi}$ is generically generated by its global sections. 
\end{lem}

The main point in the proof of this lemma is to rule out the case when the global sections of  $\TT_{\xi}$ generate a subsheaf of rank $2$. Our observation is that in such a situation $X$ admits a foliation. By that we mean that the holomorphic tangent bundle $\Theta_X$ fits into the following exact sequence
\begin{equation}\label{foliationseq}
 \xymatrix{
0\ar[r]& \OO_X (D-\KX) \ar[r]& \TET  \ar[r]& {\cal I}_A (-D) \ar[r]&0,  }
\end{equation}
where $\OO_X (D)$ is a line bundle intrinsically attached to $\xi$ and ${\cal I}_A$ is the sheaf of ideals of some $0$-dimensional subscheme $A$ of $X$. Furthermore, we observe that the cohomology class $\xi$ must come from some cohomology class $\xi'$ under the homomorphism
$$
H^1( \OO_X (D-\KX)) \longrightarrow H^1(\TET)
$$
induced by the monomorphism in the above exact sequence. It turns out that even more is true:
the class $\xi'$ comes from a cohomology class of a sufficiently negative line bundle.
Namely, we show that for every integer $m\geq 1$ and any nonzero global section
$\psi \in H^0 (\OO_X (mK_X))$, the cohomology class $\xi'$ lies in the image  
of the homomorphism
$$
H^1( \OO_X (D-(m+1)\KX)) \stackrel{\psi}{\longrightarrow} H^1( \OO_X (D-\KX))
$$
induced by the multiplication by $\psi$. Of course, for $m$ sufficiently big
the cohomology group on the left vanishes
. Hence the vanishing of $\xi$ and the conclusion that $\TT_{\xi}$ is generically globally generated.

We wish to point out that though the argument seems to turn on the unraveling the cohomological properties of the class $\xi$, the essential observation is that the failure of generic global generation of $\TT_{\xi}$ produces a nontrival subsheaf of $\TT_{\xi}$. That subsheaf is a part of an extension sequence which produces a destabilizing subobject of the extension sequence \eqref{ext}.
 It is the properties of this subobject that impose the cohomological restrictions on $\xi$.

There is even stronger generic generation criterion for $\TT_{\xi}$.
\begin{pro} \label{pro:gencr}
Let $W$ be a subspace of $\HKX$ of dimension at least $3$. Then the subspace
$\alpha(W) \subset H^0 (\TT_{\xi})$ generically generates $\TT_{\xi}$, provided
the linear subsystem $|W|\subset |\KX|$ has at most $0$-dimensional base locus.
\end{pro}
The result gives a clear indication of the relation between the geometry of $X$
and the properties of the extension \eqref{ext}. But it also tells us that it takes more than abundance of global sections of $\TT_{\xi}$ to construct  geometrically interesting subsheaves of smaller rank in $\TT_{\xi}$.

\vspace{0.2cm}
\noindent
{\it Step 2.}  
A new ingredient for  producing geometrically interesting subsheaves of $\TT_{\xi}$ and its second exterior power $\bigwedge^2 \TT_{\xi}$ is the condition of isotropy of $\xi$.
 To explain this we use the parametrization $\alpha: \HKX \longrightarrow H^0 (\TT_{\xi})$ in \eqref{alpha} to define the `higher' order products. Namely, for $k=2,3,$ we consider the linear maps
$$
\mbox{$\alpha^{(k)}: \bigwedge^k \HKX \longrightarrow  H^0 (\bigwedge^k \TT_{\xi})$}
$$
defined by the composition
$$
\mbox{$\bigwedge^k \HKX \stackrel{\wedge^k \alpha}{\longrightarrow} \bigwedge^k H^0 ( \TT_{\xi})\longrightarrow H^0 (\bigwedge^k \TT_{\xi}),$}
$$
 i.e. for any $k$-tuple $\phi_1,\ldots \phi_k \in \HKX$ one defines
$\alpha^{(k)} (\phi_1,\ldots \phi_k)$ as the image of the exterior product 
$\bigwedge^k_{s=1} \alpha(\phi_s)$ under the natural homomorphism
$\bigwedge^k H^0 ( \TT_{\xi})\longrightarrow H^0 (\bigwedge^k \TT_{\xi})$.
The two products are related by the formula
\begin{equation}\label{detformula}
\phi''\alpha^{(2)}(\phi,\phi') -\phi'\alpha^{(2)}(\phi,\phi'')+\phi \alpha^{(2)}(\phi',\phi'')=\alpha^{(3)}(\phi,\phi',\phi'') j, \,\,\forall \phi,\phi',\phi'' \in \HKX,
\end{equation}
where $j$ is the global section of $\bigwedge^2 \TT_{\xi} (-\KX)$ corresponding
to the exterior product
$\wedge^2 i: \OO_X (\KX) \longrightarrow \bigwedge^2 \TT_{\xi} $ of the monomorphism $i$ in the extension sequence \eqref{ext}. The above formula is just the expansion of the determinant of a $3\times 3$ matrix with respect to one of its columns. 
The formula also tells us that
 $\alpha^{(3)}: \bigwedge^3 \HKX \longrightarrow  H^0 (\bigwedge^3 \TT_{\xi}) = H^0 (\OO_X (2\KX))$ is a Koszul cocycle. Furthermore, it represents the cohomology class $\delta_X (\xi^2) \in H^2 (\OO_X(-\KX))$. So the assumption that $\xi$ is isotropic, i.e., $\delta_X (\xi^2)=0$, means that $\alpha^{(3)}$ is a Koszul coboundary: there exists a linear map
$$
\mbox{$l: \bigwedge^2 \HKX \longrightarrow \HKX,$}
$$
such that
{\small
\begin{equation}\label{coboundary-intro}
\alpha^{(3)}(\phi,\phi',\phi'')=d_{Kosz} (l)(\phi,\phi',\phi'')= \phi'' \, l(\phi,\phi') -\phi' \,l(\phi,\phi'')+\phi
\,l(\phi',\phi''), \,\,\forall \phi,\phi',\phi'' \in \HKX.
\end{equation}
}
Combining this with the determinantal formula in \eqref{detformula} gives 
\begin{equation}\label{Koszrel-intro}
\phi''(\alpha^{(2)}(\phi,\phi') -l(\phi,\phi')j)-\phi'(\alpha^{(2)}(\phi,\phi'')-l(\phi,\phi'')j)+\phi (\alpha^{(2)}(\phi',\phi'') -l(\phi',\phi'')j)=0, 
\end{equation}
for all $\phi,\phi',\phi'' \in \HKX$. This means that for any triple
$\phi,\phi',\phi'' \in \HKX$, the global sections $\{\alpha^{(2)}(\phi,\phi') -l(\phi,\phi')j, \alpha^{(2)}(\phi,\phi'')-l(\phi,\phi'')j, \alpha^{(2)}(\phi',\phi'') -l(\phi',\phi'')j\}$ of $\bigwedge^2 \TT_{\xi} $ fail to generate that sheaf everywhere. This is a source of producing  nontrivial subsheaves of $\bigwedge^2 \TT_{\xi} $ whose properties turn out to contradict the generic global generation criterion in Proposition \ref{pro:gencr}.

Indeed, under  
the additional technical assumption that $\TT_{\xi}$ is generated by its global sections outside of a subscheme of codimension $2$
  it is easy to arrive to a contradiction. Namely, that assumption guaranties that the nonzero divisors
$(\alpha^{(3)}(\phi,\phi',\phi'')=0)$, as $\phi,\phi',\phi''$ vary in $\HKX$, move in a linear system with no fixed part and a general divisor of that form is irreducible.
On the other hand the relation 
\eqref{Koszrel-intro} implies that the global section
$$
\tau(\phi,\phi',\phi'')=l(\phi,\phi')\alpha(\phi'') -l(\phi,\phi'')\alpha(\phi') +l(\phi',\phi'')\alpha(\phi)
$$
of $\TT_{\xi} (\KX)$ must have $1$-dimensional zero locus. In view of the coboundary relation \eqref{coboundary-intro} that $1$-dimensional zero locus must be a proper component of the divisor $(\alpha^{(3)}(\phi,\phi',\phi'')=0)$ and this contradicts the irreducibility of such a divisor for a general triple $\phi,\phi',\phi''$. 

It should be remarked that it is not so much the isotropy condition on
$\delta_X (\xi^2)$, but the Koszul coboundary relation \eqref{coboundary-intro} 
that really matters in our argument. Of course, in the situation just discussed the two conditions are equivalent.

To deal with the general situation, we produce a modified extension sequence
\begin{equation}\label{ext-modif}
\xymatrix{
0\ar[r]&{\cal P} \ar[r]& \TT'\ar[r]& \OO_X (\KX)\ar[r]&0
}
\end{equation} 
which is a subobject of the extension \eqref{ext} and such that the sheaf $\TT'$ in the middle has the following properties:

- $\TT'$ is a modification of $\TT_{\xi}$ along an effective divisor $E'$,

- $H^0 (\TT') \cong \HKX$,

- $\TT'$ is generated by its global sections outside of a codimension $2$ subscheme.


\vspace{0.2cm}
\noindent
 So replacing our original extension by \eqref{ext-modif} would seem to bring our argument to a conclusion.
However, there is an additional difficulty: the new extension corresponds to a cohomology class in $H^1 ({\cal P} (-\KX))$ via the identification
$Ext^1 (\OO_X (\KX), {\cal P}) \cong H^1 ({\cal P} (-\KX))$ and the isotropy condition imposed on $\xi$ may no longer hold for that class. But we recall that what really matters to us is the Koszul coboundary relation. The situation here is as follows: the extension \eqref{ext-modif} also comes with a parametrization
$\alpha': \HKX \longrightarrow H^0 (\TT')$ which in its own turn provides the higher order products
$$
\mbox{$\alpha'^{(k)}: \bigwedge^k \HKX \longrightarrow H^0 (\bigwedge^k \TT'), \,\,k=2,3.$}
$$
In particular, $\alpha'^{(3)} : \bigwedge^3 \HKX \longrightarrow H^0 (\bigwedge^3 \TT')= H^0 (\OO_X (2\KX -E'))$ is a Koszul cocycle and it is related to $\alpha^{(3)}$ by the formula
$$
\alpha^{(3)} =e'\alpha'^{(3)},
$$
where $e'$ is a global section defining the divisor $E'$, i.e., $E'= (e'=0)$.
Hence, the coboundary condition \eqref{coboundary-intro} for $\alpha^{(3)}$ gives
$$
d_{Kosz} (l) =e'\alpha'^{(3)},
$$
a `twisted' coboundary condition for $\alpha'^{(3)}$. This relation indicates that the Koszul cochain $l: \bigwedge^2 \HKX \longrightarrow \HKX$ should `see' the divisor $E'$. 

At this point one should take into account that $l$ is not unique and can be deformed by adding any Koszul cochain $d_{Kosz} (f)$, where $f:\HKX\longrightarrow \CC$ is a linear form. Ideally, one would wish to deform $l$ to a cochain  `divisible' by $e'$, i.e., one wishes to have a linear function $f:\HKX\longrightarrow \CC$ such that
  $$
l =e' m +d_{Kosz} (f),
$$
 where $m:\bigwedge^2 \HKX \longrightarrow  H^0 (\OO_X (\KX -E'))$. In such a situation $\alpha'^{(3)} = d_{Kosz} (m)$ becomes a Koszul coboundary and then we are done. It turns out that the above relation holds, up to multiple components of
$E'$. In other words 
 $$
l =e'_1 m +d_{Kosz} (f),
$$
where the divisor $E'_1 =(e'_1 =0)$ is a component of $E'$ such that the corresponding reduced divisors $(E'_1)_{red}$ and $E'_{red}$ coincide, and the Koszul cochain $m$
takes its values in $H^0 (\OO_X (\KX -E'_1))$ and its image $Im(m)$ determines the linear subsystem $|Im(m)| \subset |\KX -E'_1|$ without fixed part.

  The above deformation does not seem to come for free. In our argument it is
 the condition $(iii)$ in
\eqref{cond} that permits it.
Once this is done, the new coboundary relation becomes
$$
d_{Kosz} (m) =e'_0 \alpha'^{(3)}
$$
where $E'_0=(e'_0=0)=E'-E'_1$. It turns out that the presence of the residual divisor (if nonzero) $E'_0$ interferes mildly\footnote{though some technical issues need to be overcome - this is the contents of \S6.} and we can complete our argument as though $\alpha'^{(3)}$ is a coboundary.
\vspace{0.2cm}

\noindent
{\bf Concluding remarks and speculations.}
Conceptually, our proof consists of realizing Kodaira-Spencer classes lying in the kernel of the cup-product (\ref{cp}) on the level of the category of complexes of coherent sheaves on $X$. In this our argument is an instance of the general theme so aptly summarized by R.P. Thomas in \cite{T} by the slogan
\begin{center}
{\it `` Complexes good, (Co)homology bad.''}
\end{center}
We believe that the power of this approach resides in the consideration of pairs $(\TT_{\xi}, \alpha(\phi))$ of the extension sheaf $\TT_{\xi}$ with its nonzero global sections $\alpha(\phi)$ and this should relate to such topics as Hilbert schemes of points of $X$ and nonabelian Hodge theory in the spirit of Simpson.  This happens because the family of pairs $(\TT_{\xi}, \alpha(\phi))$ give rise to a family of sheaves of rank $2$ on $X$ having regular global sections\footnote{global sections with $0$-dimensional zero locus.}  - an instance of a so called nonabelian Jacobian of $X$ whose study was initiated in \cite{R3}, see \cite{R4} for an overview. This will be elaborated elsewhere. Closer to the theme of the paper,
we suggest that this approach could be useful in revealing the geometry of canonically polarized surfaces
hidden in the cup-products \eqref{cp} and \eqref{q-form}. Namely, to a nonzero Kodaira-Spencer class $\xi$ one naturally attaches the subspace 
\begin{equation}\label{Wxi}
W_{\xi} := ker\left( \HKX \stackrel{\xi}{\longrightarrow} H^1 (\OM) \right)
\end{equation}
of $\HKX$. The proof of Theorem \ref{th-tech} 
should give a precise criterion for
the linear subsystem $|W_{\xi}| \subset |\KX|$ to have base points, whenever $\xi$ is isotropic. Such a criterion could be useful for a better understanding of the canonical map of surfaces and for studying the moduli of canonical surfaces. It is also a step toward recovering a canonical surface $X$ from its IVHS, since the criterion would tell us that IVHS of $X$ recovers all multi-secant linear subspaces in the canonical embedding of $X$, unless $X$ has a special geometry.
These matters will  be discussed in the sequel to this paper. But what should be hopefully clear already is that our approach gives a substantial evidence of a rich geometry contained in Kodaira-Spencer classes.

\vspace{0.2cm}
To summarize the above discussion, one could say that the Kodaira-Spencer classes
which are isotropic with respect to the quadratic cup-product
\begin{equation}\label{qf-intro}
\delta_X: Sym^2(H^1 (\TET) )\longrightarrow H^2 (\det(\TET))=H^2(\OO_X (-K_X))
\end{equation}
carry geometric information. The `visible' part of this geometric information is either the base locus of the linear subsystem $|W_{\xi}|$ or a very particular geometry of $X$. As we pointed out in the outline of the proof, this geometric part emerges through constructions of subobjects of the extension sequence \eqref{ext} (or its modified version \eqref{ext-modif}).
So we suggest that the `hidden', perhaps more substantial, part as well as the proper context of our approach and its possible generalizations are in the realm of the derived category of coherent sheaves on $X$ with some kind of stability conditions \`a la Bridgeland, inducing a `wall' structure on $\PP(H^1(\TET))$. From this point of view the   
 quadratic cup-product $\delta_X$ seems to be of independent interest. Our results suggest that
 the totally isotropic\footnote{with respect to the quadratic form $\delta_X$ in \eqref{qf-intro}.} subspaces of $H^1(\TET)$ might be a part of this hypothetical `wall' structure of $\PP(H^1(\TET))$. More generally, the subvariety of 
$\PP(H^1(\TET))\times \PP(H^1(\TET))$ parametrizing orthogonal (with respect to $\delta_X$) pairs should be of interest from algebro-geometric as well as derived categorical points of view. It should be noticed that from the diffeo-geometric perspective, the interest of this orthogonality condition in
$H^1(\TET)$ was pointed out long time ago by Y.-T. Siu in \cite{S}, where he proves that the Weil-Petersson metric on the moduli space of canonically polarized complex compact manifolds  has negative holomorphic bisectional curvature in the directions $\xi, \eta \in H^1(\TET)$ which are orthogonal with respect to $\delta_X$.
 
\vspace{0.2cm} 
Finally, we would like to comment on the set of assumptions \eqref{cond} under which Theorem \ref{th-tech} holds.

The condition of the canonical bundle $\OO_X (K_X)$ being ample and globally generated is indispensable in the proof of Lemma \ref{lem-glgen} in the arguments
involving 

\vspace{0.2cm}
- the  construction of the destabilizing extension \eqref{ext-G} in the proof of Lemma \ref{ext-incl} and the modified extension \eqref{ext-modif},

- the semistability of $\OM$ with respect to $\KX$,

- the vanishing of $H^1 (\OO_X (D-mK_X))$, for $m$ sufficiently large.

The condition of very ampleness of $\KX$ is only used in the last stages of the proof, in 
Lemma \ref{l:zeroloc-tauPi}. So we expect
 Theorem \ref{th-tech} to hold for regular surfaces with $\OO_X (\KX)$ ample, globally generated and $(iii)$ in \eqref{cond} replaced by the condition

\vspace{0.2cm}
$(iii)'$ the canonical image of $X$ contains no rational normal curves of degree $\leq (p_g -1)$.

\vspace{0.2cm}
\noindent

The condition of the irregularity $q(X) = 0$ in  (\ref{cond}), $(ii)$, is important in establishing the isomorphism $\alpha$ in (\ref{alpha}) as well as at numerous points where the surjectivity of the homomorphism
$\xymatrix{
\HKX \ar[r]& H^0 (\OO_C (K_X))}
$
of the restriction of global sections of $\OO_X (K_X)$ to divisors $C$ in the
canonical linear system $|K_X |$ is used. 

In fact, for irregular surfaces with very ample canonical bundle, the work \cite{Ga-Z} gives
examples where the injectivity of the cup-product (\ref{cp}) fails.\footnote{we are grateful to F. Catanese for pointing out this article.} On the conceptual level, the (Infinitesimal) Torelli problem for irregular surfaces of general type requires the consideration
of the variation of Hodge structure of weight 2 as well as of weight 1. This is the approach taken in our earlier work \cite{R1} which already
has placed an accent on a use of higher rank bundles in the context of Torelli problem(s). In that work, under the assumption that $\OM$ is generated by its global sections, we were able to trace the failure of the Infinitesimal Torelli theorem to a special geometric property of zero-loci of global sections of $\OM$, a sort of hyperellipticity phenomenon for irregular surfaces of general type. From this point of view, our present work follows a similar logic by making use of the higher rank bundle $\TT_{\xi}$ with enough well-behaved sections (this is the essence of Lemma \ref{lem-glgen}). Furthermore, the choice of $\TT_{\xi}$ is quite natural since if $X$ is a fibre of a one-parameter deformation corresponding to $\xi$ and $\mathfrak{X}$ is the total space of the deformation, then $\TT_{\xi}$ is the restriction of $\bigwedge^2 \Omega_{\mathfrak{X}}$ to $X$.

The assumption $(iii)$ in \eqref{cond} is purely technical: it allows to deform the cochain $l$ in the coboundary condition for $\alpha^{(3)}$, see \eqref{coboundary-intro},
to take into account the geometry stemming from replacing the original extension sequence \eqref{ext} by its subobject - the extension \eqref{ext-modif}.

The condition of isotropy of $\xi$ is a cohomological substitute of the Noether's property of quadratic normality of the canonical bundle, the assumption $(iv)$ in \eqref{cond}. In our argument it has a flavor of a second order condition on the extension sequence \eqref{ext} which allows to `destabilize' that sequence. So the proper understanding of this condition seems to be related to stability conditions in the sense of Bridgeland.

We hope that our approach and some of the ideas of this paper could be used for a better understanding of the moduli spaces of surfaces of general type, higher dimensional manifolds fibred in polarized surfaces and for the Infinitesimal Torelli problem for smooth compact canonical varieties of dimension $\geq 3$.

\vspace{0.2cm}
\noindent
{\bf Organization of the paper.} 
The paper is organized as follows.

In \S1 we prove Lemma \ref{lem-glgen} and Lemma \ref{lem-glgenWxi}, a somewhat more general version of Proposition \ref{pro:gencr}.

\S2 studies the locus where $\TT_{\xi}$ fails to be globally generated. As a result one defines a modification of $\TT_{\xi}$ along a certain effective divisor $E$, where all global sections of $\TT_{\xi}$ are proportional. That new sheaf  is defined by the extension corresponding to a cohomology class $\eta$ in $H^1 (\Theta_X (-E))$ and it is shown that its general global section has no zeros. This gives rise to a family of rank $2$ bundles $\FF_{[\phi]}$ parametrized by points $[\phi]$ of a Zariski dense open subset of $\PP(\HKX)$. 

One then further modifies to obtain the extension \eqref{ext-modif} (this is the extension \eqref{ext-T'} in Lemma \ref{l:T'eta}). This new modification provides a new family of rank $2$ bundles $\FF'_{[\phi]}$ which are modifications of  $\FF_{[\phi]}$ along a divisor. This modification achieves additional properties      (see Proposition \ref{pro:F'phi-gen-nopencil}):

-- $\FF'_{[\phi]}$   is globally generated outside of a $0$-dimensional subscheme of $X$,

-- for a general nonzero global section $f$ of $\FF'_{[\phi]}$, its locus of zeros $Z_f=(f=0)$ is $0$-dimensional , the linear system
$|{\cal J}_{Z_f} \otimes \det(\FF'_{[\phi]})|$ is fixed part
free and its general member is irreducible ( ${\cal J}_{Z_f}$ denotes the ideal sheaf of $Z_f$).

\S3 is devoted to the study of sheaves $\FF_{[\phi]}$ and $\FF'_{[\phi]}$.

In \S4 the condition of isotropy of the class $\xi$ is given in terms of the properties of global sections of the sheaf $\TT_{\xi}$ - the Koszul coboundary relation \eqref{coboundary-intro}.

In \S5 is taken up the question of deforming the Koszul cochain in the coboundary relation 
as discussed in the outline; the main result of the section is Lemma \ref{l:l-modif}.

\S6 translates the deformed  Koszul coboundary condition obtained in \S5 into geometrically meaningful subsheaves of $\bigwedge^2 \TT_{\xi}$.

\vspace{0.2cm}
\noindent
{\bf Acknowledgments.}  
It is a pleasure to express my gratitude to S. Mori, S. Mukai and to the whole staff of RIMS at Kyoto University, where a part of this work has been completed. The final stages of this work have been done while `en d\'el\'egation' (semester of research) at C.N.R.S. and I  
   thank
that institution
for this opportunity.

\section{Proof of Lemma \ref{lem-glgen}}

We consider the following morphisms of sheaves
\begin{equation}\label{ev}
\xymatrix{
\HKX \otimes \OO_X \ar[r]^(0.55){\alpha}& H^0 (\TT_{\xi})\otimes \OO_X \ar[r]^(0.7){ev}& \TT_{\xi} }
\end{equation}
where the first map is the isomorphism defined by $\alpha$ in (\ref{alpha}) and the second is the evaluation morphism. Our task will be to show that
$\TT_{\xi}$ is generically generated by its global sections or, equivalently, that the evaluation morphism in (\ref{ev}) is generically surjective.

Let $\cal G$ be the subsheaf of $\TT_{\xi}$ defined as the saturation of the image of $ev$ in (\ref{ev}). Since $\OO_X (\KX)$ is generated by its global sections and $\alpha$ in (\ref{ev}) is an isomorphism, it follows that the inclusion
$ {\cal G} \hookrightarrow \TT_{\xi} $ composed with the epimorphism $p$ in 
(\ref{ext}) give an epimorphism

\begin{equation}\label{epim}
\xymatrix{
{\cal G} \ar[r]&\OO_X (\KX) \ar[r]&0}.
\end{equation} 
From this it follows that the rank of ${\cal G}$ is at least $2$ (if 
$\mbox{rk} ({\cal G}) =1$, then the epimorphism above, since  $\GG$ is torsion free, must be an isomorphism;
this isomorphism provides the splitting of the extension sequence (\ref{ext}) and hence
the vanishing of $\xi$).

Assume $\mbox{rk} ({\cal G}) =2$. Then the definition of ${\cal G}$ gives rise to the following exact sequence
\begin{equation}\label{ext-G-1}
\xymatrix@R=12pt@C=12pt{
0\ar[r]& {\cal G} \ar[r]&{\cal T}_{\xi} \ar[r]&{\cal I}_{A} (L) \ar[r]&0},
\end{equation}
where ${\cal I}_{A}$ is the ideal sheaf of a $0$-dimensional subscheme $A$ of $X$ and $\OO_X (L)$ is a line bundle on $X$. The above exact sequence implies that ${\cal G}$ is locally free. This is due to the fact that (\ref{ext-G-1}) exhibits $\GG$ as a second syzygy sheaf and those are locally free on a surface, see \cite{[O-S-S]}, Theorem 1.1.6.
 
Combining the exact sequence (\ref{ext-G-1}) with the defining sequence (\ref{ext}) of $\TT_{\xi}$ gives the diagram
\begin{equation}\label{inc-seq-0}
\xymatrix@R=12pt@C=12pt{
&&0\ar[d]& & \\ 
& &{\cal G} \ar[dr] \ar[d]&
&\\
0\ar[r]& \OM \ar[r]^{i} \ar[dr]&\TT_{\xi} \ar[r]^(0.35){p} \ar[d]& \OO_X (\KX)\ar[r]&0\\
& & {\cal I}_A (L)\ar[d]& &\\ 
   &  & 0& &}
\end{equation}
where the slanted arrow on the top is the epimorphism in (\ref{epim}). In particular, the kernel of this epimorphism is a line subbundle, call it  $\OO_X (D)$, of $\GG$. Thus the above diagram can be completed as follows
\begin{equation}\label{inc-seq}
\xymatrix@R=12pt@C=12pt{
&0\ar[d]&0\ar[d]& & \\ 
0\ar[r]& \OO_X (D) \ar[r] \ar[d]&{\cal G} \ar[r] \ar[d]&
\OO_X (\KX) \ar[r] \ar@{=}[d]&0\\
0\ar[r]& \OM \ar[r]^{i} \ar[d]&\TT_{\xi} \ar[r]^(0.35){p} \ar[d]& \OO(\KX)\ar[r]&0\\
& {\cal I}_A (K_X-D) \ar@{=}[r] \ar[d]& {\cal I}_A (K_X-D)\ar[d]& &\\ 
   & 0 & 0& &}
\end{equation} 
 
The above diagram implies the following properties.
\begin{lem}\label{ext-incl}
1) The cohomology class $\xi$ lies in the image of the homomorphism
$$
H^1 (\OO_X (D-K_X))\longrightarrow H^1 (\TET)
$$
induced by the monomorphism $\OO_X (D-K_X) \longrightarrow \TET$, the dual of the epimorphism in the left column of the diagram \eqref{inc-seq}.

2) The homomorphism in $1)$ admits the factorization
\begin{equation}\label{diag:xi}
\xymatrix@R=12pt@C=12pt{
H^1 (\OO_X (D-K_X)) \ar[rr] \ar[dr]& &H^1 (\TET) \\
& H^1 (\TT^{\ast}_{\xi})\ar[ur]&
}
\end{equation}
In particular, $\xi$ is an isotropic element with respect to the quadratic form
\begin{equation}\label{qf}
\delta_X: Sym^2(H^1 (\TET)) \longrightarrow H^2 (det(\TET))=H^2(\OO_X (-K_X)),
\end{equation}
in \eqref{q-form}, i.e., $ \delta_X(\xi^2)=0$.
\end{lem} 
\begin{pf}
 Dualizing the left vertical sequence  in \eqref{inc-seq} gives
\begin{equation}\label{TET-seq}
\xymatrix@R=12pt@C=12pt{
0\ar[r]&\OO_X (D-K_X) \ar[r]&\TET\ar[r]&{\cal I}_A (-D)\ar[r]&0}
\end{equation}
 and hence the homomorphism
\begin{equation}\label{ext-hom}
\xymatrix@R=12pt@C=12pt{
H^1 (\OO_X (D-K_X)) \ar[r]\ar@{=}[d] &H^1 (\TET)\ar@{=}[d]\\
Ext^1 (\OO_X (K_X),\OO_X (D))& Ext^1 (\OO_X (K_X), \Omega_X)}
\end{equation}
relating two groups of extensions.
Furthermore, the morphism of two horizontal extension sequences in (\ref{inc-seq}) tells us that
the cohomology class $\xi \in H^1 (\TET)$ is the image of the cohomology class
in $H^1 (\OO_X (D-K_X))$ corresponding to the extension sequence
\begin{equation}\label{ext-G}
\xymatrix@R=12pt@C=12pt{
0\ar[r]& \OO_X (D) \ar[r]&{\cal G} \ar[r]&\OO_X (\KX) \ar[r]&0}
\end{equation}
which is the horizontal sequence on the top of the diagram \eqref{inc-seq}. This proves the assertion $1)$ of the lemma.

We now turn to the part 2).  From the diagram \eqref{inc-seq} it follows that the dual of the epimorphism of the left column factors through $\TT^{\ast}_{\xi}$, the dual of $\TT_{\xi}$, to give the diagram
$$
\xymatrix@R=12pt@C=12pt{
\OO_X (D-K_X) \ar[rr] \ar[dr]& &\TET \\
& \TT^{\ast}_{\xi}\ar[ur]&
}
$$
Hence the first assertion in $2)$. To see the assertion that $\xi$ is isotropic
we observe that the slanted arrow on the right in the diagram \eqref{diag:xi}
is part of the long exact sequence of cohomology groups arising from the dual
of our extension sequence, the dual of the middle row in \eqref{inc-seq}. Namely, we have
\begin{equation}\label{xi-ort}
\xymatrix@R=12pt@C=12pt{
0\ar[r]&H^1( \TT^{\ast}_{\xi})\ar[r]&H^1 (\TET) \ar[r]^(0.35){\xi}& H^2(\OO_X (-K_X))
}
\end{equation}
where the coboundary map above is the cup-product with $\xi$, i.e.,
the map $H^1 (\TET) \stackrel{\xi}{\longrightarrow}  H^2(\OO_X (-K_X))$ above
is the quadratic form $\delta_X$ in \eqref{qf} restricted to the subspace
$\xi\cdot H^1 (\TET) \subset Sym^2(H^1 (\TET))$. From \eqref{xi-ort} we deduce
$$
H^1( \TT^{\ast}_{\xi}) \cong \{\xi\}^{\perp}:= \{\sigma \in H^1 (\TET) | \delta_X(\xi \cdot \sigma)=0 \}.
$$
From the first part of the assertion $2)$ it follows that $\xi$ lies in the image of the monomorphism in \eqref{xi-ort}. Hence $\delta_X(\xi^2)=0$. 
\end{pf}

\begin{rem} The fact that $\xi$ is isotropic with respect to the quadratic form
\eqref{qf} is imposed here by the presence of the extension sequence \eqref{ext-G} which is a subobject of \eqref{ext}. This is the meaning of the diagram \eqref{inc-seq}. The situation could be thought of as the object \eqref{ext} being
 destabilized by the subobject \eqref{ext-G} . In the later part of the proof the logic will reversed: the condition of isotropy of $\xi$ will be used to produce auxiliary, geometrically meaningful subobjects of \eqref{ext}. 
%
%
\end{rem}

Next we establish some geometric properties of global sections of $\GG$.

\begin{cl}\label{cl-geom-G}
Under the identifications
$$
\HKX \stackrel{\alpha}{ \cong} H^0 (\TT_{\xi})\stackrel{\tau}{ \cong} H^0 ({\cal G}),
$$
where a section $\phi \in \HKX$ goes to $\alpha (\phi) \in H^0 (\TT_{\xi})$ and $g(\phi)\stackrel{\rm{def}}{=}\tau(\alpha (\phi)) \in  H^0 ({\cal G})$, one has the following.

 Let $C_{\phi} =(\phi =0)$ be the divisor corresponding to a nonzero $\phi \in \HKX$. Then a section $g(\phi)$ gives rise to a section 
$\sigma (\phi) \in H^0 (\OO_{C_{\phi}} (D))$ and the zero-locus
$Z_{g(\phi)} =(g(\phi) =0) = (\sigma (\phi) =0)$ is a subscheme of $C_{\phi}$.
In particular, $ Z_{g(\phi)}$ is $0$-dimensional, for every $C_{\phi}$ reduced and irreducible, and its degree $deg( Z_{g(\phi)} ) =D.\KX$.
\end{cl} 
\begin{pf}
Let $g(\phi)$ be a nonzero section of ${\cal G}$ and view it as a monomorphism
$\OO_X \longrightarrow {\cal G}$. Putting it together with the extension sequence in (\ref{ext-G}) gives
the diagram
\begin{equation}\label{ext-G-sec}
\xymatrix@R=12pt@C=12pt{
& &\OO_X \ar[d] \ar[dr]^(0.45){\phi}& &\\
0\ar[r]& \OO_X (D) \ar[r]&{\cal G} \ar[r]&\OO_X (\KX) \ar[r]&0}
\end{equation}
where the slanted arrow is the multiplication by $\phi$. In particular, the
restriction of the above diagram to the divisor $C_{\phi} =(\phi =0)$ implies that the vertical arrow (which must be nonzero in view of $H^0 ({\cal G}(-K_X)) =0$, an immediate consequence of Lemma \ref{ext-incl}, 1)) factors through $\OO_{C_{\phi}} (D)$, thus giving a section,
call it
$\sigma({\phi})$, of $\OO_{C_{\phi}} (D)$. Furthermore, if $C_{\phi}$ is irreducible
the quotient of the vertical arrow in (\ref{ext-G-sec}) is torsion-free sheaf of rank $1$ and hence of the form ${\cal I}_{Z_{g(\phi)}} (K_X +D)$, where $Z_{g(\phi)} = (g(\phi) =0)$ is $0$-dimensional and ${\cal I}_{Z_{g(\phi)}}$ is its sheaf of ideals. In particular,
$deg (Z_{g(\phi)} ) = c_2 ({\cal G}) =K_X \cdot D.$ The section
$\sigma({\phi}) : \OO_{C_{\phi}} \longrightarrow \OO_{C_{\phi}} (D)$ vanishes
on the $0$-dimensional subscheme $D_{\phi} =(\sigma({\phi}) =0) \subset Z_{g(\phi)}$.
Since $degD_{\phi} =D\cdot K_X =deg (Z_{g(\phi)} )$, it follows that
 $D_{\phi} = Z_{g(\phi)}$.
\end{pf}

To analyze the situation further we will need the following general observation about
$\xi$.

\begin{lem}\label{lem-xi-log}
Let $C$ be a smooth irreducible curve on $X$ and let $\TET (-logC)$ be the sheaf of germs of holomorphic vector fields on $X$ tangent along $C$ . Then $\xi$ lies in the image of the map
\begin{equation}\label{log-xi}
H^1 (\TET (-logC)) \longrightarrow H^1 (\TET)
\end{equation}
induced by the natural inclusion
$\TET (-logC) \hookrightarrow \TET$, provided the restriction homomorphism
$\HKX \longrightarrow H^0 (\OO_C (\KX))$ is surjective.
\end{lem}
 \begin{pf}
By definition $\TET (-logC)$ is related to $\TET$ by the following exact sequence
\begin{equation}\label{log-seq}
\xymatrix@R=12pt@C=12pt{
0\ar[r]& \TET (-logC)\ar[r]& \TET \ar[r]& \OO_C (C) \ar[r]& 0.}
\end{equation}
So the assertion is equivalent to showing that $\xi$ goes to zero under the
homomorphism
$H^1 (\TET) \longrightarrow H^1 (\OO_C (C) )$ in the long exact sequence of the 
cohomology groups associated to (\ref{log-seq}). This can be seen by examining the multiplication by sections of $\OO(K_X)$. Namely, for every $\phi \in \HKX$ we have a commutative square
$$
\xymatrix@R=12pt@C=12pt{
H^1 ( \TET) \ar[r] \ar[d]^{\phi}& H^1 (\OO_C (C) )\ar[d]^{\overline{\phi}} \\
H^1 ( \OM) \ar[r]\ar[r]& H^1 (\OO_C (C+K_X)) }
$$
where $\overline{\phi}$ is the restriction of $\phi$ to $C$.
Since $\xi \cdot \phi =0$ it follows that the image $\overline{\xi}$ of $\xi$ in
$H^1 (\OO_C (C) )$ is annihilated by $\overline{\phi}$, for every $\phi \in \HKX$. Since $\HKX \longrightarrow  H^0 (\OO_C (\KX))$ is assumed to be onto, we deduce that the linear map
$$
\overline{\xi} :H^0 (\OO_C (K_X) ) \longrightarrow H^1 (\OO_C (C+K_X)) \cong \CC
$$
is identically zero. But this linear map is identified with $\overline{\xi}$ under the Serre
duality isomorphism
$H^0 (\OO_C (K_X) )^{\ast} \cong H^1 ( \OO_C (C) )$. Hence the assertion of the claim.
\end{pf}

The main point in ruling out the rank $2$ case consists of showing that $\xi$
comes from the first cohomology group of a sufficiently negative line bundle.
 Namely, for every nonzero $\psi \in H^0 (\OO_X (mK_X))$, for any $m \geq 1$, one has the following commutative diagram
\begin{equation}
\label{cohsquare}
\xymatrix@R=12pt@C=14pt{
H^1 \left(\OO_X (D-(m+1) \KX) \right)\ar[d] \ar[r]^(0.58){\psi}& H^1 (\OO_X (D-\KX))\ar[d] \\
H^1 (\TET (-m \KX)) \ar[r]^{\psi}& H^1 (\TET)
}
\end{equation}
 where the horizontal arrows are the multiplication by $\psi$ and the vertical arrows come from the monomorphisms in \eqref{TET-seq} and \eqref{TET-seq} tensored with $\OO_X (-mK_X)$ respectively. According to Lemma \ref{ext-incl}, the cohomology class $\xi \in  H^1 (\TET)$ comes from the class $\xi' \in H^1 (\OO_X (D-\KX))$ defining the extension sequence \eqref{ext-G}. We claim the following.
\begin{lem}
\label{l:xi-K}
The cohomology class $\xi'$ lies in the image of the homomorphism
$$
\xymatrix@R=12pt@C=12pt{
H^1 (\OO_X (D-(m+1)\KX)) \ar[r]^(0.58){\psi}& H^1 (\OO_X (D-\KX)),}
$$
for any nonzero $\psi \in H^0 (\OO_X (mK_X))$. In particular, the cohomology class $\xi$ comes from some cohomology class in $H^1 (\OO_X (D-(m+1)\KX))$ via the diagram \eqref{cohsquare}.
\end{lem}

\noindent
{\it Proof.} Let $C$ be a smooth curve in the linear system $\left|mK_X \right|$, for some $m\geq 1$. The assertion comes down to showing that the restriction to $C$ of the extension
sequence \eqref{ext-G} splits. This will be done by a careful examination of the relation between the extension classes $\xi'$ and $\xi$ encapsulated in the diagram
\eqref{inc-seq}.
 
From Lemma \ref{lem-xi-log} we know that $\xi$ comes from a cohomology class in
$H^1 (\TET (-log C))$. Choose such a class and call it $\eta$. Consistent with our approach we view it as an extension sequence 
{\small
$$
\xymatrix@R=12pt@C=12pt{ 
0\ar[r]& \TET (- logC)\ar[r]&{\cal E}_{\eta} \ar[r]& \OO_X\ar[r]& 0. }
$$
}
The fact that $\eta$ goes to $\xi$ as described in Lemma \ref{lem-xi-log} 
means that the extension sequences (twisted by $\OO_X (K_X)$) of those cohomology classes are related as follows.
{\small 
\begin{equation}\label{log-TET}
\xymatrix@R=12pt@C=12pt{
&0 \ar[d]&0 \ar[d]& & \\
0 \ar[r]&\TET (-log C) (K_X) \ar[r] \ar[d]& {\cal E}_{\eta} (K_X)\ar[r] \ar[d]&\OO_X (K_X) \ar[r] \ar@{=}[d]&0 \\
0 \ar[r]&\OM \ar[r] \ar[d]& \TT_{\xi} \ar[r] \ar[d]&\OO_X  (K_X)\ar[r]&0 \\
&\OO_C (C+K_X) \ar[d] \ar@{=}[r]& \OO_C (C+K_X)\ar[d]& &\\
& 0&0& &}
\end{equation}
}

\noindent
 We begin by the following observation.
\begin{cl}
\label{cl-Txi-C}  
The linear map 
\begin{equation}\label{Txi-C}
H^0 (\TT_{\xi}) \longrightarrow H^0 (\OO_C (C+K_X))
\end{equation}
 induced by the epimorphism
$\TT_{\xi} \longrightarrow \OO_C (C+K_X)$ in (\ref{log-TET})
 is nonzero.
\end{cl}
\begin{pf}
 Assume that the map in question is zero. Then all sections of
$\TT_{\xi}$ come from the global sections of ${\cal E}_{\eta} (K_X)$ and we have
the parametrization
$$
\alpha^{\prime} :\HKX \longrightarrow H^0 ({\cal E}_{\eta} (K_X))
$$
lifting the parametrization $\alpha$ in (\ref{alpha}).

Observe that for all $\phi_1, \phi_2 \in \HKX$ the global sections
\begin{equation}\label{b-pr}
\beta (\phi_1,\phi_2) = \phi_2 \alpha^{\prime} (\phi_1) - \phi_1 \alpha^{\prime} (\phi_2)
\end{equation}
of ${\cal E}_{\eta} (2K_X)$ go to zero under the homomorphism
$H^0 ({\cal E}_{\eta} (2K_X)) \longrightarrow H^0 (\OO_X (2\KX))$ induced by the 
epimorphism in the top exact sequence (tensored with $\OO_X (\KX)$) in (\ref{log-TET}). This implies that $\beta (\phi_1,\phi_2)$
 must come from the sections of 
$\TET (-log C) (2K_X) = \OM (log C)(\KX-C)=\OM (log C)(-(m-1)\KX)$. But it is well-known that 
$H^0 (\OM (log C))\cong H^0 (\OM ) $ which is zero by our assumption.
Thus we obtain
$$
\beta (\phi_1,\phi_2) = \phi_2 \alpha^{\prime} (\phi_1) - \phi_1 \alpha^{\prime} (\phi_2) =0,
$$
for all  $\phi_1, \phi_2 \in \HKX$. This implies that there is a section $\sigma \in H^0 ({\cal E}_{\eta} )$ such that sections $\alpha^{\prime} (\phi)$ have the form
$$
\alpha^{\prime} (\phi) = \phi \sigma,\,\, \forall \phi\in \HKX.
$$ 
 The section $\sigma$ delivers a monomorphism
$$
\OO_X (K_X) \longrightarrow {\cal E}_{\eta} (K_X)
$$
which gives the splitting of the top (and hence the middle) row in (\ref{log-TET}).
\end{pf}

We now bring in the subsheaf $\GG\hookrightarrow \TT_{\xi}$ in (\ref{inc-seq}) and consider the morphism
\begin{equation}\label{G-C}
\GG \longrightarrow \OO_C (C+K_X)=\OO_C ((m+1)K_X).
\end{equation}
induced by the epimorphism $\TT_{\xi} \longrightarrow \OO_C (C+K_X)$ in (\ref{log-TET}), where a smooth curve $C \in |mK_X|$ is chosen {\it not} to pass through any of the points of the $0$-dimensional subscheme $A$ in (\ref{inc-seq}).

Since all global sections of $\TT_{\xi}$ come from global sections in $\GG$,
we are assured, by Claim \ref{cl-Txi-C}, that the morphism in (\ref{G-C})
 is nonzero. The image of that morphism has the form
$ \OO_C ((m+1)K_X|_C -M)$ for some effective divisor $M$ on $C$. Thus we obtain
{\small
\begin{equation}\label{GG-Txi-C}
\xymatrix@R=12pt@C=12pt{
0\ar[r]&\GG \ar[r] \ar[d]&\TT_{\xi} \ar[r]\ar[d]& {\cal I}_{A} (K_X -D)\ar[r]\ar[d]&0\\
0\ar[r]&\OO_C ((m+1)K_X|_C -M) \ar[r]\ar[d]&\OO_C ((m+1)K_X) \ar[r]\ar[d]&\OO_M ((m+1)K_X)\ar[r]\ar[d]&0\\
&0&0&0}
\end{equation}
Our task will be to understand the subscheme $M$. For this we observe that all global sections of $\TT_{\xi}$ must vanish along the subsheaf
$\OO_X (D-K_X ) \hookrightarrow \TT^{\ast}_{\xi}$ (this is seen by dualizing the middle row in (\ref{GG-Txi-C}) and using $H^0 (\TT_{\xi} ) \cong H^0 (\GG))$.
On the other hand the restriction to $C$ of the {\it dual} of the left column in (\ref{inc-seq})
combined with the normal sequence of $C \subset X$ gives the following diagram
{\small
\begin{equation}\label{rest-to-C}
\xymatrix@R=12pt@C=12pt{
&&0 \ar[d]& & \\
&& \OO_C (D-K_X) \ar[d]\ar[rd]& \\
0 \ar[r]&\OO_C (-C -K_X)\ar[r] \ar[rd]& \TET\otimes \OO_C \ar[r] \ar[d]&\OO_C  (C)\ar[r]&0 \\
&&\OO_C ( -D) \ar[d] & &\\
& &0& &}
\end{equation}
}

\noindent
where the slanted arrows give rise to a nonzero global section, call it
$\tau_C $, of  $\OO_C (C+K_X -D)$ (the slanted arrow on the top in the above diagram is nonzero since otherwise it gives a nonzero section of $\OO_C (-C -D)=\OO_C (-mK_X -D)$ which, in view of the positivity of $K^2_X$ and Claim \ref{cl-geom-G}, is impossible). Let $\Delta_C = (\tau_C = 0)$ be the zero-divisor of $\tau_C $. It is important to observe that this is precisely the divisor where the curve $C$ is tangent to the foliation of $X$ defined by the subsheaf
$\OO_X (D-K_X) \hookrightarrow \TET$ (this inclusion is the dual of the epimorphism in the left column in (\ref{inc-seq})). This together with the previous
observation about global sections of $\TT_{\xi}$ vanishing along the subsheaf $\OO_X (D-K_X)  \hookrightarrow \TT^{\ast}_{\xi} $ imply that the image $\OO_C ((m+1)K_X |_C -M)$ of the morphism $\GG\longrightarrow \OO_C ((m+1)K_X )$ in (\ref{G-C}) factors through $\OO_C ((m+1)K_X |_C -\Delta_C) = \OO_C (D)$. Thus we obtain a nonzero morphism 
$$
\GG \longrightarrow \OO_C (D).
$$
Putting this together with the restriction to $C$ of the extension sequence (\ref{ext-G}) gives the following.
\begin{equation}\label{G-D}
\xymatrix@R=12pt@C=12pt{
0 \ar[r]&\OO_C (D) \ar[r] \ar[dr]& \GG \otimes \OO_C \ar[r] \ar[d] &\OO_C (K_X) \ar[r] &0\\
& &\OO_C (D)& &}
\end{equation}
The slanted arrow in the above diagram must be nonzero, since otherwise one has a nonzero morphism
$\OO_C (K_X) \longrightarrow \OO_C (D)$ implying the inequality $D.K_X \geq K^2_X$
contradicting the condition $D.K_X \leq \frac{1}{2} K^2_X$ of semistability of $\OM$ with respect to $K_X$, see \cite{Ts}. 

Once the slanted arrow in (\ref{G-D}) is nonzero it must be an isomorphism and it gives a splitting of the horizontal sequence in (\ref{G-D}).  This completes the proof of the lemma. $\Box$
 
\vspace{0.2cm}
The above result shifts our attention to the cohomology
 $H^1 (\OO_X (D-(m+1)K_X))$ which obviously vanishes for all $m$ sufficiently large. Hence $\xi=0$, thus ruling out the possibility
for ${\cal G}$ to be of rank $2$. This completes the proof of Lemma \ref{lem-glgen}.

\vspace{0.2cm}
The argument of the proof of  Lemma \ref{lem-glgen} applies without changes for the following more general situation.
For a nonzero $\xi \in H^1 (\Theta_X)$ we set
\begin{equation}\label{def:Wxi}
W_{\xi}:= ker (\HKX \stackrel{\xi}{\longrightarrow} H^1 (\Omega_X).
\end{equation} 
Then the extension sequence
$$
\xymatrix@R=12pt@C=12pt{
0\ar[r]& \Omega_X \ar[r]& \TT_{\xi} \ar[r]& \OO_X (\KX)\ar[r]&0.
}
$$
 as in \eqref{ext} gives rise to an identification
$$
\alpha_{\xi} :W_{\xi} \cong H^0 (\TT_{\xi}).
$$

\begin{lem}\label{lem-glgenWxi}

Assume that a nonzero Kodaira-Spencer class $\xi$ satisfies the following properties.

a) For some ample line bundle $\OO_X (H)$ and any divisor $C\in |mH|$, for $m>>0$, the class $\xi$ lies in the kernel of the obvious homomorphism
$$
 H^1 (\Theta_X) \longrightarrow H^1 (\OO_C (C)).
$$

b) $dim(W_{\xi}) \geq 3$ and the linear subsystem $|W_{\xi} | \subset |\KX|$ has at most $0$-dimensional base locus.

\vspace{0.2cm}
\noindent
Then the locally free sheaf $\TT_{\xi}$ is generically generated by the subspace $\alpha_{\xi} (W) \subset H^0 (\TT_{\xi})$, for any subspace $W \subset W_{\xi}$ of dimension at least $3$ with the linear subsystem $|W|\subset |\KX|$ having at most  $0$-dimensional base locus.
\end{lem}
\begin{pf}
Let $W$ be a subspace of $W_{\xi}$ of dimension at least $3$ such that the linear subsystem $|W|\subset |\KX|$ has at most  $0$-dimensional base locus. Assume that the subspace $\alpha_{\xi} (W)$ of global sections of $\TT_{\xi}$ fails to generate that sheaf everywhere. Then, as in the beginning of \S1, we define the subsheaf ${\cal G} \subset \TT_{\xi}$ as the saturation of
the image of the composition morphism
$$
W\otimes \OO_X \stackrel{\alpha_{\xi}}{\longrightarrow} H^0 (\TT_{\xi})\otimes \OO_X  \stackrel{ev}{\longrightarrow} \TT_{\xi}.
$$
As before ${\cal G}$ is locally free. The image of the morphism
$$
\xymatrix@R=12pt@C=12pt{
{\cal G} \ar[d] \ar[dr]&\\
 \TT_{\xi} \ar[r]& \OO_X (\KX)
}
$$
given by the slanted arrow has the form ${\cal I}_B (\KX)$, where $B$ is the base locus of the subsystem $|W|$ and ${\cal I}_B $ is its sheaf of ideals. This implies that ${\cal G}$ is of rank $2$ and fits into the diagram
\begin{equation}\label{inc-diag1}
\xymatrix@R=12pt@C=12pt{
&0\ar[d]&0\ar[d]& 0\ar[d]& \\ 
0\ar[r]& \OO_X (D) \ar[r] \ar[d]&{\cal G} \ar[r] \ar[d]&
{\cal I}_B (\KX) \ar[r] \ar[d]&0\\
0\ar[r]& \OM \ar[r]^{i} \ar[d]&\TT_{\xi} \ar[r]^(0.35){p} \ar[d]& \OO_X(\KX)\ar[r]
\ar[d]&0\\
0\ar[r]& {\cal I}_{A'} (K_X-D) \ar[r] \ar[d]& {\cal I}_{A''} (K_X-D)\ar[d]\ar[r]&\OO_B (K_X) \ar[r] \ar[d]&0\\ 
   & 0 & 0&0 &}
\end{equation}
which is an analogue of the diagram \eqref{inc-seq}. From the bottom row we deduce that the $0$-dimensional subschemes in the above diagram are subject to the relation 
$$
A'=B+A''.
$$

Next we take a smooth irreducible curve $C \in |mH|$, for some $m>>0$, and not passing through any of the points in $A'$. From the relation above it follows that $C$ is disjoint from $B$. Hence the restriction to $C$ of the top row in \eqref{inc-diag1} gives an exact sequence
\begin{equation}\label{G-Cext}
\xymatrix@R=12pt@C=12pt{
0\ar[r]& \OO_C (D) \ar[r] &{\cal G}\otimes \OO_C \ar[r] &
\OO_C (\KX) \ar[r] &0}
\end{equation}
and we claim that it splits. This is seen as follows.

We know that $\TT_{\xi}$ admits an epimorphism
$
\TT_{\xi} \longrightarrow \OO_C (C+\KX),
$
see \eqref{log-TET}. Thus restricting to $C$ we obtain an exact sequence
\begin{equation}\label{epi-C}
\xymatrix@R=12pt@C=12pt{
 0\ar[r]& {\cal E} \ar[r] &\TT_{\xi}\otimes \OO_C \ar[r]&\OO_C (C+\KX)
 \ar[r] &0
}
\end{equation}
where ${\cal E}$ is a rank $2$ bundle whose determinant $\det({\cal E})=\OO_C (\KX -C)=\OO_C (\KX -mH) $. Hence for $m$ sufficiently large its degree is negative. We now bring the subbundle ${\cal G}\otimes \OO_C \subset \TT_{\xi}\otimes \OO_C$. Combining this inclusion with the epimorphism in \eqref{epi-C} gives a morphism
\begin{equation}\label{Gto(C+K)}
{\cal G}\otimes \OO_C \longrightarrow \OO_C (C+\KX)
\end{equation}
which must be nonzero, since otherwise ${\cal G}\otimes \OO_C \cong {\cal E}\otimes \OO_C$ which is impossible because of the degree consideration.
  
 To understand the image of the above morphism we use the global sections
of ${\cal G}$ (they generically generate ${\cal G}\otimes \OO_C$) and the observation that their image in $\OO_C (C+\KX)$ must be
sections vanishing along the divisor $\Delta_C \in |\OO_C (C+\KX-D)$, for details see the argument just below the diagram \eqref{rest-to-C}. Hence we deduce that the morphism in \eqref{Gto(C+K)} factors through 
$\OO_C (C+\KX) \otimes \OO_C (-\Delta_C )=\OO_C (D)$. Therefore ${\cal G}\otimes \OO_C$ admits a nonzero morphism 
$$
{\cal G}\otimes \OO_C \longrightarrow \OO_C (D).
$$ 
This together with the exact sequence \eqref{G-Cext} gives the diagram 
$$
\xymatrix@R=12pt@C=12pt{
0\ar[r]& \OO_C (D) \ar[r] \ar[dr]&{\cal G}\otimes \OO_C \ar[r]\ar[d] &
\OO_C (\KX) \ar[r] &0\\
& &\OO_C (D)&& }
$$
with the slanted arrow being nonzero, see the argument just below \eqref{G-D}.
This implies the splitting of the horizontal sequence above.

The splitting also implies that $H^0 ({\cal G}\otimes \OO_C (-K_X)) \neq 0$, for
all $C\in |mH|$. But for $m>>0$ we clearly have an isomorphism
$H^0 ({\cal G}\otimes \OO_C (-K_X) \cong H^0 ({\cal G}(-K_X))$. This and the middle column in \eqref{inc-diag1} gives a nonzero global section of $\TT_{\xi} (-\KX)$ and hence the splitting of our extension sequence in the middle row of that diagram.  
\end{pf}

\section{Study of the degeneracy locus of $\TT_{\xi}$}

We know now that $\TT_{\xi}$ is generically generated by its global sections. This means that the evaluation morphism
$$
ev: H^0 (\TT_{\xi}) \otimes \OO_X \longrightarrow \TT_{\xi}
$$
is generically surjective. Let $D_{\xi}$ be the subscheme of $X$ where that morphism fails to be onto.
  The subscheme $D_{\xi}$ admits the rank stratification: for $k\in \{0,1,2\}$, we let $D^k_{\xi}$ to be the subscheme of $D_{\xi}$ whose closed points are defined as follows:
$$
D^k_{\xi}=\{ x\in X | rk (ev_x) \leq k \}.
$$
This gives the stratification
\begin{equation}\label{xi-strat}
\emptyset =D^0_{\xi} \subset D^1_{\xi}\subset D^2_{\xi} = D_{\xi}.
\end{equation}
The main purpose of this section is to show that if $D_{\xi}$ has the nonzero divisorial part, then $\TT_{\xi}$ can be modified along that divisor to obtain a new locally free sheaf with the same space of global sections and whose degeneracy locus is at most $0$-dimensional. Such a modification is achieved in two stages:

\vspace{0.2cm}
1) a modification which eliminates the divisorial part (if nonzero) of $D^1_{\xi}$,

2) a modification producing a locally free sheaf with at most $0$-dimensional degeneracy locus.

\vspace{0.2cm}
It should be stressed that the modified sheaves are part of extension sequences
which are subsequences of the initial one in \eqref{ext}. So what we are producing here is a particular filtration of the object \eqref{ext} by its subobjects in the category of exact complexes of coherent sheaves of $X$. That filtration could be viewed as a categorical counterpart of the stratification \eqref{xi-strat}.
\vspace{0.2cm}
\subsection{The divisorial part of $D^1_{\xi}$} 
We will be concerned here with the stratum $D^1_{\xi}$
 and, in particular, with its divisorial component which will be denoted $D^{1,1}_{\xi}$. We begin by establishing a relation between the cohomology class $\xi$ and $D^{1,1}_{\xi}$.
\begin{pro}\label{pro:eta}
1) The cohomology class $\xi$ lies in the image of the homomorphism
$$
H^1 (\Theta_X (-D^{1,1}_{\xi} )) \longrightarrow H^1 (\Theta_X )
$$
induced by the multiplication morphism
$\Theta_X (-D^{1,1}_{\xi} ) \longrightarrow \Theta_X  $.

\vspace{0.2cm}
2) There is a cohomology class $\eta \in H^1 (\Theta_X (-D^{1,1}_{\xi} ))$  mapping to $\xi$ under the homomorphism in 1) and lying in the kernel of the
homomorphism
$$
H^1 (\Theta_X (-D^{1,1}_{\xi} )) \longrightarrow \HKX^{\ast} \otimes H^1 (\Theta_X (K_X - D^{1,1}_{\xi} )) \cong H^1 (\Omega_X ( - D^{1,1}_{\xi} )).
$$ 

3) $H^0 (\OO_X (2K_X -2D^{1,1}_{\xi})) \neq 0$.
\end{pro}

\begin{pf}
To simplify the notation we set $D^{1,1}_{\xi} =D$ and consider the restriction of the extension sequence \eqref{ext} to $D$.
$$
\xymatrix@R=12pt@C=12pt{
0 \ar[r]&\Omega_X \otimes \OO_D \ar[r] & \TT_{\xi} \otimes \OO_D \ar[r]  &\OO_D (K_X) \ar[r] &0.}
$$
The fact that the evaluation morphism has rank $1$ on $D$ implies that the image of the evaluation morphism is an invertible subsheaf of $\TT_{\xi} \otimes \OO_D$. Call that subsheaf ${\cal L}$. By definition we have
\begin{equation}\label{extonD}
\xymatrix@R=12pt@C=12pt{
&&0\ar[d]& &\\
& &{\cal L} \ar[d] \ar[dr]& &\\
0 \ar[r]&\Omega_X \otimes \OO_D \ar[r] & \TT_{\xi} \otimes \OO_D \ar[r]  &\OO_D (K_X) \ar[r] &0.}
\end{equation}
Furthermore, the slanted arrow must be surjective everywhere and hence an isomorphism since both sheaves are of rank $1$. The surjectivity of the slanted arrow
comes from the fact that the homomorphism
$$
H^0 ({\cal L} ) \longrightarrow H^0(\OO_D (K_X))
$$
 induced on the global sections has the image which generates the sheaf
$\OO_D (K_X)$.

Once the slanted arrow in \eqref{extonD} is an isomorphism, the short exact sequence in that diagram splits. The short exact sequence in \eqref{extonD} corresponds to the image
of $\xi$ in $H^1(\Theta_X \otimes \OO_D)$ under the homomorphism
\begin{equation}\label{extonD-H1}
H^1(\Theta_X ) \longrightarrow H^1(\Theta_X \otimes \OO_D)
\end{equation}
induced by the restriction morphism of sheaves
$\Theta_X \longrightarrow \Theta_X \otimes \OO_D$. Hence $\xi$ lies in the kernel of \eqref{extonD-H1} or, equivalently, in the image
$$ 
H^1 (\Theta_X (-D)) \longrightarrow H^1(\Theta_X ).
$$ 

We now turn to the part 2) of the proposition. From the splitting of the exact sequence in \eqref{extonD} it follows that there is a surjective morphism
$$
\TT_{\xi} \longrightarrow \Omega_X \otimes \OO_D
$$
whose kernel, call it $\TT$, is a subsheaf of $\TT_{\xi} $ with 
$H^0 (\TT) \cong H^0(\TT_{\xi})$. This implies that the morphism
$\TT \longrightarrow \OO_X (\KX)$, the composition of the inclusion 
$\TT \subset \TT_{\xi}$ with the epimorphism in \eqref{ext}, is still surjective
and it gives rise to the following commutative diagram:
\begin{equation}\label{ext-extD}
\xymatrix@R=12pt@C=12pt{
&0\ar[d]&0\ar[d]&&\\
0 \ar[r]&{\cal K} \ar[r] \ar[d]&\TT \ar[r] \ar[d]& \OO_X (\KX) \ar[r] \ar@{=}[d]&0\\ 
0 \ar[r]&\Omega_X \ar[r]\ar[d]&\TT_{\xi} \ar[r]\ar[d]& \OO_X (\KX) \ar[r]&0 \\
&\Omega_X \otimes \OO_D \ar[d]\ar@{=}[r]&\Omega_X \otimes \OO_D \ar[d]&&\\
&0&0& &
}
\end{equation}
where the sheaf ${\cal K}$ is the kernel of the epimorphism in the top row. From the column on the left that sheaf is identified with $\Omega_X (-D)$ and the top horizontal sequence becomes
\begin{equation}\label{ext-D}
\xymatrix@R=12pt@C=12pt{
0 \ar[r]&\Omega_X (-D) \ar[r] &\TT \ar[r] & \OO_X (\KX) \ar[r] &0.
}
\end{equation}
This is an extension sequence corresponding to a class in
$Ext^1 (\OO_X (K_X), \Omega_X (-D)) \cong H^1 (\Theta_X (-D))$. Furthermore, the
morphism between the rows in \eqref{ext-extD} means that this cohomology class in $ H^1 (\Theta_X (-D))$ goes to $\xi \in H^1 (\Theta_X)$ under the homomorphism
$$
H^1 (\Theta_X (-D)) \longrightarrow H^1 (\Theta_X)
$$
 defined by the multiplication by a section defining the divisor $D$.

Let $\eta$ be the cohomology class in $H^1 (\Theta_X (-D))$ corresponding to the extension in \eqref{ext-D}. By construction we have
$$
H^0 (\TT) \cong H^0(\TT_{\xi}) \cong H^0 (\OO_X (\KX)).
$$
 This implies that the coboundary map
$H^0 (\OO_X (\KX)) \longrightarrow H^1 (\Omega_X (-D))$ in the cohomology sequence of \eqref{ext-D} is identically zero. Since that coboundary map is given by the cup-product with the class $\eta$ we deduce that $\eta$ lies in the kernel of
the homomorphism
$$
H^1 (\Theta_X (-D)) \longrightarrow \HKX^{\ast} \otimes H^1 (\Omega_X (-D))
$$
as asserted in 2) of the proposition.

\vspace{0.2cm}
The part 3) of the proposition uses the fact that $\TT_{\xi}$ is generically generated by global sections, Lemma \ref{lem-glgen}. This amounts to saying that
for a general triple $\phi,\phi', \phi'' \in \HKX$ the exterior product 
$\alpha(\phi)\wedge \alpha(\phi')\wedge \alpha(\phi'') \in \bigwedge^3 H^0 (\TT_{\xi})$ goes to a nonzero global section of $\OO_X (2K_X)$ under the homomorphism
$$
\mbox{$\bigwedge^3 H^0 (\TT_{\xi})$} \longrightarrow H^0 (\det(\TT_{\xi}))=H^0 (\OO_X(2K_X)).
$$ 
Let $\alpha^{(3)}(\phi,\phi', \phi')$ be the image of $\alpha(\phi)\wedge \alpha(\phi')\wedge \alpha(\phi'')$ under the above map. The fact that the evaluation morphism drops its rank by $2$ along $D$ means that all global sections lying in the image of the homomorphism above vanish along $D$ with multiplicity $2$, i.e. the section $\alpha^{(3)}(\phi,\phi', \phi')$ has the form
$$
\alpha^{(3)}(\phi,\phi', \phi')=t^2_D \,\alpha'(\phi,\phi', \phi'),
$$
where $ \alpha'(\phi,\phi', \phi')$ is a nonzero section of $\OO_X (2K_X-2D)$ and $t_D$ is a global section of $\OO_X (D)$ corresponding to $D$. 
\end{pf}

The set of effective divisors of $X$ subject to Proposition \ref{pro:eta} is partially ordered by the inclusion and the part 3) of the proposition implies that the set admits maximal elements. We choose one and denote it by $E$.
By definition the cohomology group $H^1 (\Theta_X (-E ))$ has a class, call it 
$\eta$, which maps to $\xi$ under the homomorphism

$$
H^1 (\Theta_X (-E )) \longrightarrow  H^1 (\Theta_X ).
$$

Using the identification

$$
H^1 (\Theta_X (-E )) =Ext^1 (\OO_X (K_X), \Omega_X ( - E ))
$$

we interpret that class as the corresponding extension
\begin{equation}\label{ext-eta}
\xymatrix@R=12pt@C=12pt{
0\ar[r]& \Omega_X ( - E )  \ar[r]&{\cal T}_{\eta} \ar[r]&\OO_X (\KX) \ar[r]&0}
\end{equation}
The property 1) in Proposition \ref{pro:eta} implies that this extension sequence is related to the one defined by $\xi$, see \eqref{ext}, by the following diagram

\begin{equation}\label{diag:ext-eta-xi}
\xymatrix@R=12pt@C=12pt{
&0\ar[d]& 0\ar[d]& & &\\
0\ar[r]& \Omega_X ( - E)  \ar[r] \ar[d]&{\cal T}_{\eta} \ar[r] \ar[d]&\OO_X (\KX) \ar[r] \ar@{=}[d]&0\\
0\ar[r]& \Omega_X  \ar[r] \ar[d]&{\cal T}_{\xi} \ar[r]\ar[d]&\OO_X (\KX) \ar[r] &0\\
& {\Omega_X \otimes \OO_ {E } } \ar@{=}[r] \ar[d]&{\Omega_X \otimes \OO_ {E}} \ar[d] & &\\
&0& 0& & &}
\end{equation}
The meaning of this diagram is that the exact sequence \eqref{ext-eta} is a subobject of the sequence \eqref{ext} in the category of exact complexes of coherent sheaves on $X$.
Furthermore, by the property 2) of Proposition \ref{pro:eta} we have the isomorphisms
\begin{equation}\label{iso-glsec}
H^0 (\TT_{\eta}) \cong \HKX \cong H^0 (\TT_{\xi}).
\end{equation}
This immediately implies
\begin{pro}\label{pro:gglgTeta}
The sheaf $\TT_{\eta}$ is generically generated by its global sections.
\end{pro}
\begin{pf}
From the middle column of the diagram \eqref{diag:ext-eta-xi} it follows that
the sheaves $\TT_{\eta}$ and ${\cal T}_{\xi}$ are isomorphic outside of the divisor $E$. This together with the isomorphism $H^0 (\TT_{\eta})\cong H^0 (\TT_{\xi})$ provided by  \eqref{iso-glsec} and the generic global generation of  ${\cal T}_{\xi}$ implies the assertion.
\end{pf}

As for the sheaf $\TT_{\xi}$, we define the degeneracy locus $D_{\eta}$ of $\TT_{\eta}$ as the support of the cokernel of the evaluation morphism
$$
H^0 (\TT_{\eta}) \otimes \OO_X \longrightarrow \TT_{\eta}.
$$
Its rank stratification
\begin{equation}\label{Deta}
\emptyset=D^0_{\eta} \subset D^1_{\eta}  \subset D^2_{\eta}=D_{\eta}
\end{equation}
becomes simpler due to the following.
\begin{lem}\label{D1eta-dim}
The stratum $D^1_{\eta} $ is at most $0$-dimensional.
\end{lem}
\begin{pf}
Let $F$ be the $1$-dimensional part of $D^1_{\eta}$. Then arguing as in the proof of Proposition \ref{pro:eta}, we deduce that the cohomology class $\eta $

\vspace{0.2cm}
--  comes from a cohomology class $H^1 (\Theta_X (-E-F))$ under the 
homomorphism $H^1 (\Theta_X (-E-F))\longrightarrow H^1 (\Theta_X (-E))$ induced by the the multiplication with a global section of $\OO_X (F)$ corresponding to $F$,

\vspace{0.2cm}
--  the above cohomology class can be chosen to lie in the kernel of the cup-product
$$
H^1 (\Theta_X (-E-F))\longrightarrow \HKX^{\ast}\otimes H^1 (\Theta_X (K-E-F)) = H^1 (\Omega_X (-E-F)).
$$
But this means that the cohomology class $\xi$ comes from a cohomology class 
$H^1 (\Theta_X (-E-F))$ lying in the kernel of the cup-product above. The condition of maximality imposed on $E$ forces $F$ to be $0$. 
\end{pf}
The first isomorphism in \eqref{iso-glsec} will be recorded as 
\begin{equation}\label{alphaeta}
\alpha_{\eta} : \HKX \longrightarrow H^0 (\TT_{\eta}),
\end{equation}
a parametrization of global sections of $\TT_{\eta}$ by the space $\HKX$.
One of the advantages of working with the sheaf $\TT_{\eta}$ is the following property.
\begin{pro}\label{pro:glsecTeta}
A general global section of $\TT_{\eta}$ has no zeros.
\end{pro}
\begin{pf}
Consider the incidence correspondence
$$
{\cal Z}= \{(x, [\phi]) \in X\times \PP(\HKX) | \alpha_{\eta} (\phi) (x) =0 \},
$$
where $\alpha_{\eta}(\phi) (x)$ stands for the value of the global section $\alpha_{\eta}(\phi)$ at a point $x\in X$. The correspondence comes with two projections
$$
\xymatrix@R=12pt@C=12pt{
&{\cal Z} \ar[dr]^{p_1} \ar[dl]_{p_2}&\\
X& &\PP(\HKX)
}
$$
and the assertion of the proposition comes down to showing that $p_1 ({\cal Z} )$ is a proper subscheme of $\PP(\HKX)$. This will be done by counting the dimension of ${\cal Z}$.

Over the open surface $X'=X\setminus D_{\eta}$ the scheme ${\cal Z}$ is
$\PP^{p_g-4}$- bundle. Hence ${\cal Z}'=p^{-1}_2 (X')$ has dimension $(p_g-2)$.

Over the stratum $D'_{\eta}=D_{\eta} \setminus D^1_{\eta}$, the projection $p_2$ is
a $\PP^{p_g-3}$- bundle. Since $ dim (D'_{\eta}) \leq 1$, we have 
$$
dim(p^{-1}_2 (D'_{\eta}))\leq p_g-2.
$$

Over the stratum $D^1_{\eta}$ the projection $p_2$ is a $\PP^{p_g-2}$-bundle. By
Lemma \ref{D1eta-dim} that stratum is at most $0$-dimensional. Hence
$$
 dim(p^{-1}_2 (D^1_{\eta}))\leq p_g-2.
$$
The above considerations imply that all irreducible components of ${\cal Z}$ have dimension at most $p_g-2$. Hence $p_1 ({\cal Z})$ is a subscheme of codimension
at least $1$.
\end{pf}

We know now that the locus of $\PP(\HKX)$ parametrizing the projectivized global sections of $\TT_{\eta}$ with zeros is a proper subscheme of $\PP(\HKX)$.
We will also need the following.
\begin{lem}\label{lem:1dimzeroloc}
The subscheme of $\PP(\HKX)$ parametrizing the projectivized global sections of $\TT_{\eta}$ with $1$-dimensional zero locus has the codimension at least $2$.
\end{lem}
\begin{pf}
Let $\Sigma$ be a reduced irreducible subscheme of  $\PP(\HKX)$ whose closed
points parametrize the projectivized global sections of  of $\TT_{\eta}$ having  $1$-dimensional zero locus and assume it to be of codimension $1$ in $\PP(\HKX)$, i.e. $dim(\Sigma) =p_g-2$. Consider the incidence
$$
I_{\Sigma} =\{ (x,[\phi]) \in X\times \Sigma\, |\, \alpha(\phi) (x)=0 \}.
$$
This is a scheme of dimension $p_g-1$. Choose a reduced irreducible component
$\widetilde{\Sigma}$ of $I_{\Sigma}$ having dimension $(p_g-1)$ and consider
the morphism
$$
p_1:\widetilde{\Sigma} \longrightarrow X
 $$
induced by the projection of $I_{\Sigma} $ onto the first factor. Observe that $p_1$ is surjective. Indeed, otherwise the image is a reduced irreducible curve $\Gamma \subset X$ contained in the zero-locus of all sections $\alpha(\phi)$ with
$[\phi] \in \Sigma$. In particular, $\Gamma$ is in the zero-locus of all $[\phi] \in \Sigma$. But then $\Gamma$ is in the zero-locus of all 
$[\phi]$ in the linear span of $ \Sigma$ in $\PP(\HKX)$.
That linear span is a linear subspace of $\PP(\HKX)$ of codimension at least $1$ and, i.e. $h^0 (\OO_X (\KX-\Gamma) \geq p_g-1$. But this means that $\Gamma $ is either in the base locus of $|\KX|$ or is mapped onto a point under the canonical map. Either is in contradiction with the condition $(i)$ in \eqref{cond}.

Once the morphism $p_1$ is surjectve, the dimension of its fibre
$p^{-1}_1 (x)$ for a general $x\in X$ has dimension $(p_g -3)$. Under the projection $p_2:\widetilde{\Sigma} \longrightarrow \Sigma$ that fibre is identified with the sections $\alpha(\phi)$ vanishing at $x$ and $[\phi] \in \Sigma$. But the generic global generation of $\TT_{\eta}$, see Proposition \ref{pro:gglgTeta}, tells us that the projectivization of the subspace of global sections of $\TT_{\eta}$ vanishing at a general point on $X$ has dimension $(p_g-4)$. An obvious contradiction.  
\end{pf}

We now return to the parametrization \eqref{alphaeta}
 of global sections of $\TT_{\eta}$ by the space $\HKX$. It will be used to construct a family of sheaves of rank $2$ on $X$. Namely, for every nonzero $\phi \in \HKX $ we view the section $\alpha_{\eta} (\phi)$
as a (nonzero) morphism
$$
 \alpha_{\eta} (\phi): \OO_X \longrightarrow \TT_{\eta}
$$
and we set
\begin{equation}\label{Fphi}
{\cal F}_{[\phi]} :=coker( \alpha_{\eta} (\phi)),
\end{equation}
where $[\phi]$ is the point of the projective space $\PP(\HKX)$ underlying a vector $\phi \in \HKX$. 
\begin{pro}\label{pro:Fphi}
1) For every  $[\phi]\in \PP(\HKX)$ the sheaf ${\cal F}_{[\phi]}$ has rank $2$ and the Chern invariants
$$
c_1({\cal F}_{[\phi]})=2K_X -2E, \hspace{0.2cm} c_2 ({\cal F}_{[\phi]})=c_2 +K^2_X-
3K_X \cdot E +E^2.
$$

\vspace{0.2cm}
2) $H^0({\cal F}_{[\phi]}) \cong \HKX /\CC\phi$.

\vspace{0.2cm}
3) The sheaf ${\cal F}_{[\phi]}$ is locally free if and only if the section $\alpha_{\eta} (\phi)$ has no zeros. A general section of such ${\cal F}_{[\phi]}$ has at most $0$-dimensional scheme of zeros.
\end{pro}
\begin{pf}
By definition of $\FF_{[\phi]}$ we have the exact sequence
\begin{equation}\label{Fphi-ext}
\xymatrix@R=12pt@C=12pt{
0\ar[r]&\OO_X \ar[r]^{\alpha_{\eta} (\phi)}&\TT_{\eta} \ar[r]&{\cal F}_{[\phi]}\ar[r]&0
}
\end{equation}
from which it follows that ${\cal F}_{[\phi]}$ has rank $2$ and the Chern invariants equal the ones of $\TT_{\eta}$. The Chern invariants of the latter are easily computed from the extension sequence \eqref{ext-eta}.

The assertion 2) follows from \eqref{Fphi-ext} since that sequence gives
rise to the exact sequence on the level of global sections
$$
\xymatrix@R=12pt@C=12pt{
0\ar[r]&H^0 (\OO_X) \ar[r]^{\alpha_{\eta} (\phi)}&H^0(\TT_{\eta}) \ar[r]&H^0({\cal F}_{[\phi]})\ar[r]&0.
}
$$
This together with the parametrization \eqref{alphaeta} implies the identifications
$$
\HKX/\CC \phi \cong H^0(\TT_{\eta})/\CC\alpha_{\eta} (\phi)\cong H^0 ({\cal F}_{[\phi]}).
$$ 

The first assertion of part 3) is again  immediate from \eqref{Fphi-ext}, since by definition $\FF_{[\phi]}$ is locally free if and only if the monomorphism in that sequence is a monomorphism of vector bundles and this occurs precisely when 
$\alpha_{\eta} (\phi)$ has no zeros.  

For the second assertion of 3), we take $\phi' \in \HKX$ linearly independent of $\phi$ and consider the diagram
\begin{equation}\label{diag:phi-phi'}
\xymatrix@R=12pt@C=12pt{
& & \OO_X \ar[d]^{\alpha_{\eta} (\phi')}&&\\
0\ar[r]&\OO_X \ar[r]^{\alpha_{\eta} (\phi)}&\TT_{\eta} \ar[r]&{\cal F}_{[\phi]}\ar[r]&0
}
\end{equation}
The composition of the vertical arrow and the epimorphism of the horizontal sequence gives a morphism $ \OO_X \longrightarrow {\cal F}_{[\phi]}$ corresponding to the global section $f_{\phi'}$ of ${\cal F}_{[\phi]}$ given by the equivalence class of $\phi'$ in
$\HKX /\CC\phi$ under the isomorphism established in 2) of the proposition.
Observe that $f_{\phi'}$ has $1$-dimensional zero-locus iff there is a reduced irreducible curve $C \subset X$ such that the vertical arrow in \eqref{diag:phi-phi'} restricted to $C$ factors through
$$
\alpha_{\eta} (\phi) |_C : \OO_C \longrightarrow \TT_{\eta} \otimes \OO_C.
$$
This means that there is a constant $\lambda \in \CC$ such that
$$
\alpha_{\eta} (\phi') |_C =\lambda \alpha_{\eta} (\phi) |_C.
$$
Equivalently, the above relation means that the section
$\alpha_{\eta} (\phi'-\lambda \phi)$ vanishes on $C$. By Lemma \ref{lem:1dimzeroloc} a general line $l$ in $\PP(\HKX)$ passing through $[\phi]$ does not meet
the subscheme of $\PP(\HKX)$ parametrizing the global sections of  $\TT_{\eta}$ having $1$-dimensional zero-locus. Hence for a general choice of $\phi' \in \HKX$
the global sections
$\alpha_{\eta} (\psi)$ have at most $0$-dimensional scheme of zeros for all $[\psi]$ in the pencil spanned by $[\phi]$ and $[\phi']$.
\end{pf}

The parametrization $\alpha_{\eta}: \HKX \longrightarrow H^0 (\TT_{\eta})$ also gives rise to `higher products' 
\begin{equation}\label{alpha(k)}
\mbox{$\alpha^{(k)}_{\eta}: \bigwedge^k \HKX \longrightarrow H^0 (\bigwedge^k \TT_{\eta})$}
\end{equation}
for $k=2,3$, where $\alpha^{(k)}_{\eta}(\phi_1,\ldots,\phi_k)$ is defined to be the image of the exterior product 
$\displaystyle{\bigwedge^k_{i=1} \alpha_{\eta}(\phi_i)}$ under the natural homomorphism
$$
\mbox{$ \bigwedge^k H^0 (\TT_{\eta}) \longrightarrow H^0 (\bigwedge^k \TT_{\eta}).$}
$$
Notice that the map
$$
\mbox{$\alpha^{(3)}_{\eta}: \bigwedge^3 \HKX \longrightarrow H^0 (\bigwedge^3 \TT_{\eta}) =H^0 (\OO_X (2\KX-2E))$}
$$
is related to the analogous product $\alpha^{(3)}: \bigwedge^3 \HKX \longrightarrow H^0 (\bigwedge^3 \TT_{\xi}) =H^0 (\OO_X (2\KX))$ by the formula
\begin{equation}\label{alpha-alphaeta}
\alpha^{(3)} =\tau^2_E \alpha^{(3)}_{\eta},
\end{equation}
where $\tau_E$ is a global section of $\OO_X (E)$ corresponding to the divisor $E$.
 The following property follows easily from the previous results. 
\begin{cor}\label{cor:irred}
 1) For every $\phi \in \HKX$ with the divisor $C_{\phi} =(\phi=0)$ reduced and irreducible, the section $\alpha^{(3)}_{\eta} (\phi,\phi', \phi'') \neq 0$, for any
linearly independent triple $\phi,\phi', \phi'' \in \HKX$.

2) For any subsheaf ${\cal T} \subset \TT_{\eta}$ of rank $1$, one has
$h^0 ({\cal T}) \leq 1$. In particular, $\alpha^{(2)}_{\eta} (\phi,\phi')=0$ if and only if $\phi,\phi' \in \HKX$ are linearly dependent.
\end{cor} 
\begin{pf}
 Take $W:= \CC\{\phi, \phi', \phi''\}$ a $3$-dimensional subspace of $\HKX$. By the assumption on $\phi$ the base locus of the linear subsystem $|W|\subset |\KX|$
is at most $0$-dimensional. Applying Lemma \ref{lem-glgenWxi} to $W$, we deduce
that the subspace $\alpha(W) \subset H^0(\TT_{\xi})$ generically generates $\TT_{\xi}$ or, equivalently, $\alpha^{(3)} (\phi,\phi', \phi'') \neq 0$. From the relation
$$
\alpha^{(3)} (\phi,\phi', \phi'')= \tau^2_E \alpha^{(3)}_{\eta} ( \phi, \phi', \phi'')
$$
in \eqref{alpha-alphaeta} it follows that $\alpha^{(3)}_{\eta} ( \phi, \phi', \phi'')\neq 0$.

For the second assertion, assume that there is  a subsheaf ${\cal T} \subset \TT_{\eta}$ of rank $1$ with $h^0 ({\cal T}) \geq 2$. Then a pair $t', t''$ of linearly independent sections of ${\cal T}$ corresponds to two global sections $\alpha_{\eta} ( \phi'), \alpha_{\eta} (\phi'')$ of $\TT_{\eta}$ which are proportional, i.e.
$ \alpha^{(2)}_{\eta} ( \phi', \phi'')=0$. But then
 $\alpha^{(3)}_{\eta} ( \phi, \phi', \phi'') =0$, for all $\phi \in \HKX$, and this is in contradiction with the first assertion proved above.
\end{pf}

We will be concerned now with the image $Im(\alpha^{(3)}_{\eta})$ of $\alpha^{(3)}_{\eta}$. In particular, we are interested in the geometry of the linear subsystem $|Im(\alpha^{(3)}_{\eta})|$ of $|2\KX-2E|$. This linear subsystem may have the nonzero fixed part - the divisorial part of the degeneracy locus $D_{\eta}$ in \eqref{Deta}. In the next subsection we show that it is always possible to modify the extension
\eqref{ext-eta} so that the resulting linear system has no fixed part. 

\subsection{The divisorial part of $D_{\eta}$ and a modification of \eqref{ext-eta}}

 Let
$\phi, \phi'\in \HKX$ be a linearly independent pair and consider the linear map
$$
\alpha^{(3)}_{\eta} ( \phi, \phi',\bullet): \HKX \longrightarrow H^0 (\OO_X (2\KX-2E)).
$$
From Corollary \ref{cor:irred} this is nonzero. Its image $Im(\alpha^{(3)}_{\eta} ( \phi, \phi',\bullet))$ gives rise to a nonzero subspace of $H^0 (\OO_X (2\KX-2E))$.
Assume that the corresponding linear
subsystem $|Im(\alpha^{(3)}_{\eta} ( \phi, \phi',\bullet))|$ of $|2\KX-2E|$ has the nonzero fixed part (otherwise there is nothing to modify), call it $E_{\phi, \phi'}$. As $\CC\{\phi, \phi'\}$ vary in some Zariski dense open subset of the Grassmannian of lines in $\PP(\HKX)$, these divisors are rationally equivalent and we denote by $E_1$ the corresponding rational equivalence class in the Picard group $Pic(X)$ of $X$.
Thus $Im(\alpha^{(3)}_{\eta} ( \phi, \phi',\bullet))$ lies in the image of the multiplication map
\begin{equation}\label{fixed-phiphi'}
e_{\phi,\phi'}:H^0 (\OO_X (2\KX-2E-E_1))=  H^0 (\OO_X (2\KX-2E-E_{\phi, \phi'})) \longrightarrow H^0 (\OO_X (2\KX-2E)),
\end{equation}
where $e_{\phi,\phi'}$ is a global section of $\OO_X (E_1)$ corresponding to $E_{\phi, \phi'}$.

We will now interpret $Im(\alpha^{(3)}_{\eta} ( \phi, \phi',\bullet))$ from the point of view of the family of sheaves $\{\FF_{[\psi]}\}$. For this   
we choose two linearly independent $\phi, \phi' \in \HKX$ so that the sections $\alpha_{\eta} ( \phi)$ and $\alpha_{\eta} ( \phi')$ have no zeros. Those sections define the locally free sheaves
$\FF_{[\phi]}$ and $\FF_{[\phi']}$ respectively. Furthermore, we know that those sheaves are both subbundles of $\bigwedge^2\TT_{\eta}$. This is recorded in the following diagram.
\begin{equation}\label{diag:FphiFphi'}
\xymatrix@R=12pt@C=12pt{
&&0\ar[d]& &\\
&&\FF_{[\phi']}\ar[d]\ar[dr]& &\\
0\ar[r]& \FF_{[\phi]}\ar[r]\ar[dr]& \bigwedge^2\TT_{\eta}\ar[r]^(0.4){\wedge \alpha_{\eta}(\phi)} \ar[d]^{\wedge \alpha_{\eta}(\phi')} &\OO_X (2\KX-2E)\ar[r]&0
\\
&&\OO_X (2\KX-2E)\ar[d]&&\\
&&0&&}
\end{equation}
We also know that $H^0 (\FF_{[\phi']})$ (resp. $H^0 (\FF_{[\phi]})$ ) is isomorphic to the image $Im (\alpha^{(2)}_{\eta} (\bullet, \phi'))$ (resp. $Im (\alpha^{(2)}_{\eta} (\phi,\bullet,))$ of the linear map
$\alpha^{(2)}_{\eta} (\bullet, \phi'): \HKX \longrightarrow  H^0(\bigwedge^2\TT_{\eta})$
(resp., $ \alpha^{(2)}_{\eta} (\phi,\bullet,):\HKX \longrightarrow  H^0(\bigwedge^2\TT_{\eta})$). Hence the slanted arrow in the lower left (resp. upper right) corner of the above diagram
gives rise to the linear map
$$
H^0 (\FF_{[\phi]}) \longrightarrow H^0 (\OO_X (2\KX-2E)) \hspace{0.2cm} (\mbox{resp. $H^0 (\FF_{[\phi']}) \longrightarrow H^0 (\OO_X (2\KX-2E))$})
$$
whose image is precisely $Im(\alpha^{(3)}_{\eta} ( \phi, \phi',\bullet))$. From what was said above that map factors through $H^0 (\OO_X (2\KX-2E -E_1))$:
\begin{equation}\label{Fphiglobsect}
\xymatrix@R=12pt@C=12pt{
&H^0 (\OO_X (2\KX-2E -E_1))\ar[d]^{e_{\phi,\phi'}}\\
H^0 (\FF_{[\phi]}) \ar[r]\ar[ur]& H^0 (\OO_X (2\KX-2E))}
\end{equation}
where the vertical arrow is the map in \eqref{fixed-phiphi'}. 

Consider the composition
\begin{equation}\label{morph:cphiphi'}
c_{\phi,\phi'} :\FF_{[\phi]} \longrightarrow \OO_{E_{\phi,\phi'}} (2\KX-2E)
\end{equation}
of the slanted arrow in the lower left corner in the diagram \eqref{diag:FphiFphi'} with the restriction morphism 
$\OO_X (2\KX-2E) \longrightarrow \OO_{E_{\phi,\phi'}} (2\KX-2E)$.
\begin{lem}\label{lem:cokercphiphi'}
Let $\phi$ and $\phi'$ in $\HKX$ be two linear independent sections chosen so that $\alpha_{\eta} (\phi)$ and  $\alpha_{\eta} (\phi')$ are nowhere vanishing, and the global sections $\alpha_{\eta} (\psi)$, for all nonzero $\psi \in \CC\{\phi,\phi'\}$, have at most $0$-dimensional zero locus.  Then the cokernel of $c_{\phi,\phi'}$ is supported on at most $0$-dimensional subscheme.
\end{lem}
\begin{pf}
Assume that the support of the cokernel of $c_{\phi,\phi'}$ is $1$-dimensional and let $\Gamma$ be a reduced, irreducible component in the support of $coker(c_{\phi,\phi'})$. Restricting the diagram \eqref{diag:FphiFphi'} to $\Gamma$ implies that
$\FF_{[\phi]} \otimes \OO_{\Gamma} =\FF_{[\phi']} \otimes \OO_{\Gamma}$ or, equivalently, $\alpha_{\eta} (\phi)$ and $\alpha_{\eta} (\phi')$ are proportional along $\Gamma$. Arguing as in the proof of Proposition \ref{pro:Fphi}, 3), we deduce that
for some constant $\lambda \in \CC$ the global section
$\alpha_{\eta} (\phi'-\lambda \phi)$ vanishes on $\Gamma$ and this is contrary to the assumption on the pencil $\CC\{\alpha_{\eta} (\phi), \alpha_{\eta} (\phi')\}$.
\end{pf}
 
  Let $\phi, \phi'$ in $\HKX$ be subject to the condition of Lemma \ref{lem:cokercphiphi'} and define
$$
\FF_{[\phi],[\phi']}:= ker \Big(\FF_{[\phi]} \longrightarrow 
\OO_{E_{\phi,\phi'}} (2\KX-2E) \Big).
$$
From  this it follows that $\FF_{[\phi],[\phi']}$ is a torsion free subsheaf of $\FF_{[\phi]}$ of rank $2$ fitting into the exact sequence
\begin{equation}\label{seqFphiphi'}
\xymatrix@R=12pt@C=12pt{
0\ar[r]&\FF_{[\phi],[\phi']}\ar[r]& \FF_{[\phi]} \ar[r]& Im(c_{\phi,\phi'} ) \ar[r]&0,}
\end{equation}
where $Im(c_{\phi,\phi'} )$ is the image of the morphism $c_{\phi,\phi'}$ in \eqref{morph:cphiphi'}. From Lemma \ref{lem:cokercphiphi'} we have
\begin{equation}\label{Imcphiphi'}
\xymatrix@R=12pt@C=12pt{
0\ar[r]&Im(c_{\phi,\phi'} )\ar[r]& \OO_{E_{\phi,\phi'}} (2\KX-2E) \ar[r]& coker(c_{\phi,\phi'} ) \ar[r]&0,}
\end{equation}
where $coker(c_{\phi,\phi'} )$ is supported on at most $0$-dimensional subscheme of $X$.  
 Furthermore, the above exact sequence tells us that $Im(c_{\phi,\phi'} )$ has no subsheaves supported on $0$-dimensional subschemes. This in turn means that $\FF_{[\phi],[\phi']} $ is locally free (since otherwise taking the double dual $\FF^{\ast\ast}_{[\phi],[\phi']}$ of $\FF_{[\phi],[\phi']}$, we obtain the cokernel of the canonical
inclusion $coker(\FF_{[\phi],[\phi']} \longrightarrow  \FF^{\ast\ast}_{[\phi],[\phi']})$, a sheaf supported on a $0$-dimensional subscheme, injecting into $Im(c_{\phi,\phi'} )$). In addition,
from \eqref{Fphiglobsect} we deduce
$$
H^0 (\FF_{[\phi],[\phi']}) \cong H^0 (\FF_{[\phi]}).
$$


Recall that $\FF_{[\phi]}$ is also a quotient of $\TT_{\eta}$, see \eqref{Fphi}. Combining that sequence with the one in \eqref{seqFphiphi'} we obtain the diagram
\begin{equation}\label{diag:T'eta}
\xymatrix@R=12pt@C=12pt{
&&0\ar[d]&0\ar[d]&\\
0\ar[r]&\OO_X \ar[r]\ar@{=}[d]&\TT'_{\eta}  \ar[r]\ar[d]&\FF_{[\phi],[\phi']} \ar[r] \ar[d]&0\\
0\ar[r]&\OO_X \ar[r]&\TT_{\eta}  \ar[r]\ar[d]&\FF_{[\phi]} \ar[r]\ar[d]&0\\
&&{Im(c_{\phi,\phi'} )} \ar@{=}[r]\ar[d]&{Im(c_{\phi,\phi'} )}\ar[d]&\\
&&0&0}
\end{equation}
where the sheaf $\TT'_{\eta}$ is defined as the kernel of the epimorphism
$\TT_{\eta} \longrightarrow {Im(c_{\phi,\phi'} )}$.
\begin{lem}\label{l:T'eta}
1) The sheaf $\TT'_{\eta}$ in \eqref{diag:T'eta} is locally free and
$H^0(\TT'_{\eta}) \cong H^0(\TT_{\eta})$.

2) The parametrization $\alpha_{\eta} : \HKX \cong  H^0(\TT_{\eta})$ induces a parametrization   
$$
\alpha'_{\eta} : \HKX \cong  H^0(\TT'_{\eta}).
$$
Under this parametrization the global sections $\alpha_{\eta} (\psi)$ with no zeros go to sections $\alpha'_{\eta} (\psi)$ having no zeros as well.

3) The sheaf $\TT'_{\eta}$ fits into the extension sequence
\begin{equation}\label{ext-T'}
\xymatrix@R=12pt@C=12pt{
0\ar[r]& {\cal P} \ar[r]&\TT'_{\eta}  \ar[r] & \OO_X (\KX)\ar[r] & 0,
}
\end{equation}
where ${\cal P}$ is locally free sheaf. That extension sequence is related to the one in \eqref{ext-eta} by the following commutative diagram
\begin{equation}\label{diag:ext-T'}
\xymatrix@R=12pt@C=12pt{
&0\ar[d]&0\ar[d]&&\\
0\ar[r]& {\cal P} \ar[r]\ar[d]&\TT'_{\eta}  \ar[r]\ar[d] & \OO_X (\KX)\ar[r]\ar@{=}[d] & 0\\
0\ar[r]& {\Omega_X (-E)} \ar[r]\ar[d]&\TT_{\eta}  \ar[r]\ar[d] & \OO_X (\KX)\ar[r] & 0\\
&{Im(c_{\phi,\phi'} )}\ar@{=}[r]\ar[d]&{Im(c_{\phi,\phi'} )}\ar[d]& & \\
&0&0&&
}
\end{equation}

\end{lem}
\begin{pf}
The first two parts, 1) and 2), are immediate from the construction of $\TT'_{\eta} $, the diagram \eqref{diag:T'eta} and the properties of $\FF_{[\phi],[\phi']} $. To see the part 3), we combine the middle column
in \eqref{diag:T'eta} with the extension sequence \eqref{ext-eta}:
\begin{equation}\label{diag1:ext-T'}
\xymatrix@R=12pt@C=12pt{
&&0\ar[d]&&\\
& &\TT'_{\eta}\ar[d]\ar[dr]&&\\
0\ar[r]& {\Omega_X (-E)} \ar[r]&\TT_{\eta}  \ar[r]\ar[d] & \OO_X (\KX)\ar[r] & 0\\
&&{Im(c_{\phi,\phi'} )}\ar[d]&&\\
&&0&&
}
\end{equation}
where the slanted arrow is the composition of the monomorphism of the column with the epimorphism of the row. From the isomorphism $H^0 (\TT'_{\eta}) \cong \HKX$ and the global generation of $\OO_X (\KX)$ it follows that the slanted arrow is an epimorphism. This implies that the above diagram can be completed
to the diagram \eqref{diag:ext-T'}, where the sheaf ${\cal P}$ is defined to be the kernel of the slanted arrow $\TT'_{\eta}\longrightarrow \OO_X (\KX)$ in \eqref{diag1:ext-T'}. Thus ${\cal P}$ is a second syzygy sheaf and hence locally free.
\end{pf}
 
\begin{rem}\label{rem:eta-eta'}
1) Under the identification
$H^1 ({\cal P}(-\KX)) \cong Ext^1 (\OO_X (\KX),{\cal P})$ the extension sequence \eqref{ext-T'} corresponds to a cohomology class $\eta'\in H^1 ({\cal P}(-\KX))$. Tensoring the diagram \eqref{diag:ext-T'} with $\OO_X (-\KX)$, we deduce that this cohomology class maps to the class $\eta \in H^1 (\Theta(-E))$
under the homomorphism
$$
H^1 ({\cal P}(-\KX))\longrightarrow H^1 (\Theta(-E))
$$
induced by the monomorphism ${\cal P}(-\KX) \longrightarrow \Theta(-E)$ in the left column in \eqref{diag:ext-T'} tensored with $\OO_X (-\KX)$.

2) From Lemma \ref{l:T'eta}, 1), it follows that the cohomology class $\eta'$ lies in the kernel of the cup-product
$$
H^1 ({\cal P}(-\KX))\longrightarrow \HKX^{\ast} \otimes H^1 ({\cal P}).
$$
\end{rem}
 
\vspace{0.2cm} 
The parametrization $\alpha'_{\eta}$ in Lemma \ref{l:T'eta}, 2), similarly to the parametrization $\alpha_{\eta}$, is used to define the higher order products
\begin{equation}\label{hprod-alpha'}
\mbox{$\alpha'^{(k)}_{\eta}: \bigwedge^k \HKX \longrightarrow H^0 (\bigwedge^k\TT'_{\eta})$, for $k=2,3$.} 
\end{equation}
By definition
$$
\mbox{$
\alpha'^{(3)}_{\eta}:\bigwedge^3 \HKX \longrightarrow H^0 (\bigwedge^3\TT'_{\eta})=H^0 (\OO_X (2\KX -2E -E_1))$.}
$$ 
By construction of $\TT'_{\eta}$, the linear maps $\alpha'^{(3)}_{\eta}$ and $\alpha^{(3)}_{\eta}$ are related by the formula
\begin{equation}\label{rela3-a'3}
\alpha^{(3)}_{\eta} = e_{\phi,\phi'} \alpha'^{(3)}_{\eta} ,
\end{equation}
where $e_{\phi,\phi'}$ is a global section of $\OO_X (E_1)$ corresponding to $E_{\phi,\phi'}$. This relation implies the following.

\begin{lem}\label{lem:fixed}
1) Let $\phi, \phi'$ be as in the diagram \eqref{diag:T'eta}. Then the linear subsystem $|Im(\alpha'^{(3)}_{\eta}(\phi,\phi',\bullet)|$ has no fixed part.
In particular, the linear subsystem $|Im(\alpha'^{(3)}_{\eta})| \subset |2\KX -2E -E_1|$ has no fixed part as well, or, equivalently, the sheaf $\TT'_{\eta}$ is generated by its global sections outside of a subscheme of codimension at least $2$.

2) The divisor $E_0:=E_{\phi,\phi'}$ is the fixed part of the linear subsystem $|Im(\alpha^{(3)}_{\eta})|$. In particular, that divisor is independent
of the choice of $\phi, \phi'$ used in defining the extension sequence \eqref{ext-T'}.
\end{lem}
\begin{pf}
The assertion 2) follows immediately from the formula \eqref{rela3-a'3} and the part 1) of the lemma. 

For the part 1), we go back to the diagram \eqref{diag:T'eta}. From the second exterior power of its top row it follows that the sheaf
$\FF_{[\phi],[\phi']}$ also fits into the exact sequence
$$
\xymatrix@R=12pt@C=26pt{
0\ar[r]&\FF_{[\phi],[\phi']} \ar[r]&\bigwedge^2 \TT'_{\eta}\ar[r]^(0.3){\wedge\alpha'_{\eta}(\phi) }& \OO_X (2\KX -2E -E_1) \ar[r]&0
}
$$
and hence the global sections of $\FF_{[\phi],[\phi']} $ admit the following identification
\begin{equation}\label{H0Fphi-phi'}
H^0 (\FF_{[\phi],[\phi']}) \cong \{ \alpha'^{(2)}_{\eta} (\phi, \psi) | \psi \in \HKX \} = Im (\alpha'^{(2)}_{\eta} (\phi, \bullet)).
\end{equation}
Viewing the sheaf $\FF_{[\phi],[\phi']} $ (resp. $\FF_{[\phi]}$) as a subsheaf of 
$\bigwedge^2 \TT'_{\eta}$ (resp. $\bigwedge^2 \TT_{\eta}$), the triangle in \eqref{Fphiglobsect} can be completed as follows.
$$
\xymatrix@R=12pt@C=26pt{
H^0 (\FF_{[\phi],[\phi']} )\ar[r]^(0.4){\wedge\alpha'_{\eta}(\phi) } \ar[d]& H^0 (\OO_X (2\KX -2E-E_1) \ar^{e_{\phi,\phi'}}[d]\\
H^0 (\FF_{[\phi]} )\ar[r]^(0.4){\wedge\alpha_{\eta}(\phi) }& H^0 (\OO_X (2\KX -2E)
}
$$
This and the identification in \eqref{H0Fphi-phi'} (resp. $H^0 (\FF_{[\phi]}) \cong Im (\alpha^{(2)}_{\eta} (\phi, \bullet))$) imply 
$$
\alpha^{(3)}_{\eta} (\phi,\phi',\psi)= e_{\phi,\phi'} \alpha'^{(3)}_{\eta} (\phi,\phi',\psi),
$$
 for every $\psi\in \HKX$. Furthermore, by definition the divisor
$E_{\phi,\phi'} =(e_{\phi,\phi'}=0)$ is the fixed part of $|Im(\alpha^{(3)}_{\eta} (\phi,\phi',\bullet)|$. Therefore $|Im(\alpha'^{(3)}_{\eta} (\phi,\phi',\bullet)|$ is the moving part of $|Im(\alpha^{(3)}_{\eta} (\phi,\phi',\bullet)|$ and hence has no fixed part.
\end{pf}

We summarize the above results in the following.
\begin{pro}\label{pro:modif-eta}
1) For a general pair of linearly independent sections $\phi, \phi' \in \HKX$, the fixed part
$E_{\phi,\phi'}$  of the linear system $|Im(\alpha^{(3)}_{\eta} (\phi,\phi',\bullet)|$ is independent of $\phi, \phi'$ and coincides with the fixed part $E_0$ of $|Im(\alpha^{(3)}_{\eta})|$. 

2) There is a locally free sheaf ${\cal P}$ fitting into the exact sequence
\begin{equation}\label{P-oMEGA-l}
\xymatrix@R=12pt@C=12pt{
0\ar[r]& {\cal P} \ar[r]& \Omega_X (-E)\ar[r]& {\cal L}  \ar[r]&0}
\end{equation}
where ${\cal L}=Im(c_{\phi,\phi'})$, see \eqref{diag:ext-T'}, is a sheaf of rank $1$ supported on $E_0$ and which is part of the following exact sequence
\begin{equation}\label{sheafL}
\xymatrix@R=12pt@C=12pt{
0\ar[r]&{\cal L} \ar[r]& \OO_{E_{0}} (2\KX-2E) \ar[r]& {\cal J} \ar[r]&0}
\end{equation}
with the sheaf  ${\cal J}$ ($= coker(c_{\phi,\phi'})$, see \eqref{Imcphiphi'}) supported on at most $0$-dimensional subscheme.
 
3) There is a distinguished cohomology class $\eta'\in H^1 ({\cal P}(-\KX))$
which maps to the class $\eta$ under the homomorphism
$H^1 ({\cal P}(-\KX)) \longrightarrow H^1 (\Theta_X (-E))$ induced by the monomorphism in the exact sequence \eqref{P-oMEGA-l} tensored with $\OO_X (-\KX)$. Furthermore, $\eta'$ lies in the kernel of the cup-product
$$
H^1 ({\cal P}(-\KX)) \longrightarrow \HKX^{\ast} \otimes H^1 ({\cal P}).
$$

4) The extension sequence
\begin{equation}\label{ext-eta'}
\xymatrix@R=12pt@C=12pt{
0\ar[r]& {\cal P} \ar[r]& \TT'_{\eta}\ar[r]& \OO_X (\KX) \ar[r]&0}
\end{equation}
defined by $\eta'$ under the identification $H^1 ({\cal P}(-\KX))\cong Ext^1 (\OO_X (\KX),{\cal P})$,
gives rise to a locally free sheaf $\TT'_{\eta}$ with the isomorphism
$H^0 (\TT'_{\eta}) \cong \HKX$ induced by the epimorphism of the extension sequence. Furthermore, the extension sequence above is related to the one in \eqref{ext-eta} by the commutative diagram in \eqref{diag:ext-T'}.

5) Let $\alpha'_{\eta}: \HKX\longrightarrow H^0 (\TT'_{\eta}) $ be the inverse of the isomorphism in 4) and let
$\alpha'^{(k)}_{\eta} : \bigwedge^k \HKX \longrightarrow H^0 (\bigwedge^k\TT'_{\eta}) $, $k=2,3$ be the higher products defined in \eqref{hprod-alpha'}. Then 

\vspace{0.2cm}
(i) for a general $\phi \in \HKX$ the global section $\alpha'_{\eta} (\phi)$ of $\TT'_{\eta}$ has no zeros,

(ii) $\alpha'^{(3)}_{\eta} (\phi,\phi',\phi'') \neq 0$, for any linearly independent triple $\phi,\phi',\phi'' \in \HKX$ for which the linear subsystem of $|\KX|$
generated by $\phi,\phi',\phi'' $ has at most $0$-dimensional base locus,

(iii) any subsheaf of rank $1$ of $\TT'_{\eta}$ has at most $1$ dimensional space of global sections; in particular, $\alpha'^{(2)}_{\eta} (\phi,\phi')=0$
iff $\phi,\phi'$ are linearly dependent,

(iv)  for a general pair $\phi, \phi' \in \HKX$ the linear system $|Im(\alpha'^{(3)}_{\eta} (\phi,\phi',\bullet))|$ is the moving part of the linear system $|Im(\alpha^{(3)}_{\eta} (\phi,\phi',\bullet))|$,

(v) $\TT'_{\eta}$ is generated by its global sections outside of a subscheme of codimention $\geq 2$.  
\end{pro}

As for $\TT_{\eta}$, we associate with the sheaf $\TT'_{\eta}$ and the parametrization $\alpha'_{\eta}$ in Proposition \ref{pro:modif-eta}, 5), the family $\{\FF'_{[\phi]} \}_{[\phi]}$ of rank $2$ sheaves by setting
\begin{equation}\label{shF'phi}
\FF'_{[\phi]}:=coker(\OO_X \stackrel{a'_{\eta} (\phi)}{\longrightarrow} \TT'_{\eta}),
\end{equation}
for every nonzero $\phi \in \HKX$. These sheaves are related to $\FF_{[\phi]}$'s by a modification along the divisor $E_0$ in Proposition \ref{pro:modif-eta}, 1). Namely, we have the following commutative diagram
$$
 \xymatrix@R=12pt@C=12pt{
& &0\ar[d]& 0\ar[d]\\
0\ar[r]& {\OO_X} \ar[r]^{\alpha'_{\eta}(\phi)} \ar@{=}[d]& \TT'_{\eta}\ar[r] \ar[d]& \FF'_{[\phi]} \ar[r] \ar[d]&0\\
0\ar[r]& {\OO_X} \ar[r]^{\alpha_{\eta}(\phi)}& \TT_{\eta}\ar[r] \ar[d]& \FF_{[\phi]} \ar[r] \ar[d]&0\\
 & & {\cal L}\ar[d] \ar@{=}[r]& {\cal L}\ar[d]\\
& &0 &0}
$$
for every nonzero $\phi\in \HKX$, and where the middle column is as in the diagram \eqref{diag:ext-T'}.  

The following properties of the sheaves $\FF'_{[\phi]} $'s are similar to the ones
found for $\FF_{[\phi]} $'s in Proposition \ref{pro:Fphi}.
\begin{pro}\label{pro:F'phi}
 1) For every nonzero $\phi\in \HKX$ the sheaf $\FF'_{[\phi]} $ has rank $2$ and the Chern invariants
$$
c_1 (\FF'_{[\phi]})=2\KX-2E-E_0, \hspace{0.2cm} c_2(\FF'_{[\phi]})=c_2 (\FF_{[\phi]})+\KX \cdot E_0 -deg({\cal J}),
$$
where ${\cal J}$ is as in \eqref{sheafL}.
  
2) $H^0 (\FF'_{[\phi]}) \cong \HKX / \CC\phi$.

3) The sheaf $\FF'_{[\phi]} $ is locally free if and only if the global section $\alpha'_{\eta}(\phi)$ has no zeros. For such a sheaf $\FF'_{[\phi]} $ and a general $\psi \in  \HKX$, the global section $f_{\psi}$ of $\FF'_{[\phi]}$ corresponding to $\psi$ under the isomorphism in 2) has $0$-dimensional scheme of zeros $Z_{f_{\psi}}$ and gives rise to the Koszul sequence
$$
 \xymatrix@R=12pt@C=35pt{
0\ar[r]& \OO_X \ar[r]^{f_{\psi}}&\FF'_{[\phi]}  \ar[r]^(0.26){\wedge f_{\psi}}& {\cal I}_{Z_{f_{\psi}}} (2\KX-2E-E_0)\ar[r]&0.}
$$
The homomorphism $\wedge f_{\psi}: H^0 (\FF'_{[\phi]}) \longrightarrow H^0 ({\cal I}_{Z_{f_{\psi}}} (2\KX-2E-E_0))$ induced by the epimorphism in the Koszul sequence above is given by the map $\alpha'^{(3)}_{\eta} (\phi,\psi,\bullet)$. More precisely, one has a commutative diagram
$$
\xymatrix@R=12pt@C=37pt{
H^0 (\FF'_{[\phi]}) \ar[r]^(0.3){\wedge f_{\psi}} \ar@{=}[d]&  H^0 ({\cal I}_{Z_{f_{\psi}}} (2\KX-2E-E_0))\ar@{^{(}{-}{>}}[d]\\
\HKX /\CC\phi \ar^(0.45){\mbox{\tiny$\alpha'^{(3)}_{\eta} \!(\phi,\!\psi,\!\bullet)$}}[r]& H^0 (\OO_X(2\KX-2E-E_0))
}
$$
where the identification on the left is the isomorphism in 2). Thus, for a general linearly independent pair $\phi, \psi \in \HKX$, one has an identification 
$$
 H^0 ({\cal I}_{Z_{f_{\psi}}} (2\KX-2E-E_0))=Im (\alpha'^{(3)}_{\eta} \!(\phi,\!\psi,\!\bullet)),
$$
 implying that the linear system $|{\cal I}_{Z_{f_{\psi}}} (2\KX-2E-E_0)|$ has at most $0$-dimensional base locus.
\end{pro}

In addition to the above, the locally free sheaves $\FF'_{[\phi]}$ are subject to the following. 
\begin{pro}\label{pro:F'phi-gen-nopencil}
1) For $\phi \in \HKX$ with $\alpha'_{\eta} (\phi)$ nowhere vanishing, the locally free sheaf $\FF'_{[\phi]}$ is generated by its global sections outside of at most $0$-dimensional subscheme of $X$. 

2) Set $Y:= \PP(\FF'^{\ast}_{[\phi]})$, the projectivization of $\FF'^{\ast}_{[\phi]}$, and define $\OO_Y (1)$ so that the direct image $\pi_{\ast} (\OO_Y (1)) =\FF'_{[\phi]}$, where $\pi: Y= \PP(\FF'^{\ast}_{[\phi]}) \longrightarrow X$ is the structure projection. Then $\OO_Y (1)$ defines the (rational) map
$$
\widetilde{\gamma}: Y --\rightarrow \PP(H^0 (\FF'_{[\phi]})^{\ast}) \cong \PP^{p_g-2}
$$
whose locus of indeterminacy is contained in a finitely many fibres of $\pi$.
Furthermore, if $p_g \geq 5$, then the image of the above map is $3$-dimensional. 

3) Under the natural identifications
\begin{equation}\label{ident-sections}
H^0 (\OO_{\PP(H^0 (\FF'_{[\phi]})^{\ast})} (1) ) \cong H^0 (\OO_Y (1)) \cong H^0 (\FF'_{[\phi]})\cong \HKX /\CC\phi 
\end{equation}
every nonzero $\psi \in \HKX$ linearly independent of $\phi$ corresponds to a unique hyperplane $\widetilde{H}_{\psi}$ in 
$\PP(H^0 (\FF'_{[\phi]})^{\ast})$. For a general such $\psi$, with the assumption $p_g \geq 5$, the hyperplane section
$H_{\psi} =\widetilde{H}_{\psi} \bigcap \widetilde{\gamma} (Y)$ is a reduced irreducible surface which is the image of $X$ under the rational map defined by the linear system $|Im (\alpha'^{(3)}_{\eta} \!(\phi,\!\psi,\!\bullet))|$.
In particular, for general $\psi \in \HKX$ the linear system $|Im (\alpha'^{(3)}_{\eta} \!(\phi,\!\psi,\!\bullet))|$ is not composed of a pencil and hence its general member is irreducible.
\end{pro}
\begin{pf}
The sheaf $\FF'_{[\phi]}$ is defined by the exact sequence
$$
\xymatrix@R=12pt@C=12pt{
0\ar[r]&\OO_X \ar[r]^{\alpha'_{\eta} (\phi)} \ar[r]& \TT'_{\eta} \ar[r]&\FF'_{[\phi]}\ar[r]&0.
}
$$
From this we deduce the following  commutative diagram
$$
\xymatrix@R=12pt@C=12pt{
0\ar[r]&\OO_X  \ar[r] \ar@{=}[d]&H^0 ( \TT'_{\eta})\otimes \OO_X \ar[r]\ar[d]&H^0 (\FF'_{[\phi]})\otimes \OO_X \ar[d] \ar[r] &0\\
0\ar[r]&\OO_X \ar[r]^{\alpha'_{\eta} (\phi)} \ar[r]& \TT'_{\eta} \ar[r]&\FF'_{[\phi]}\ar[r]&0.
}
$$
This implies that the two vertical arrows have isomorphic cokernels. By Proposition \ref{pro:modif-eta}, $(v)$, the cokernel of the middle vertical arrow is supported on at most $0$-dimensional subscheme. Hence the same holds for the right vertical arrow. This proves the assertion 1) of the proposition.

Set $Z$ to be the support of the cokernel of the evaluation morphism
\begin{equation}\label{evF'phi}
H^0 (\FF'_{[\phi]})\otimes \OO_X \longrightarrow \FF'_{[\phi]}.
\end{equation}
According to 1), the subscheme $Z$ is at most $0$-dimensional and  $\FF'_{[\phi]}$ is globally generated outside of $Z$. By definition $\OO_Y (1)$ is globally generated outside $\pi^{\ast}(Z)$, a subscheme supported on a finite number of fibres of the structure projection $\pi: Y=\PP(\FF'^{\ast}_{[\phi]})\longrightarrow X$.

We now turn to the statement about the image of the map $\widetilde{\gamma}$
defined by $\OO_Y (1)$.
This map admits the following description. The dual of the evaluation morphism
defines the morphism
\begin{equation}\label{gamma-mor}
\gamma: X\setminus Z \longrightarrow Gr(1,\PP(H^0 (\FF'_{[\phi]})^{\ast})),
\end{equation} 
where $Gr(1,\PP(H^0 (\FF'_{[\phi]})^{\ast}))$ is the Grassmannian of lines in $\PP(H^0 (\FF'_{[\phi]})^{\ast})$ and $\gamma$ sends a point $x \in  X\setminus Z$ to the line 
$\PP(\FF'^{\ast}_{[\phi]}(x))$ in $\PP(H^0 (\FF'_{[\phi]})^{\ast})$, where
$\FF'^{\ast}_{[\phi]}(x)$ denotes the fibre of $\FF'^{\ast}_{[\phi]}$ at $x$. 
This tells us that the universal bundle ${\cal U}$ on $Gr(1,\PP(H^0 (\FF'_{[\phi]})^{\ast}))$
pulls back by $\gamma$ to $\FF'^{\ast}_{[\phi]}$, i.e., $\gamma^{\ast} ({\cal U})=\FF'^{\ast}_{[\phi]}$. The map $\widetilde{\gamma}$ defined by  $\OO_Y (1)$ is the composition of the arrows in the top row of the following diagram
\begin{equation}\label{diag:grass}
\xymatrix@R=12pt@C=12pt{
\PP(\FF'^{\ast}_{[\phi]}) \ar@{-->}[r] \ar[d]& \PP({\cal U}|_{\overline{\Sigma}} ) \ar@{^{(}->}[r] \ar[d]& Gr(1,\PP(H^0 (\FF'_{[\phi]})^{\ast})) \times \PP(H^0 (\FF'_{[\phi]})^{\ast}) \ar[r] \ar[d]& \PP(H^0 (\FF'_{[\phi])})^{\ast})\\
X \ar@{-->}[r]^{\gamma}& \overline{\Sigma}  \ar@{^{(}->}[r]&Gr(1,\PP(H^0 (\FF'_{[\phi]})^{\ast})) &
}
\end{equation}
where $\overline{\Sigma}$ is the closure of the image of the morphism $\gamma$
in \eqref{gamma-mor}. 
\begin{cl}\label{cl:im-gamma}
The image $\Sigma$ of $\gamma$ in \eqref{gamma-mor} is $2$-dimensional.
\end{cl}
\begin{pf}
The vector bundle $\FF'_{[\phi]}$ is obviously nontrivial, since
$h^0 (\FF'_{[\phi]}) =p_g-1 \geq 4-1=3$. So $\gamma$ is not a constant map.
Assume that the image $\Sigma$ of $\gamma$  is a curve. Then, in view of \eqref{cond}, $(ii)$, it must be a rational curve. So after resolving the indeterminacy locus $Z$ we obtain a surjective morphism $f:\widetilde{X} \longrightarrow \PP^1$ with connected fibres, where $b:\widetilde{X} \longrightarrow X$ is a sequence of of blowing down maps. Let $\widetilde{F}$ denote the rational equivalence class of the fibres of $f$. Then $\OO_{\widetilde{X}}(\widetilde{F}) =f^{\ast} (\OO_{\PP^1} (1))$. Furthermore, the induced map 
$\widetilde{\gamma}:\widetilde{X}\longrightarrow Gr(1,\PP(H^0 (\FF'_{[\phi]})^{\ast}))$ is constant on the fibres of $f$. This means that the pull back
$b^{\ast} (\FF'_{[\phi]})$ must be trivial on the fibres of $f$. In particular,
 $h^0 (b^{\ast} (\FF'_{[\phi]}) \otimes \OO_{\widetilde{F}_t}) =2$, for every reduced irreducible fibre $\widetilde{F}_t =f^{-1} (t)$, $t\in \PP^1$, of $f$. From the exact sequence
$$
\xymatrix@R=12pt@C=12pt{
0\ar[r]&b^{\ast} (\FF'_{[\phi]})(-\widetilde{F}) \ar[r]&b^{\ast} (\FF'_{[\phi]})\ar[r]&b^{\ast} (\FF'_{[\phi]}) \otimes \OO_{\widetilde{F}_t} \ar[r]&0
}
$$
follows the estimate
$$
h^0 ( b^{\ast} (\FF'_{[\phi]})(-\widetilde{F}) )\geq h^0 ( b^{\ast} (\FF'_{[\phi]})) -h^0 (b^{\ast} (\FF'_{[\phi]}) \otimes \OO_{\widetilde{F}_t})=
h^0 (\FF'_{[\phi]})-2=p_g-1-2=p_g -3 \geq 1.
$$ 
Setting $F:=b_{\ast} (\widetilde{F})$ we have the identification $b^{\ast} ({\cal J}_Z (F))=\OO_{\widetilde{X}} (\widetilde{F})$, where ${\cal J}_Z$ is the ideal sheaf of the indeterminacy locus $Z$ of $\gamma$. Hence we obtain
 $$
\begin{gathered}
H^0 ( (b^{\ast} (\FF'_{[\phi]})(-\widetilde{F}))=Hom_{\OO_{\widetilde{X}}}(\OO_{\widetilde{X}} (\widetilde{F}), b^{\ast} (\FF'_{[\phi]}))=
 Hom_{\OO_{\widetilde{X}}}(b^{\ast} ({\cal J}_Z (F)), b^{\ast} (\FF'_{[\phi]}))=\\Hom_{\OO_{X}}({\cal J}_Z (F), b_{\ast} (b^{\ast} (\FF'_{[\phi]})))=Hom_{\OO_{X}}({\cal J}_Z (F), \FF'_{[\phi]}) =Hom_{\OO_{X}}(\OO_X (F), \FF'_{[\phi]}),
\end{gathered}
$$
where the third equality uses the adjoint property of the functors $b^{\ast},b_{\ast}$. Thus we obtain a nonzero morphism

$$
\OO_X ({F}) \longrightarrow \FF'_{[\phi]}.
$$
Since $h^0 (\OO_X (F)) \geq h^0 ({\cal J}_Z (F) ) =
h^0 (b_{\ast}(\OO_{\widetilde{X}} (\widetilde{F}))) =h^0 (\OO_{\widetilde{X}}(\widetilde{F}))=2$, the above morphism gives rise to two linearly independent global sections of $\FF'_{[\phi]}$ which are proportional. In view of Proposition \ref{pro:F'phi}, 3), this means that there exist linearly independent $\psi,\psi' \in \HKX$ such that $\phi,\psi,\psi'$ are linearly independent and $\alpha'^{(3)}_{\eta} (\phi,\psi,\psi')=0$ and this is contrary to Proposition \ref{pro:modif-eta}, 5), $(ii)$.
\end{pf}
Once we know that $\overline{\Sigma}$ in the diagram \eqref{diag:grass} is a surface, the $\PP^1$-bundle $\PP({\cal U}|_{\overline{\Sigma}} )$ over it is a $3$-fold. The image of $\PP({\cal U}|_{\overline{\Sigma}} ) \subset Gr(1,\PP(H^0 (\FF'_{[\phi]})^{\ast})) \times \PP(H^0 (\FF'_{[\phi]})^{\ast})$ under the projection to $\PP(H^0 (\FF'_{[\phi]})^{\ast}) \cong \PP^{p_g -2}$ is irreducible and spans $\PP^{p_g -2}$. Hence
it must be at least of dimension $2$. Furthermore, the dimension $2$ can 
only occur if $\PP^{p_g -2}$ is a plane, i.e. $p_g =4$. Thus the image of
$\widetilde{\gamma}: Y=\PP(\FF'^{\ast}_{[\phi]}) --\rightarrow \PP(H^0 (\FF'_{[\phi]})^{\ast})$ is a $3$-dimensional irreducible variety,
provided $p_g \geq 5$. This completes the proof of part 2).

In the part 3) of the proposition
the first identification from the left is obvious, the second comes from
$\pi_{\ast} (\OO_Y (1)) =\FF'_{[\phi]}$ and the third is Proposition \ref{pro:F'phi}, 2).

We now assume $p_g \geq 5$. By 2) of the proposition the image $\widetilde{\gamma} (Y)$ of $\widetilde{\gamma}$ is a $3$ dimensional irreducible subvariety of $\PP(H^0 (\FF'_{[\phi]})^{\ast})$. Hence for a general hyperplane $\widetilde{H}_{\psi} $, the hyperplane section $H_{\psi} =\widetilde{H}_{\psi} \bigcap \widetilde{\gamma} (Y)$ is a reduced irreducible surface. Our task is to identify that surface as the image of the map defined by the linear system
$|Im (\alpha'^{(3)}_{\eta} \!(\phi,\!\psi,\!\bullet))|$.

 By Proposition \ref{pro:F'phi}, 3), for a general $\psi \in \HKX$ linearly independent of $\phi$, we have the identification
\begin{equation}\label{im(a3)=If-psi}
Im (\alpha'^{(3)}_{\eta} \!(\phi,\!\psi,\!\bullet))\cong 
H^0 ({\cal I}_{Z_{f_{\psi}}} (2\KX-2E-E_0)),
\end{equation}
  where 
$f_{\psi}$ is the global section of $\FF'_{[\phi]}$ corresponding to $\psi$ under the first isomorphism from the right in \eqref{ident-sections},
$Z_{f_{\psi}} =(f_{\psi} =0)$ and ${\cal I}_{Z_{f_{\psi}}}$ is the ideal sheaf of 
$Z_{f_{\psi}}$. By the identification
$H^0 (\FF'_{[\phi]} ) \cong H^0 (\OO_Y (1))$ in \eqref{ident-sections}, the section $f_{\psi}$ corresponds to a global section of $\OO_Y (1)$ which we denote by $s_{\psi}$. The divisor
$D_{s_{\psi}} =(s_{\psi}=0)$ is a reduced irreducible surface\footnote{$D_{s_{\psi}} $ is birationally isomorphic to $X$ - the projection $\pi$ restricted to $D_{s_{\psi}} $ exhibits that surface as the blow-up of $X$ along $Z_{f_{\psi}}$.} in $Y$ which is mapped by
$\widetilde{\gamma}$ onto the hyperplane section $H_{\psi}$, i.e.,
 the rational map
$\widetilde{\gamma}|_{D_{s_{\psi}}} :D_{s_{\psi}} --\rightarrow H_{\psi}$ is defined by $\OO_{D_{s_{\psi}}} (1)$, the restriction of $\OO_Y (1)$ to $D_{s_{\psi}}$.  Furthermore, the invertible sheaf $\OO_{D_{s_{\psi}}} (1)$ is related to ${\cal I}_{Z_{f_{\psi}}} (2\KX-2E-E_0)$ by the following isomorphism
\begin{equation}\label{O(1)=I}
\pi_{\ast} (\OO_{D_{s_{\psi}}} (1)) \cong {\cal I}_{Z_{f_{\psi}}} (2\KX-2E-E_0).
\end{equation} 
This is seen by taking the exact sequence
$$
\xymatrix@R=12pt@C=12pt{
0\ar[r]& \OO_Y \ar[r]^(0.4){s_{\psi}} &\OO_Y (1)\ar[r]& \OO_{D_{s_{\psi}}} (1)\ar[r]&0
}
$$
on $Y$ and applying to it the direct image $\pi_{\ast}$ to obtain the exact sequence
$$
\xymatrix@R=12pt@C=12pt{
0\ar[r]& \OO_X \ar[r]^{f_{\psi}} &\FF'_{[\phi]} \ar[r]& \pi_{\ast} (\OO_{D_{s_{\psi}}} (1))\ar[r]&0,
}
$$
which is the Koszul sequence of $(\FF'_{[\phi]}, f_{\psi})$. Hence the identification \eqref{O(1)=I}. This identification combined with \eqref{im(a3)=If-psi}
 implies that the linear system $|H^0 ({\cal I}_{Z_{f_{\psi}}} (2\KX-2E-E_0))|$ defines a rational map
$X --\rightarrow H_{\psi}$. Since $H_{\psi}$ is $2$-dimensional and the linear system $|H^0 ({\cal I}_{Z_{f_{\psi}}} (2\KX-2E-E_0))|$ has no fixed part, we deduce that $|H^0 ({\cal I}_{Z_{f_{\psi}}} (2\KX-2E-E_0))|$ is not composed of pencil and its general member is irreducible.
\end{pf}

\section{Sheaves ${\cal F}_{[\phi]}$ and ${\cal F'}_{[\phi]}$}
Set $U$ to be the subset of $\PP(\HKX)$ parametrizing global sections of $\TT_{\eta}$ without zeros. By Proposition \ref{pro:glsecTeta} this is a Zariski dense open subset of  $\PP(\HKX)$. In view of Lemma \ref{l:T'eta}, 2), the same open subset parametrizes global sections of $\TT'_{\eta}$ without zeros.  By Proposition \ref{pro:Fphi}, 3), (resp. Proposition \ref{pro:F'phi}, 3)) that set also parametrizes the sheaves ${\cal F}_{[\phi]}$ (resp. ${\cal F}'_{[\phi]}$) which are locally free. We will need some further properties of these sheaves.

\begin{pro}\label{pro:FphiOmega}
For $[\phi]\in U$ the locally free sheaf ${\cal F}_{[\phi]}$ (resp. ${\cal F'}_{[\phi]}$)  admits the following inclusions
$$
\Omega_X(-E) \hookrightarrow {\cal F}_{[\phi]} \hookrightarrow \Omega_X(K_X-E) \hspace{0.2cm} (\mbox{ resp. \,\,${\cal P} \hookrightarrow {\cal F'}_{[\phi]} \hookrightarrow {\cal P}(K_X)$}),
$$
where ${\cal P}$ is as in the extension sequence \eqref{ext-T'}.
Furthermore, those inclusions can be completed to the following exact sequences:

\begin{equation}\label{OmFphi-1}
\xymatrix@R=12pt@C=12pt{ 
0\ar[r] &\Omega_X(-E)\ar[r] &{\cal F}_{[\phi]} \ar[r]&\OO_{C_{\phi}} (K_X) \ar[r]&0\\
0\ar[r]&{\cal F}_{[\phi]} \ar[r] &\Omega_X(K_X-E)\ar[r] &\OO_{C_{\phi}} (2K_X -2E) \ar[r]&0,\\
0\ar[r] &{\cal P} \ar[r] &{\cal F'}_{[\phi]} \ar[r]&\OO_{C_{\phi}} (K_X) \ar[r]&0\\
0\ar[r]&{\cal F'}_{[\phi]} \ar[r] &{\cal P} (\KX) \ar[r] &\OO_{C_{\phi}} (2K_X -2E-E_0) \ar[r]&0,
}
\end{equation}
where $C_{\phi} =(\phi=0)$. In addition, the two pairs of exact sequences are related by the following commutative diagrams
\begin{equation}\label{OmFphi-2}
\xymatrix@R=12pt@C=12pt{
&0\ar[d] &0\ar[d] & &\\ 
0\ar[r] &{\cal P} \ar[r] \ar[d] &{\cal F'}_{[\phi]} \ar[r]\ar[d]&\OO_{C_{\phi}} (K_X) \ar[r] \ar@{=}[d]&0\\
0\ar[r] &\Omega_X(-E)\ar[r] \ar[d]&{\cal F}_{[\phi]} \ar[r] \ar[d]&\OO_{C_{\phi}} (K_X) \ar[r]&0\\ 
& {\cal L}\ar@{=}[r]\ar[d] &{\cal L}\ar[d] & &\\
&0 &0& &}
\end{equation}
\begin{equation}\label{OmFphi-3}
\xymatrix@R=12pt@C=12pt{
&0\ar[d] &0\ar[d] & 0\ar[d] &\\ 
0\ar[r] &{\cal F'}_{[\phi]} \ar[r] \ar[d] &{\cal P}(\KX) \ar[r]\ar[d]&\OO_{C_{\phi}} (2K_X -2E-E_0) \ar[r] \ar[d]&0\\
0\ar[r] &{\cal F}_{[\phi]}\ar[r] \ar[d]&\Omega_X(\KX-E) \ar[r] \ar[d]&\OO_{C_{\phi}} (2K_X -2E) \ar[r] \ar[d]&0\\ 
0\ar[r]& {\cal L}\ar[r]\ar[d] &{\cal L}\otimes \OO_X (\KX)\ar[d]\ar[r] &\OO_{C_{\phi} \cdot E_0} (2K_X -2E) \ar[r] \ar[d]&0\\
&0 &0&0 &}
\end{equation}
\end{pro}
\begin{pf}
Consider the Koszul sequence defined by the section $\alpha_{\eta} (\phi)$
$$
\xymatrix@R=12pt@C=28pt{ 
0\ar[r] &\OO_X\ar[r]^{\alpha_{\eta} (\phi)} &\TT_{\eta} \ar[r]^{\wedge\alpha_{\eta} (\phi)}& \bigwedge^2 \TT_{\eta}\ar[r]^(.4){\wedge\alpha_{\eta} (\phi)}&\OO_X (2K_X -2E) \ar[r]&0.
}
$$
For $\alpha_{\eta} (\phi)$ with no zeros the above complex is exact. From this it follows that ${\cal F}_{[\phi]} $, the cokernel of the first morphism on the left, is isomorphic to the kernel of the second arrow from the right. Namely, the above Koszul complex breaks into the following two short exact sequences
{\scriptsize
\begin{equation}\label{Kosz-twoseq}
\xymatrix@R=12pt@C=20pt{
0\ar[r] &\OO_X\ar[r]^{\alpha_{\eta} (\phi)} &\TT_{\eta} \ar[r]&{\cal F}_{[\phi]} \ar[r]&0\\
0\ar[r] &{\cal F}_{[\phi]} \ar[r]& \bigwedge^2 \TT_{\eta}\ar[r]^(.4){\wedge\alpha_{\eta} (\phi)}&\OO_X (2K_X -2E) \ar[r]&0.
}
\end{equation}
}
Combining the first sequence with the defining extension sequence in \eqref{ext-eta} gives the diagram
{\scriptsize
\begin{equation}\label{diag-twoseq-1}
\xymatrix@R=12pt@C=20pt{
& &0\ar[d]& &\\
& &\Omega_X(-E) \ar[d] \ar[dr]& &\\
0\ar[r] &\OO_X\ar[r]^{\alpha_{\eta} (\phi)} \ar[dr]_(.4){\phi} &\TT_{\eta} \ar[r] \ar[d]&{\cal F}_{[\phi]} \ar[r]&0\\
& &\OO_X (K_X) \ar[d]& &\\
& &0& &
}
\end{equation}
}
where the slanted arrow in the lower left corner of the diagram is the multiplication by the section $\phi$ of $\OO_X (K_X)$. Hence the slanted arrow in the upper right corner of the diagram can be completed to the first exact sequence in
\eqref{OmFphi-1}.

Combining the second exact sequence in \eqref{Kosz-twoseq} with the second exterior power of \eqref{ext-eta} gives rise to the diagram
{\scriptsize
\begin{equation}\label{diag-twoseq-2}
\xymatrix@R=12pt@C=20pt{
& &0\ar[d]& &\\
& &\OO_X(K_X-2E) \ar[d] \ar[dr]^{\phi}& &\\
0\ar[r] &{\cal F}_{[\phi]} \ar[r] \ar[dr] &\bigwedge^2 \TT_{\eta} \ar[r]^(.4){\wedge\alpha_{\eta} (\phi)} \ar[d]&{\cal O}_X (2K_X -2E) \ar[r]&0\\
& &\Omega_X (K_X -E) \ar[d]& &\\
& &0& &
}
\end{equation}
}
where the slanted arrow in the upper part of the diagram is the multiplication by $\phi$. Hence the slanted arrow in the lower part can be completed to the second exact sequence in \eqref{OmFphi-1}.

Repeating the same argument for the global section $\alpha'_{\eta} (\phi)$ of $\TT'_{\eta}$
one obtains the assertions involving $\FF'_{[\phi]}$.

The diagrams \eqref{OmFphi-2} and \eqref{OmFphi-3} are obtained by exploiting the diagram \eqref{diag:ext-T'} and the exact sequences in \eqref{OmFphi-1}.
\end{pf}

\begin{cor}\label{Fphi-Cphi}
Let  ${\cal F}_{[\phi]}$ and $\FF'_{[\phi]}$ be as in Proposition \ref{pro:FphiOmega}. Then the following hold.
 
1) The restriction of ${\cal F}_{[\phi]}$ to $C_{\phi}$ fits into the exact sequence
{\scriptsize
$$
\xymatrix@R=12pt@C=12pt{
0\ar[r]& \OO_{C_{\phi}} (K_X -2E) \ar[r]&{\cal F}_{[\phi]} \otimes \OO_{C_{\phi}}\ar[r]&  \OO_{C_{\phi}}(K_X) \ar[r]&0.
}
$$
}
2) The restriction of ${\cal F'}_{[\phi]}$ to $C_{\phi}$ fits into the exact sequence
{\scriptsize
$$
\xymatrix@R=12pt@C=12pt{
0\ar[r]& \OO_{C_{\phi}} (K_X -2E-E_0) \ar[r]&{\cal F'}_{[\phi]} \otimes \OO_{C_{\phi}}\ar[r]&  \OO_{C_{\phi}}(K_X) \ar[r]&0.
}
$$
}

3) The two exact sequences above are related by the commutative diagram
{\scriptsize
$$
\xymatrix@R=12pt@C=12pt{
 &0\ar[d] &0\ar[d] & &\\
0\ar[r]& \OO_{C_{\phi}} (K_X -2E-E_0) \ar[r] \ar[d]&{\cal F'}_{[\phi]} \otimes \OO_{C_{\phi}}\ar[r] \ar[d]&  \OO_{C_{\phi}}(K_X) \ar[r] \ar@{=}[d]&0\\
0\ar[r]& \OO_{C_{\phi}} (K_X -2E) \ar[r] \ar[d]&{\cal F}_{[\phi]} \otimes \OO_{C_{\phi}}\ar[r] \ar[d]&  \OO_{C_{\phi}}(K_X) \ar[r]&0\\
&\OO_{C_{\phi} \cdot E_0} (K_X -2E) \ar@{=}[r] \ar[d]& {\cal L} \otimes \OO_{C_{\phi}} \ar[d]\\
&0 &0&&
}
$$
}
\end{cor}
\begin{pf}
The sequence in 1) (resp. 2)) follows from the epimorphism in the first (resp. third) exact sequence in \eqref{OmFphi-1}  restricted to $C_{\phi}$ and the fact that 
$\det({\cal F}_{[\phi]}) =\OO_X (2K_X-2E)$ (resp. $\det({\cal F'}_{[\phi]})=\OO_X (2K_X-2E-E_0)$, see Proposition \ref{pro:Fphi}, 1) (resp. Proposition \ref{pro:F'phi}, 1)).

The diagram in 3) is obtained from restricting \eqref{OmFphi-2} to $C_{\phi}$ and
and then using 1) and 2).
\end{pf}

\section{The condition of isotropy of $\xi$}

We return to the cohomology class $\xi$ in the kernel of the cup-product \eqref{cp} fixed at the outset and explain how the isotropy condition with respect to the quadratic form
$\delta_X$ in \eqref{q-form} enters the study of the extension \eqref{ext}. For this we recall the commutative triangle
$$
\xymatrix@R=12pt@C=8pt{
Sym^2(H^1 (\TET)) \ar[rr]^{\delta_X}  \ar[rd]_(0.4){p^{(2)}}& &
H^2(\OO_X (-K_X)\ar[dl]^(0.45){m^{\ast}_2}\\
& ( Sym^2(\HKX))^{\ast} &
}
$$
discussed in the introduction, see \eqref{diag-comtriang} for notation.

The fact that $\xi$ lies in the kernel of the cup-product $p$ implies that 
$p^{(2)} (\xi^2)=0$. The commutativity of the diagram above implies that 
$\delta_X (\xi^2) \in H^2 (\OO_X (-\KX))$ is in the kernel of $m^{\ast}_2$.
For our purposes it is better to write that map in the following form
\begin{equation}\label{m2dual-Koszul}
H^2(\OO_X (-K_X)) \longrightarrow \HKX^{\ast} \otimes\HKX^{\ast} \cong \HKX^{\ast}  \otimes H^2 (\OO_X),
\end{equation}
where the second identity comes from the Serre duality $\HKX^{\ast}  \cong H^2 (\OO_X)$. The kernel of this map can now be identified with a certain Koszul cohomology. Namely, set $W=\HKX$ and consider the exact Koszul complex
$$
\xymatrix@R=12pt@C=8pt{ 
0\ar[r]&\OO_X (-K_X) \ar[r]&W^{\ast} \otimes\OO_X  \ar[r]&\bigwedge^2 W^{\ast} \otimes\OO_X (\KX)  \ar[r]&\bigwedge^3 W^{\ast} \otimes\OO_X (2\KX)\ar[r]& \cdots
}
$$
 From the spectral sequence associated to this complex one deduces the isomorphism 
\begin{equation}\label{d3}
d_3: E^{0,2}_3 \longrightarrow E^{3,0}_3,
\end{equation}
where the term  $E^{0,2}_3 =E^{0,2}_2 =E^{0,2}_1$ is the kernel of the map in  \eqref{m2dual-Koszul}, while the term on the right is the Koszul cohomology group alluded to above.
More precisely, the group $E^{3,0}_3=E^{3,0}_2 =E^{3,0}_1$ is the quotient of the space of Koszul cocycles, i.e. linear maps
$$
\gamma:\mbox{$ \bigwedge^3 W \longrightarrow H^0 (\OO_X (2\KX))$}
$$
subject to the cocycle relation
$$
\phi''' \gamma(\phi,\phi', \phi'') - \phi''\gamma(\phi,\phi', \phi''') +\phi'\gamma(\phi,\phi'', \phi''')-\phi \gamma(\phi',\phi'', \phi''')=0, \,\forall \phi,\phi', \phi'',\phi''' \in W=\HKX,
$$
 by the space of Koszul coboundaries, the cocycles $\gamma$ of the form
$\gamma=d_{Kosz} (b)$, where $b$ is a linear map
$b: \bigwedge^2 W \longrightarrow \HKX$ and $d_{Kosz}$ is the Koszul differential given by the formula
$$
d_{Kosz} (b) (\phi,\phi',\phi'')=\phi''\, b(\phi,\phi') - \phi'\, b(\phi,\phi'')+\phi\, b(\phi',\phi''),
$$
for every $\phi,\phi', \phi'' \in \HKX$.

Let us now see how to attach a Koszul cocycle as above to $\xi$ by using the extension construction. This exploits the parametrization $\alpha$ in \eqref{alpha}. Namely, for
every triple  $\phi,\phi', \phi'' \in \HKX$ we set
$\alpha^{(3)}(\phi,\phi', \phi'')$ to be the image of the exterior product
$\alpha(\phi)\wedge \alpha(\phi')\wedge \alpha(\phi'') \in \bigwedge^3 H^0 (\TT_{\xi} )$ under the natural homomorphism
$$
\mbox{$\bigwedge^3 H^0 (\TT_{\xi} ) \longrightarrow H^0 (\bigwedge^3 \TT_{\xi})= H^0 (\OO_X (2\KX)).$}
$$
This way we obtain the linear map
$$
\mbox{$\alpha^{(3)}:\bigwedge^3 \HKX  \longrightarrow H^0 (\OO_X (2\KX)).$}
$$
To see that this is a Koszul cocycle we take the second exterior power of the extension sequence \eqref{ext} tensored with $\OO_X (\KX)$
$$
\xymatrix@R=12pt@C=12pt{
0\ar[r]&\OO_X (2\KX)\ar[r]^{\wedge^2 i}& \bigwedge^2 \TT_{\xi} (\KX) \ar[r]&\Omega_X (2\KX)\ar[r]&0
}
$$
and observe the following relation of global sections of $\bigwedge^2 \TT_{\xi} (\KX)$
\begin{equation}\label{det-formula}
\phi''\alpha^{(2)}(\phi,\phi') -\phi'\alpha^{(2)}(\phi,\phi'')+ \phi\alpha^{(2)}(\phi',\phi'') = \wedge^2 i (\alpha^{(3)}(\phi,\phi', \phi'') )
\end{equation}
for every triple $\phi,\phi', \phi'' \in \HKX$ and where $\alpha^{(2)} (\psi,\psi')$ stands for the image of the exterior product
$\alpha (\psi)\wedge \alpha(\psi')$ under the natural homomorphism
$$
\mbox{$\bigwedge^2 H^0 (\TT_{\xi} ) \longrightarrow H^0 (\bigwedge^2 \TT_{\xi}).$}
$$
The above formula is essentially the standard expansion of the determinant of
a $3\times 3$ matrix with respect to one of its columns. From the point of view of Koszul cohomology, the formula exhibits the cochain $\alpha^{(3)} \in \bigwedge^3 W^{\ast}\otimes H^0 (\OO_X (2\KX)) \stackrel{\wedge^2 i}{\hookrightarrow} \bigwedge^3  W^{\ast}\otimes H^0 (\bigwedge^2\TT_{\xi} (\KX))$ as a Koszul coboundary with values in
  $H^0 (\bigwedge^2 \TT_{\xi} (\KX))$. This immediately implies that $\alpha^{(3)}$ is a cocycle (with values in $H^0 (\OO_X (2\KX))$). However, there seems to be no reason for $\alpha^{(3)}$ to be a coboundary with values in $H^0 (\OO_X (2\KX))$. This occurs precisely when $\xi$ is isotropic, i.e.
$\delta_X (\xi^2) =0$.

\begin{lem}\label{l:isotr-Koszul}
$\delta_X (\xi^2) =0$ iff the cocycle
$\alpha^{(3)} \in \bigwedge^3 W^{\ast}\otimes H^0 (\OO_X (2\KX))$ is a coboundary.
\end{lem}
\begin{pf}
This is a restatement of the fact that $d_3$ in \eqref{d3} is an isomorphism and the identification
$$
 E^{3,0}_3 = E^{3,0}_1 =\frac{ker\Big(\bigwedge^3 W^{\ast}\otimes H^0 (\OO_X (2\KX))\longrightarrow \bigwedge^4 W^{\ast}\otimes H^0 (\OO_X (3\KX))\Big)}{im\Big(\bigwedge^2 W^{\ast}\otimes H^0 (\OO_X (\KX))\longrightarrow \bigwedge^3 W^{\ast}\otimes H^0 (\OO_X (2\KX))\Big)}.
$$
\end{pf}

From now on we assume that the cohomology class $\xi$ is isotropic with respect to the quadratic form $\delta_X$. According to the above lemma we can write $\alpha^{(3)}$ as a Koszul coboundary, i.e. there is a linear map
$$
l: \mbox{$\bigwedge^2 W \longrightarrow \HKX$}
$$
such that
\begin{equation}\label{cobond-rel4}
\alpha^{(3)}(\phi,\phi', \phi'')=\phi'' l(\phi,\phi') -\phi' l(\phi,\phi'')+ 
\phi \,l(\phi',\phi''),
\end{equation}
for every triple $\phi,\phi', \phi'' \in W$. Substituting this into the formula
\eqref{det-formula} gives the following relation of global sections of $\bigwedge^2 \TT_{\xi} (\KX)$

\begin{equation}\label{det-formulaKosz}
\phi''\big(\alpha^{(2)}(\phi,\phi') -l(\phi,\phi')j\big) -\phi'\big(\alpha^{(2)} (\phi,\phi'') -l(\phi,\phi')j\big)+ \phi \,\big(\alpha^{(2)}(\phi',\phi'') -l(\phi,\phi')j\big) = 0,
\end{equation}
for every triple $\phi,\phi', \phi'' \in W$ and where $j$ is the global section of $\bigwedge^2 \TT_{\xi} (-\KX)$ corresponding to the monomorphism
$\wedge^2 i: \OO_X (\KX) \longrightarrow \bigwedge^2 \TT_{\xi}$.

The above relation means that for every triple $\phi,\phi', \phi'' \in W$ the global sections 
$\{\alpha^{(2)}(\phi,\phi') -l(\phi,\phi')j, \alpha^{(2)}(\phi,\phi'') -l(\phi,\phi')j,\alpha^{(2)}(\phi,\phi') -l(\phi,\phi')j \}$ of $\bigwedge^2 \TT_{\xi} (\KX)$ fail to generate that sheaf everywhere. This eventually will be exploited to arrive to a contradiction with the generic global generation of $\TT_{\xi}$.
The direct approach, however, fails because of the possible presence of the $1$-dimensional part of the degeneracy locus $D^1_{\xi}$ as well as the fixed part of linear systems $|Im(\alpha^{(3)}(\phi,\phi', \bullet)|$. We have seen how to get rid of those loci by replacing $\xi$ by a cohomology class $\eta' \in H^1 ({\cal P}(-\KX))$, see Proposition \ref{pro:modif-eta} for notation and precise statement.
But in this replacement the isotropy condition for $\eta'$ might be lost. Our next task is to examine such a possibility. In fact, we will be concerned with a somewhat different though closely related question. Namely, as in the case of 
$\xi$, we will attach to $\eta'$ a Koszul cocycle and investigate how far it is from being a Koszul coboundary. This will lead us to the question of deforming the cochain $l$ in \eqref{cobond-rel4}.

\section{Deforming  the Koszul cochain in \eqref{cobond-rel4}}

We recall
the parametrization
$$
\alpha'_{\eta}: \HKX \longrightarrow H^0 (\TT'_{\eta})
$$ 
introduced in Lemma \ref{l:T'eta}, 2), and the higher order products $\alpha'^{(k)}_{\eta}$, for $k=2,3$, in \eqref{hprod-alpha'}. Our study will focus on the
 triple product
\begin{equation}\label{alphaeta3}
\mbox{$\alpha'^{(3)}_{\eta}: \bigwedge^3(\HKX)\longrightarrow H^0 (\bigwedge^3\TT'_{\eta})=
 H^0 (\OO_X (2\KX -2E-E_0)))$}
\end{equation}
defined by sending the exterior product 
$\phi \wedge  \phi'\wedge \phi''$ to $\alpha'^{(3)}_{\eta} (\phi ,\phi',\phi'')$, the image of the exterior product
$\alpha'_{\eta} (\phi) \wedge \alpha'_{\eta} (\phi')\wedge \alpha'_{\eta} (\phi'')$
under the natural homomorphism
$$
\mbox{$\bigwedge^3(H^0 (\TT'_{\eta}))\longrightarrow H^0 (\bigwedge^3\TT'_{\eta})=
 H^0 (\OO_X (2\KX -E')),$}
$$
where $E':=2E +E_0$. 

As in the case of $\alpha^{(3)}$, the map $\alpha'^{(3)}_{\eta}$ is a Koszul cocycle
due to the determinantal formula
$$
\alpha'^{(3)}_{\eta} (\phi ,\phi',\phi'')\cdot h=\phi'' \alpha'^{(2)}_{\eta} (\phi ,\phi') - \phi' \alpha'^{(2)}_{\eta} (\phi ,\phi'')+\phi \alpha'^{(2)}_{\eta} (\phi' ,\phi''),
$$
where $h$ stands for the global section of $\bigwedge^2 \TT'_{\eta} (-(\KX-E'))$
corresponding to the monomorphism in the exact sequence
$$
\xymatrix@R=12pt@C=12pt{
0\ar[r]&\OO_X (\KX-E')\ar[r]& \bigwedge^2 \TT'_{\eta}  \ar[r]&{\cal P}(\KX)\ar[r]&0.
}
$$
This sequence is the second exterior power of the extension sequence \eqref{ext-eta'} defined by $\eta'$, see Proposition \ref{pro:modif-eta}, 4).

We wish to understand how far the Koszul cocycle $\alpha'^{(3)}_{\eta}$ is from being a Koszul coboundary.
%
%
%
%
  Our considerations begin by comparing the cocycles $\alpha^{(3)}$ and $\alpha'^{(3)}_{\eta}$. This is given by the formula
$$
\alpha^{(3)} =\tau^2_E \alpha^{(3)}_{\eta} =\tau^2_E e \alpha'^{(3)}_{\eta},
$$
where the first and the second  equalities are \eqref{alpha-alphaeta} and \eqref{rela3-a'3} respectively and $e$ is a global section of $\OO_X (E_1)$ defining $E_0$. Denoting by $e'=\tau^2_E e$ a global section of $\OO_X (2E+E_1)$ corresponding to the divisor $E'(=2E+E_0)$ in $|2E+E_1|$, we record the above identity as follows:
\begin{equation}\label{alpha-alpheta'} 
\alpha^{(3)}=e' \alpha'^{(3)}_{\eta}.
\end{equation}

We know that the condition of isotropy of $\xi$ means that $\alpha^{(3)}$ is a Koszul coboundary
\begin{equation}\label{alpha-alpheta'-Kosz}
d_{Kosz} (l)=\alpha^{(3)} =e' \alpha'^{(3)}_{\eta} ,
\end{equation}
where $l: \bigwedge^2 \HKX \longrightarrow \HKX$ is a linear map. This linear  
map is not unique and the ambiguity is given by $d_{Kosz} (f)$, for any linear form $f:\HKX \longrightarrow H^0 (\OO_X )$. This is important since we are able to deform $l$ to take into account the divisor $E'$. Throughout the subsequent discussion we assume that $E'\neq 0$, otherwise $\alpha^{(3)} = \alpha'^{(3)}_{\eta}$ and there is nothing to modify.
\begin{lem}\label{l:l-modif}
 If the canonical image of $X$ has no rational
normal curves $\Gamma$ of degree $d_{\Gamma} \leq (p_g -1)$, then there exists $f:\HKX \longrightarrow H^0 (\OO_X ) \cong \CC$ such that
\begin{equation}\label{l-modif}
l=e'_1 \,m +d_{Kosz} (f),
\end{equation}
where $(e'_1=0)=E'_1 $ is a component of $E'$ such that the reduced divisor $(E'_1)_{red}$ equals $E'_{red}$, the reduced part of $E'$, and $m:\bigwedge^2 \HKX \longrightarrow H^0 (\OO_X (\KX -E'_1))$ is a linear map. Furthermore, the image $Im(m)$ of $m$ corresponds to the linear subsystem $|Im(m)|$ of $|\KX -E'_1|$ with at most $0$-dimensional base locus.

In addition, the cocycle $\alpha'^{(3)}_{\eta}$ is subject to the following coboundary condition
 \begin{equation}\label{Kosz-coboune'0}
e'_0 \alpha'^{(3)}_{\eta}=d_{Kosz}(m)
\end{equation}
where $e'_0$ is defined by the relation $e'=e'_1 e'_0$.
\end{lem}

\begin{pf}
We start with the Koszul relation in \eqref{alpha-alpheta'-Kosz}
\begin{equation}\label{Koszrel-l}
d_{Kosz} (l)(\phi,\phi',\phi'')=\phi'' l(\phi,\phi') -\phi' l(\phi,\phi'')+ \phi \,l(\phi',\phi'') = \alpha^{(3)} (\phi,\phi',\phi'') =e'\alpha'^{(3)}_{\eta} (\phi,\phi',\phi''),
\end{equation}
 holding for any triple $\phi,\phi',\phi'' \in \HKX$. 
To define our deformation of $l$ we 
 examine the restriction of \eqref {Koszrel-l} to the complete intersection of the divisor ${E'_{red}}$ with the divisors of the canonical linear system $|\KX|$. Namely, for every $\phi \in \HKX$ with the divisor $C_{\phi}=(\phi=0)$ having no components in common with ${E'_{red}}$, we restrict the above relation to the complete intersection subcheme $\Delta_{E'_{red}, \phi} =C_{\phi} \cdot E'_{red} $ to obtain
\begin{equation}\label{restCphi-Koszul0} 
[\phi'' \,\,l(\phi,\phi') -\phi'\, \,l(\phi,\phi'')]|_{\Delta_{E'_{red}, \phi}} =0.
\end{equation}
Write $E'_{red}=\sum^N_{k=1}  \Gamma_k$, the decomposition of $E'_{red}$ into its irreducible components $\Gamma_k$'s. 
Then 
$$
\Delta_{E'_{red}, \phi} =C_{\phi} \cdot E'_{red} =\sum^N_{k=1} C_{\phi} \cdot \Gamma_k.
$$
 We choose $C_{\phi}$ smooth so that all complete intersections $C_{\phi} \cdot \Gamma_k$
are distinct points in the smooth part of $E'_{red} =\sum^N_{k=1}   \Gamma_k$, and in general position in each linear span $\Pi_k$ of
$\Gamma_k$ in $\PP(\HKX^{\ast})$. 

Let  $\Pi_{E'_{red}}$ be the linear span of $E'_{red}$. The hyperplane $H_{\phi}$ of $\PP(\HKX^{\ast})$ corresponding to $\phi$ intersects $\Pi_{E'_{red}}$ along the hyperplane $\Pi_{E'_{red}, \phi}$ spanned by the $0$-dimensional subscheme
 $\Delta_{E'_{red}, \phi} \subset \Pi_{E'_{red}, \phi} $. Observe:
\begin{equation}\label{degComplint}
deg(\Delta_{E'_{red}, \phi})=\sum^N_{k=1} deg( \Gamma_k) > \sum^N_{k=1} dim(\Pi_{k}) \geq dim (\Pi_{E'_{red}}),
\end{equation}
where the first inequality uses the condition that none of the curves $\Gamma_k$ can be a rational normal curve in its linear span $\Pi_{k}$. The resulting inequality
\begin{equation}\label{degComplint1}
deg(\Delta_{E'_{red}, \phi}) >dim (\Pi_{E'_{red}})
\end{equation}
will be used shortly in the proof of Claim \ref{cl:endom-id} below. 

The projective space $\Pi_{E'_{red},\phi}$ is the projectivization of the vector space $\Lambda$, whose dual
$\Lambda^{\ast}$ is defined as follows:
$$
\Lambda^{\ast} =\HKX/(e'_r H^0(\OO_X (\KX-E'_{red})) +\CC\phi),
$$
where $(e'_r =0)=E'_{red}$.

From the equation \eqref{restCphi-Koszul0} it follows that $l(\phi,\bullet)$ defines an endomorphism
$$
l_{\phi}: \Lambda^{\ast} \longrightarrow  \Lambda^{\ast},
$$
where $l_{\phi} (v) =l(\phi,\widetilde{v})\,\, mod \Big(e'_rH^0(\OO_X (\KX-E'_{red})) +\CC\phi \Big)$, for every $v\in \Lambda^{\ast}$ and where $\widetilde{v}$ is any lift of $v$ to
$\HKX$. 
 The essential observation is the following.
\begin{cl}\label{cl:endom-id}
The endomorphism $l_{\phi}$ is a multiple of identity.
\end{cl}
 
Let us assume this and complete the proof of the lemma. According to Claim \ref{cl:endom-id} we have $l_{\phi} =\lambda id_{\Lambda^{\ast}}$ for some constant $\lambda$. Set $f(\phi)=-\lambda$. Then we have
$$
l(\phi,\psi) +f(\phi)\psi \in e'_rH^0(\OO_X (\KX-E'_{red})) +\CC\phi, \,\forall \psi \in \HKX.
$$
This implies that there is a unique linear function $f: \HKX \longrightarrow \CC$ and a unique linear map $m(\phi,\psi): \HKX \longrightarrow H^0(\OO_X (\KX-E'_{red})$ such that  
$$
l(\phi,\psi) +f(\phi)\psi = e'_rm(\phi,\psi) + f(\psi)\phi,  \,\forall \psi \in \HKX.
$$
 Substituting these expressions into \eqref{Koszrel-l} we deduce 
$$
\phi [l(\phi',\phi'')-f(\phi'') \phi' +f (\phi') \phi''] \in e'_r H^0 (\OO_X (\KX-E'_{red})),  \,\forall \phi',\phi'' \in \HKX.
$$
Since $\phi$ does not vanish on any irreducible component of $E'_{red}$, we deduce
$$
l(\phi',\phi'')-f(\phi'') \phi' +f (\phi') \phi'' \in e'_r H^0 (\OO_X (\KX-E'_{red})),  \,\forall \phi',\phi'' \in \HKX.
$$
Thus there is a unique linear map
$$
\mbox{$m_{E'_{red}}: \bigwedge^2 \HKX \longrightarrow H^0(\OO_X (\KX-E'_{red}))$}
$$
such that
$$
l(\phi',\phi'') = f(\phi'') \phi' -f (\phi') \phi'' +e'_r m_{E'_{red}}(\phi',\phi'') =d_{Kosz} (f)(\phi',\phi'') + e'_r m_{E'_{red}}(\phi',\phi'') ,  \,\forall \phi',\phi'' \in \HKX.
$$
Hence the coboundary $l$ has the form
\begin{equation}\label{lmE'red}
l=d_{Kosz} (f) +e'_r m_{E'_{red}}.
\end{equation}
Substituting  into \eqref{alpha-alpheta'-Kosz} gives the relation
\begin{equation}\label{lmE'red1}
e' \alpha'^{(3)}_{\eta} = d_{Kosz} (l) = d_{Kosz} (d_{Kosz} (f) +e'_r m_{E'_{red}})=e'_r d_{Kosz} (m_{E'_{red}}).
\end{equation}

 Let $B=(e'_B =0)$ be the fixed part of the linear subsystem $|Im (m_{E'_{red}})| \subset |\KX-E'_{red}|$. Then the coboundary $m_{E'_{red}} $ has the form
$$
m_{E'_{red}} =e'_B m,
$$
where $m:\bigwedge^2 \HKX \longrightarrow H^0 (\OO_X(\KX -E'_{red}-B))$ determines
 the linear subsystem $|Im(m)| \subset |\KX-E'_{red}-B| $ which is fixed part free.
Thus the Koszul relation \eqref{lmE'red1} becomes
$$
e' \alpha'^{(3)}_{\eta}=e'_r d_{Kosz} (m_{E'_{red}})=e'_r d_{Kosz} (e'_B m)=e'_r e'_B d_{Kosz} ( m).
$$
Since the linear subsystem $|Im( \alpha'^{(3)}_{\eta})|$ is fixed part free, see Lemma \ref{lem:fixed}, the equation above implies that the divisor $E'_{red} +B=(e'_r e'_B=0)$ is a component of $E'=(e'=0)$. Thus setting $e'_1:=e'_r e'_B$, the cochain $l$ in \eqref{lmE'red} becomes
$$
l=d_{Kosz} (f) +e'_r m_{E'_{red}} =d_{Kosz} (f) +e'_r e'_B m=d_{Kosz} (f) +e'_1 m.
$$
Hence the equality
$
e'\alpha'^{(3)}_{\eta}=e'_1 d_{Kosz} ( m)
$.
Dividing by $e'_1$, gives the relation
$$
e'_0\alpha'^{(3)}_{\eta}=d_{Kosz} ( m),
$$
where $e'_0$ is defined by the factorization $e'=e'_1 e'_0$.  

All statements of the lemma are now proved. It remains
 to prove Claim \ref{cl:endom-id}. For this we return to the complete intersection cycle $\Delta_{E'_{red}, \phi}$. By definition this is a set of $deg (\Delta_{E'_{red}, \phi})$ distinct points
spanning the projective space 
$\Pi_{E'_{red},\phi} =\PP(\Lambda)$. Set $\pi=dim (\Lambda)$ and choose the points
$p_1, \ldots, p_{\pi} $ in $\Delta_{E'_{red}, \phi}$ so that they span $\Pi_{E'_{red},\phi}$.
In addition, we may always choose the set $\{p_1, \ldots, p_{\pi} \} $ so that the hyperplanes in $\Pi_{E'_{red},\phi}$ spanned by any of its $(\pi-1)$ points  do not contain any other points of $\Delta_{E'_{red}, \phi}$. With such a choice made, we let $\widetilde{p_1}, \ldots, \widetilde{p_{\pi}} $ be vectors in $\Lambda$ overlying the points $p_1, \ldots, p_{\pi} $. The vectors $\{\widetilde{p_1}, \ldots, \widetilde{p_{\pi}} \}$ form a basis of $\Lambda$
and we let $\{v_1, \ldots, v_{\pi} \}$
to be the dual basis of $\Lambda^{\ast}$. We will now compute
$l_{\phi} (v_i)$, for $i=1,\ldots,\pi$.

By definition $l_{\phi} (v_i)=l(\phi,\phi_i)\,mod \Big(e'_rH^0(\OO_X (\KX-E'_{red})) +\CC\phi \Big)$, where $\phi_i \in \HKX$ is a lifting of $v_i$. Observe, by the choice of $\{p_1, \ldots, p_{\pi} \} $ made, the hyperplane in $\PP(\Lambda)$ defined by $\phi_i$ contains the points $\{p_j\}_{j\neq i}$ and no other point of $\Delta_{E'_{red}, \phi}$, i.e.,  we have
\begin{equation}\label{choice}
\mbox{$\phi_i(p_j)=0$, for all $j\neq i$, and $\phi_i (p)\neq 0$, for all $p\in \Delta_{E'_{red}, \phi} \setminus \{p_j\}_{j\neq i}$.}
\end{equation}
Substituting $\phi_i $ for $\phi'$ in the equation \eqref{restCphi-Koszul0} gives 
\begin{equation}\label{choice1}
[\phi'' l(\phi,\phi_i) -\phi_i l(\phi,\phi'')]|_{\Delta_{E'_{red}, \phi}} =0, \,\,\forall \phi'' \in \HKX.
\end{equation}
From this it follows that $l(\phi,\phi_i)$ vanishes on the set $\{p_j\}_{j\neq i}$ and hence 
$$
l(\phi,\phi_i) \equiv \lambda_i \phi_i \,mod \Big(e'_rH^0(\OO_X (\KX-E'_{red})) +\CC\phi \Big) =\lambda_i v_i,
$$
 for some scalar $\lambda_i \in \CC$. Thus
$l_{\phi} (v_i)=\lambda_i v_i$, for every $i$.
Furthermore, substituting $\phi_j$ for $\phi''$ in \eqref{choice1}, we obtain
$$
0=[\phi_j l(\phi,\phi_i) -\phi_i l(\phi,\phi_j)]|_{\Delta_{E'_{red}, \phi}}
= [\phi_j (\lambda_i \phi_i) -\phi_i (\lambda_j \phi_j)]|_{\Delta_{E'_{red}, \phi}}=(\lambda_i-\lambda_j) [\phi_i \phi_j]_{\Delta_{E'_{red}, \phi}}, \forall i,j.
$$
From \eqref{choice} it follows that $\phi_i \phi_j$ does not vanish at any point
of the set
$\Delta_{E'_{red}, \phi} \setminus \{p_k\}_{k=1,\ldots,\pi}$ and that set, according to
the inequality \eqref{degComplint1}, is nonempty. Hence $\lambda_i=\lambda_j$,
for all $i,j$ and thus $l_{\phi}$ is a multiple of identity as asserted in Claim \ref{cl:endom-id}. 
\end{pf}

We will now draw some geometric conclusions from the coboundary relation 
\begin{equation}\label{cobound-rel}
e'_0 \alpha'^{(3)}_{\eta} =d_{Kosz} (m)
\end{equation}
in Lemma \ref{l:l-modif}, \eqref{Kosz-coboune'0}. They all will result from the considerations of the
cochain
\begin{equation}\label{m}
\mbox{$m:\bigwedge^2 \HKX \longrightarrow H^0 (\OO_X (\KX-E'_1)).$}
\end{equation}
Namely, for every nonzero $\phi \in \HKX$ we consider the linear map
$$
m(\phi,\bullet): \HKX \longrightarrow H^0 (\OO_X (\KX-E'_1)).
$$
\begin{lem}\label{l:m(phi,.)}
Set $W_{\phi}:=ker(m(\phi,\bullet))$. Then the following holds.

1) $\CC\phi \subset W_{\phi}$ and $dim (W_{\phi}) \geq 3$,

2) for $\phi$ with $C_{\phi} =(\phi=0)$ having no components  in common with $E'_0 =(e'_0)$, the restriction of $m$ to the subspace
$\bigwedge^2 W_{\phi}$ takes values in the subspace $e'_0 H^0 (\OO_X (\KX-E')$
of $H^0 (\OO_X (\KX-E'_1))$. In particular, 
\begin{equation}\label{Wphi}
 h^0 (\OO_X (\KX-E')) \geq 2dim(W_{\phi})-5 \geq 1,
\end{equation}
for every $\phi$ with $C_{\phi} =(\phi=0)$ irreducible. 
\end{lem}
\begin{pf}
The inclusion $\CC\phi \subset W_{\phi}$ is obvious in view of the skew-symmetry of $m$. To see the asserted lower bound on the dimension of $W_{\phi}$, we observe that $h^0 (\OO_X (\KX-E'_1)) \leq p_g -3$, since by assumption \eqref{cond}, $(iii)$, the divisor $E'_1$ can not be a line. From this it follows
$$
dim(W_{\phi})=p_g -dim (Im(m(\phi,\bullet))) \geq p_g -h^0 (\OO_X (\KX-E'_1))\geq p_g -(p_g -3)=3.
$$

For part 2), we use the coboundary relation \eqref{cobound-rel} which gives
$$
e'_0 \alpha'^{(3)}_{\eta} (\phi,\phi', \phi'')=\phi'' m(\phi, \phi')-\phi' m(\phi, \phi'')+\phi m(\phi', \phi'')
$$
for every triple $\phi,\phi', \phi''\in \HKX$. Taking $\phi', \phi''\in W_{\phi}$
in the above equation, we obtain
\begin{equation}\label{Wphi1}
e'_0 \alpha'^{(3)}_{\eta} (\phi,\phi', \phi'')=\phi m(\phi', \phi'').
\end{equation}
If, in addition, the divisor $C_{\phi} =(\phi=0)$ has no common components with $E'_0=(e'_0=0)$, the above equation implies that $m(\phi', \phi'')$ is a multiple of $e'_0$, i.e. there is a unique global section $\sigma(\phi', \phi'') \in H^0 (\OO_X (\KX-E'))$ such that $m(\phi', \phi'')=e'_0 \sigma(\phi', \phi'') $, for every
$\phi', \phi'' \in W_{\phi}$. Thus we obtain a linear map
$$
\mbox{$ \sigma_{\phi}: \bigwedge^2 (W_{\phi} /\CC\phi) \longrightarrow H^0 (\OO_X (\KX-E'))$}
$$
subject to the relation
\begin{equation}\label{Wphi2}
 \alpha'^{(3)}_{\eta} (\phi,\phi', \phi'')=\phi  \sigma_{\phi}(\phi', \phi''), \,\,\forall \phi', \phi'' \in W_{\phi} .
\end{equation}

We now assume that $C_{\phi} $ is irreducible. Then according to Proposition \ref{pro:modif-eta}, 5), $(ii)$, the left side of the equation \eqref{Wphi2} is nonzero for any linearly independent triple $\{\phi,\phi', \phi''\}$. Hence $\sigma_{\phi}$ does not vanish on any nonzero decomposable vector $\psi \wedge \psi' \in \bigwedge^2 (W_{\phi} /\CC\phi)$. Hence the estimate
$$
h^0 (\OO_X (\KX-E')) \geq dim(Im(\sigma_{\phi})) \geq 2\Big(dim(W_{\phi} /\CC\phi) -2\Big)+1 = 2dim(W_{\phi})-5 \geq 1,
$$
where the last inequality comes from 1) of the lemma.
\end{pf}

We can now strengthen the assertion in Lemma \ref{l:l-modif} about the linear system $|Im(m)|$ having at most $0$-dimensional base locus.  
\begin{cor}\label{cor:m-3dim}
Assume $p_g \geq 5$. Then, for a general $3$-dimensional subspace $\Pi \subset \HKX$, the linear system
corresponding to the subspace
$$
V_{\Pi}:= \{ m(\psi, \psi') | \psi, \psi' \in \Pi \} \subset Im(m) 
$$
has at most $0$-dimensional base locus.
\end{cor}
\begin{pf}
Assume that every nonzero subspace $V_{\Pi}$ gives the linear system $|V_{\Pi}|$ with $1$-dimensional base locus. Set $B_{\Pi}$ to be the divisorial part of that base locus. Then the divisors $B_{\Pi}$ must move as $\Pi$ varies, since otherwise
$|Im(m)|$ has $1$-dimensional base locus and this is contrary to Lemma \ref{l:l-modif}. From the relation $e'_0 \alpha'^{(3)}_{\eta} =d_{Kosz} (m)$ in \eqref{Kosz-coboune'0}
 it follows that for a general $\Pi$ the moving part of $B_{\Pi}$, call it $B'_{\Pi}$, must be a component of the divisor  
$(\alpha'^{(3)}_{\eta} (\phi, \phi',\phi'')=0)$, where $\{\phi, \phi',\phi''\}$ is a basis of $\Pi$. By Proposition \ref{pro:F'phi-gen-nopencil}, 2), the divisor
$(\alpha'^{(3)}_{\eta} (\phi, \phi',\phi'')=0)$ is irreducible, for a general choice of $\Pi$. Hence, the equality
$$
 B'_{\Pi} =(\alpha'^{(3)}_{\eta} (\phi, \phi',\phi'')=0) =2\KX -E'.
$$
On the other hand, by definition, $B'_{\Pi} =\KX-E'_1 -M_{\Pi}$, where $M_{\Pi}$ is the moving part of the linear system $|V_{\Pi}|$. Substituting into the above equality we obtain
\begin{equation}\label{K-E'-M}
\KX-E'_1 -M_{\Pi} =2\KX -E'
\end{equation}
or, equivalently, $\KX -E' =-E'_1 -M_{\Pi}$. But by Lemma \ref{l:m(phi,.)}, 2), the divisor $(\KX -E')$ is effective. Hence $E'_1 =M_{\Pi}=0$. Since the reduced parts of $E'_1 $ and $E'$ are the same, Lemma \ref{l:l-modif}, we deduce that $E'=0$ and this is contrary to our assumption (of course, with the divisors $E'=E'_1=M_{\Pi} =0$, the equation \eqref{K-E'-M} becomes $\KX =2\KX$, i.e., $\KX =0$ which is impossible).
\end{pf}

We now return to subspaces $W_{\phi}$ in Lemma \ref{l:m(phi,.)} and show how to `lift' it to the level of the category of coherent sheaves on $X$. More precisely, we take $\phi \in \HKX$ with
$C_{\phi} $ a smooth curve and $\alpha'_{\eta} (\phi)$ with no zeros and consider the morphism
\begin{equation}\label{beta_phi}
\beta_{\phi}: \Big (W_{\phi}/ \CC\phi \Big) \otimes \OO_X \longrightarrow \FF'_{[\phi]}
\end{equation}
which is the composition of the inclusion $\Big (W_{\phi}/ \CC\phi\Big) \otimes \OO_X \hookrightarrow H^0( \FF'_{[\phi]}) \otimes \OO_X$ together with the evaluation morphism
$H^0( \FF'_{[\phi]}) \otimes \OO_X \longrightarrow \FF'_{[\phi]}$. The inclusion comes from the parametrization
\begin{equation}\label{paramF'phi}
 \HKX /\CC\phi \longrightarrow H^0(\FF'_{[\phi]})
\end{equation}
established in Proposition \ref{pro:F'phi}, 2). 

\begin{lem}\label{l:coker-beta_phi}
The cokernel of $\beta_{\phi}$ in \eqref{beta_phi} is supported on $C_{\phi}$.
In particular, the restriction $\FF'_{[\phi]} \otimes \OO_{C_{\phi}}$ fits into the following exact sequence
\begin{equation}\label{F'phionCphi-1}
\xymatrix@R=12pt@C=12pt{
0\ar[r]& \OO_{C_{\phi}} (\KX) \otimes \OO_{C_{\phi}} (-B_{\phi}) \ar[r]&{\cal F'}_{[\phi]}\otimes \OO_{C_{\phi}}\ar[r]&\OO_{C_{\phi}} (\KX -E') \otimes \OO_{C_{\phi}} (B_{\phi})\ar[r]&0,}
\end{equation}
where $B_{\phi}$ is an effective divisor on $C_{\phi}$. Furthermore, the morphism 
\begin{equation}\label{F'phionCphi-2}
\overline{\beta_{\phi}}: \Big(W_{\phi}/ \CC\phi \Big)\otimes \OO_X \longrightarrow {\cal F'}_{[\phi]}\otimes \OO_{C_{\phi}},
\end{equation} 
the composition of $\beta_{\phi}$ with the restriction morphism ${\cal F'}_{[\phi]}\longrightarrow {\cal F'}_{[\phi]}\otimes \OO_{C_{\phi}}$, factors through $\OO_{C_{\phi}} (\KX) \otimes \OO_{C_{\phi}} (-B_{\phi})$, i.e. there is a commutative diagram
\begin{equation}\label{F'phionCphi-3}
\xymatrix@R=12pt@C=12pt{ 
&& \Big(W_{\phi}/ \CC\phi \Big) \otimes \OO_X \ar[d]^{\overline{\beta_{\phi}}} \ar[ld]&&\\
0\ar[r]& \OO_{C_{\phi}} (\KX) \otimes \OO_{C_{\phi}} (-B_{\phi}) \ar[r]&{\cal F'}_{[\phi]}\otimes \OO_{C_{\phi}}\ar[r]&\OO_{C_{\phi}} (\KX -E') \otimes \OO_{C_{\phi}} (B_{\phi})\ar[r]&0.
}
\end{equation} 
In particular, $B_{\phi}$ is contained in the base locus of the linear subsystem
$|W_{\phi}/ \CC\phi  | \subset |\OO_{C_{\phi}} (\KX)|$.
\end{lem}
\begin{pf}
For $\psi \in \HKX$ denote by $f_{\psi}$ the global section of $\FF'_{[\phi]}$ corresponding to the equivalence class $\overline{\psi} \in \HKX /\CC\phi$ under the isomorphism in \eqref{paramF'phi}.
From Proposition \ref{pro:F'phi}, 3), we know that 
$\alpha'^{(3)}_{\eta} (\phi,\phi', \phi'') =f_{\phi'} \wedge f_{\phi''} $, for all $\phi',\phi'' \in \HKX$, where the expression on the right stands for the image of
the exterior product $f_{\phi'} \wedge f_{\phi''}$ under the map
$$
 \mbox{$\bigwedge^2 H^0(\FF'_{\phi})\longrightarrow H^0(\bigwedge^2 \FF'_{\phi})=H^0(\OO_X (2\KX-E')).$}
$$
This together with the formula \eqref{Wphi2} imply that all global sections in the image of $\beta_{\phi}$ are proportional along $C_{\phi}$. Hence the cokernel of the morphism $\overline{\beta_{\phi}}$ defined in \eqref{F'phionCphi-2} is a sheaf of rank $1$ on $C_{\phi}$. Factoring out its torsion part (if nonzero), we obtain a torsion free and, hence, locally free sheaf on $C_{\phi}$. Call that sheaf ${\cal M'}$. By construction we have the exact sequence
$$
\xymatrix@R=12pt@C=12pt{
0\ar[r]&{\cal M} \ar[r]&{\cal F'}_{[\phi]}\otimes \OO_{C_{\phi}}\ar[r]&{\cal M'}\ar[r]&0,}
$$
where ${\cal M}= \OO_{C_{\phi}} (2\KX-E')\otimes {\cal M'}^{-1}$, and the factorization of $\overline{\beta_{\phi}}$ through ${\cal M}$ as in the diagram \eqref{F'phionCphi-3}. It remains to show that ${\cal M}$ (resp. ${\cal M'}$) has the form asserted in \eqref{F'phionCphi-1}. For this we recall the exact sequence
$$
\xymatrix@R=12pt@C=12pt{
0\ar[r]&\OO_{C_{\phi}} (\KX-E') \ar[r]&{\cal F'}_{[\phi]}\otimes \OO_{C_{\phi}}\ar[r]&{\cal O}_{C_{\phi}} (\KX) \ar[r]&0}
$$
 in Corollary \ref{Fphi-Cphi},2). Putting this sequence together with the preceding one, we obtain the diagram
\begin{equation}\label{F'phionCphi-4}
\xymatrix@R=12pt@C=12pt{
&&0\ar[d]&&\\
& &{\cal M} \ar[d] \ar[dr]&&\\
0\ar[r]&\OO_{C_{\phi}} (\KX-E') \ar[r]&{\cal F'}_{[\phi]}\otimes \OO_{C_{\phi}}\ar[r]\ar[d]&{\cal O}_{C_{\phi}} (\KX) \ar[r]&0\\
& &{\cal M}' \ar[d]&&\\
&&0&&}
\end{equation}
The slanted arrow, obtained as the composition of the monomorphism of the column with the epimorphism of the row, is nonzero since on the level of global sections 
$$
\xymatrix@R=12pt@C=12pt{
H^0 ( {\cal F'}_{[\phi]})\ar[d] \ar[dr]&\\
H^0( {\cal F'}_{[\phi]}\otimes \OO_{C_{\phi}} ) \ar[r]& H^0 ({\cal O}_{C_{\phi}} (\KX)) \cong \HKX /\CC\phi
}
$$
the composition (slanted) arrow is the isomorphism in Proposition \ref{pro:F'phi}, 2).

Once the morphism ${\cal M}\longrightarrow {\cal O}_{C_{\phi}} (\KX)$ in the diagram \eqref{F'phionCphi-4}  is nonzero, we identify it as a nonzero global section, call it $\rho_{\phi}$, of 
${\cal O}_{C_{\phi}} (\KX)\otimes {\cal M}^{-1} =\OO_{C_{\phi}}(B_{\phi}) $, where $B_{\phi}=(\rho_{\phi}=0)$. Hence ${\cal M} ={\cal O}_{C_{\phi}} (\KX)\otimes \OO_{C_{\phi}} (-B_{\phi})$ and ${\cal M}'=\det({\cal F'}_{[\phi]}\otimes \OO_{C_{\phi}}  )\otimes{\cal M}^{-1} =\OO_{C_{\phi}}(2\KX-E')\otimes{\cal M}^{-1}=\OO_{C_{\phi}}(\KX-E')\otimes
\OO_{C_{\phi}} (B_{\phi})$ as asserted in \eqref{F'phionCphi-1}.
\end{pf} 

Using the epimorphism in \eqref{F'phionCphi-1} we can now define a modification of ${\cal F'}_{[\phi]}$ along the curve $C_{\phi}$. Namely, taking the composition
\begin{equation}\label{modifF'phi}
{\cal F'}_{[\phi]} \longrightarrow {\cal O}_{C_{\phi}} (\KX-E')\otimes \OO_{C_{\phi}} (B_{\phi})
\end{equation}
of the restriction morphism ${\cal F'}_{[\phi]} \longrightarrow{\cal F'}_{[\phi]}\otimes {\cal O}_{C_{\phi}} $ with the epimorphism in \eqref{F'phionCphi-1}, we define
\begin{equation}\label{modifF'phi1}
{\cal F''}_{[\phi]}:=ker \Big({\cal F'}_{[\phi]} \longrightarrow {\cal O}_{C_{\phi}} (\KX-E')\otimes \OO_{C_{\phi}} (B_{\phi})\Big).
\end{equation}
The following properties of ${\cal F''}_{[\phi]}$ follow immediately from the construction.
\begin{pro}\label{pro:F''phi}
1) $\FF''_{[\phi]}$ is a locally free sheaf of rank $2$ fitting into the following exact sequence
\begin{equation}\label{modifF'phi2}
\xymatrix@R=12pt@C=12pt{
0\ar[r]& \FF''_{[\phi]}\ar[r]&{\cal F'}_{[\phi]}\ar[r]&{\cal O}_{C_{\phi}} (\KX-E')\otimes \OO_{C_{\phi}} (B_{\phi})\ar[r]&0.
}
\end{equation}
In particular, $c_1 (\FF''_{[\phi]}) =\KX-E'$ and $c_2 (\FF''_{[\phi]})=c_2(\FF'_{[\phi]})-\KX^2 +deg(B_{\phi})$.

\vspace{0.2cm}
2) The morphism 
$\beta_{\phi}: \Big(W_{\phi} /\CC\phi\Big) \otimes \OO_X \longrightarrow\FF'_{[\phi]}$ in \eqref{beta_phi} factors through $\FF''_{[\phi]}$ and gives an injection
$$
 W_{\phi} /\CC\phi  \hookrightarrow H^0 (\FF''_{[\phi]} ).
$$  
\end{pro}

\vspace{0.2cm}
We can go further and `promote' $W_{\phi}$ to an object in the category of (short) exact complexes. For this we recall that ${\cal F'}_{[\phi]}$ is defined as a quotient of $\TT'_{\eta}$ in the exact sequence
$$
\xymatrix@R=12pt@C=12pt{
0\ar[r]& \OO_X \ar[r]^{\alpha'_{\eta} (\phi)}&{\cal T}'_{\eta}\ar[r]&\FF'_{[\phi]}\ar[r]&0,
}
$$
see \eqref{shF'phi}. Combining this with the morphism in \eqref{modifF'phi}
gives the epimorphism
 \begin{equation}\label{modifT'eta}
{\cal T}'_{\eta} \longrightarrow {\cal O}_{C_{\phi}} (\KX-E')\otimes \OO_{C_{\phi}} (B_{\phi})
\end{equation}
and we set ${\cal T}'_{\eta, \phi}$ to be the kernel of that epimorphism. This gives rise to the following commutative diagram.
 \begin{equation}\label{diag:modifT'eta}
\xymatrix@R=12pt@C=12pt{
&&0\ar[d]&0\ar[d]&\\
0\ar[r]& \OO_X \ar[r]\ar@{=}[d]&{\cal T}'_{\eta,\phi}\ar[r]\ar[d]&\FF''_{[\phi]}\ar[r]\ar[d]&0\\
0\ar[r]& \OO_X \ar[r]^{\alpha'_{\eta} (\phi)}&{\cal T}'_{\eta}\ar[r]\ar[d]&\FF'_{[\phi]}\ar[r]\ar[d]&0\\
&&{\cal O}_{C_{\phi}} (\KX-E')\otimes \OO_{C_{\phi}} (B_{\phi})\ar@{=}[r]\ar[d]&{\cal O}_{C_{\phi}} (\KX-E')\otimes \OO_{C_{\phi}} (B_{\phi})\ar[d]&\\
&&0&0&
}
\end{equation}
The top row of the above diagram gives a realization of $W_{\phi}$ in the category of (short) exact complexes of locally free sheaves on $X$.

\section{From the cochain $m$ to subsheaves of $\bigwedge^2 \TT'_{\eta}$}

In this section we show how to construct geometrically interesting subsheaves of
$\bigwedge^2 \TT'_{\eta} (E'_0)$ from the fact that the Koszul cocycle $e'_0 \alpha'^{(3)}_{\eta}$ is a coboundary. For this we write down the coboundary
 condition $e'_0 \alpha'^{(3)}_{\eta}=d_{Kosz} (m)$ explicitly:
\begin{equation}\label{cobound-formula}
e'_0 \alpha'^{(3)}_{\eta} (\phi,\phi',\phi'')= \phi'' m (\phi, \phi') - \phi' m (\phi, \phi'')+\phi \,m (\phi', \phi''),\,\, \forall \phi,\phi',\phi''\in \HKX.
\end{equation}
On the other hand we have the determinantal formula
$$
\phi'' \alpha'^{(2)}_{\eta} (\phi,\phi') - \phi' \alpha'^{(2)}_{\eta} (\phi,\phi'') +\phi \alpha'^{(2)}_{\eta}(\phi',\phi'')=\alpha'^{(3)}_{\eta} (\phi,\phi',\phi'') h,
$$
where $h$ denotes the global section of $\bigwedge^2 \TT'_{\eta} (-(\KX-E'))$ corresponding to the monomorphism in the exact sequence
$$
\xymatrix@R=12pt@C=12pt{
0\ar[r]& \OO_X (K_X -E') \ar[r]& \bigwedge^2 \TT'_{\eta}  \ar[r]& {\cal P}(\KX) \ar[r]&0,}
$$
the second exterior power of the extension sequence \eqref{ext-T'} defined by $\eta'$. Combining the determinantal and the coboundary formulas together
gives the following relation
$$
e'_0 \Big(\phi'' \alpha'^{(2)}_{\eta} (\phi,\phi') - \phi' \alpha'^{(2)}_{\eta} (\phi,\phi'') +\phi \alpha'^{(2)}_{\eta}(\phi',\phi'')\Big) =\phi'' m (\phi, \phi')h - \phi' m (\phi, \phi'')h+\phi \,m (\phi', \phi'')h,
$$
for every triple $\phi, \phi', \phi'' \in \HKX$. This can be rewritten in the form
\begin{equation}\label{Koszrel-eta}
\phi'' (e'_0 \alpha'^{(2)}_{\eta} (\phi,\phi')- m (\phi, \phi')h) - \phi' (e'_0 \alpha'^{(2)}_{\eta} (\phi,\phi'')- m (\phi, \phi'')h) +\phi \,(e'_0 \alpha'^{(2)}_{\eta}(\phi',\phi'')- m (\phi', \phi'')h)=0,
\end{equation}  
for every triple $\phi, \phi', \phi'' \in \HKX$. The relation means that the triple of global sections
$\{ e'_0 \alpha'^{(2)}_{\eta} (\phi,\phi')- m (\phi, \phi')h, e'_0 \alpha'^{(2)}_{\eta} (\phi,\phi'')- m (\phi, \phi'')h, e'_0 \alpha'^{(2)}_{\eta} (\phi',\phi'')- m (\phi', \phi'')h \}$
of $\bigwedge^2 \TT'_{\eta} (E'_0)$ fail to generate this sheaf everywhere. This allows to construct subsheaves of $\bigwedge^2 \TT'_{\eta}$. Namely, let $\Pi$ be a $3$-dimensional subspace of $\HKX$. Under an identification $\bigwedge^2 \Pi \cong \Pi^{\ast}$ we have the linear map
$$
\mbox{$\Pi^{\ast}\cong \bigwedge^2 \Pi \longrightarrow H^0 (\bigwedge^2 \TT'_{\eta}(E'_0))$},
$$
defined by the rule
$\pi \wedge \pi' \mapsto e'_0 \alpha'^{(2)}_{\eta} (\pi,\pi')- m (\pi, \pi')h$, for every
$\pi \wedge \pi' \in \bigwedge^2 \Pi$. This gives rise to the morphism of sheaves
$$
\mbox{$\Pi^{\ast} \otimes \OO_X \cong \bigwedge^2 \Pi  \otimes \OO_X \longrightarrow \bigwedge^2 \TT'_{\eta} (E'_0)$}.
$$
 The Koszul relation \eqref{Koszrel-eta} tells us that the saturation 
${\cal G}_{\Pi}$ of the image of the above morphism is a subsheaf of $\bigwedge^2 \TT'_{\eta}(E'_0) $ of rank at most $2$.
\begin{lem}\label{l:rk2}
The subsheaf ${\cal G}_{\Pi}$ has rank $2$, if the subspace
$\alpha'_{\eta} (\Pi) \subset H^0 (\TT'_{\eta})$ generically generates $\TT'_{\eta}$.
\end{lem}
\begin{pf}
The result follows from the following formula
\begin{equation}\label{Pi-relation}
\begin{gathered}
(e'_0 \alpha'^{(2)}_{\eta} (\phi,\phi')- m (\phi, \phi')h) \wedge (e'_0 \alpha'^{(2)}_{\eta} (\phi,\phi'')- m (\phi, \phi'')h) =\\
e'_0 \phi \Big(m (\phi, \phi')\alpha'_{\eta}(\phi'') -m (\phi, \phi'')\alpha'_{\eta}(\phi')+m (\phi', \phi'')\alpha'_{\eta}(\phi) \Big),
\end{gathered}
\end{equation}
where the expression on the left is viewed as a global section of 
$\bigwedge^2 \Big(\bigwedge^2 \TT'_{\eta} (E'_0)\Big) \cong \TT'_{\eta} \otimes \det(\TT'_{\eta}) (2E'_0) =\TT'_{\eta} (2\KX -E'+2E'_0)=\TT'_{\eta} (2\KX -E'_1+E'_0)$.
The formula \eqref{Pi-relation} is in turn  a consequence of the coboundary condition \eqref{cobound-formula} and the identities:
$$
\begin{aligned}
\alpha'^{(2)}_{\eta} (\psi,\psi')\wedge \alpha'^{(2)}_{\eta} (\psi,\psi'')=
\alpha'^{(3)}_{\eta} (\psi,\psi',\psi''') \alpha'_{\eta} (\psi),& \hspace{0.2cm}\forall \psi,\psi',\psi''' \in \HKX, \\
\alpha'^{(2)}_{\eta} (\psi,\psi')\wedge \rho h=\rho \Big(\psi \alpha'_{\eta} (\psi') -
\psi' \alpha'_{\eta} (\psi) \Big),& \hspace{0.2cm}\forall \psi,\psi'\in \HKX , \, \forall \rho \in H^0 (\OO_X (K-E')),
\end{aligned}
$$
where the expressions on the left are viewed as global sections of
$\bigwedge^2 \Big(\bigwedge^2 \TT'_{\eta} \Big) \cong \TT'_{\eta} \otimes \det(\TT'_{\eta}) =\TT'_{\eta} (2\KX -E')$. The above identities are
 easily derived from writing down the sections $\alpha'_{\eta}(\psi)$, $\psi \in \HKX$, in a local frame of $\TT'_{\eta}$ adapted to the extension sequence \eqref{ext-T'}. 
\end{pf}

For a $3$-dimensional subspace $\Pi $ in $\HKX$ we define
\begin{equation}\label{tauPi-formula}
\tau(\phi,\phi',\phi'') =m (\phi, \phi')\alpha'_{\eta}(\phi'') -m (\phi, \phi'')\alpha'_{\eta}(\phi')+m (\phi', \phi'')\alpha'_{\eta}(\phi) 
\end{equation}
a global section of $\TT'_{\eta} (\KX-E'_1)$, for some basis $\{\phi,\phi',\phi''\}$ of
$\Pi$. This section, up to a nonzero scalar multiple, is intrinsically attached to a subspace $\Pi$. Since the subsequent constructions depend only on the homothety class of $\tau(\phi,\phi',\phi'')$, we denote that section by $\tau_{\Pi}$. With this notation, the sheaf ${\cal G}_{\Pi}$ can now be defined as the kernel of the morphism
\begin{equation}\label{tauPi-morph}
\mbox{$\bigwedge^2 \TT'_{\eta} (E'_0) \stackrel{\wedge \tau_{\Pi}}{\longrightarrow} \OO_X (3\KX-E'_1 -E').$}
\end{equation}
We establish the following property of $\tau_{\Pi}$.
\begin{lem}\label{l:zeroloc-tauPi}
The zero locus $Z_{\Pi}$ of a nonzero $\tau_{\mbox{\tiny$\Pi$}}$ is $1$-dimensional.
\end{lem} 
\begin{pf}
Assume that for a general $3$-dimensional subspace $\Pi \subset \HKX$ the zero locus $Z_{\Pi}$ of $\tau_{\Pi}$ is at most $0$-dimensional.
 Choose $\Pi \subset \HKX$ so that the corresponding linear subsystem $|\Pi|$ of $|\KX|$ is base point free. Then from the formula \eqref{Pi-relation} in the proof of Lemma \ref{l:rk2}, it follows that the evaluation morphism
\begin{equation}\label{evPi}
\xymatrix@R=12pt@C=12pt{
 ev_{\Pi^{\ast}} :\Pi^{\ast} \otimes \OO_X\ar[r]& {\cal G}_{\Pi}
}
\end{equation}
has the cokernel, call it ${\cal S}$, supported on $(e'_0\tau_{\Pi}=0)=E'_0+T$, where $T$ is the subscheme $(\tau_{\Pi}=0)$ and the determinant $\det({\cal G}_{\Pi})= \OO_X (K+E'_0)$. We claim that our choice of $\Pi$ guaranties that $T=\emptyset$. Indeed, consider the composition
$$
{\cal G}_{\Pi} \longrightarrow {\cal S} \longrightarrow {\cal S}\otimes \OO_{E'_0}.
$$
Set ${\cal G}'_{\Pi}:=ker\Big({\cal G}_{\Pi} \longrightarrow {\cal S}\otimes \OO_{E'_0} \Big)$. This is a torsion free subsheaf of ${\cal G}_{\Pi}$ of determinant
$\det ({\cal G}'_{\Pi}) =\det({\cal G}_{\Pi}) (-E'_0)= \OO_X (K)$. Furthermore, the evaluation morphism $ ev_{\Pi^{\ast}}$ in \eqref{evPi} factors through ${\cal G}'_{\Pi}$
$$
\xymatrix@R=12pt@C=12pt{
\Pi^{\ast} \otimes \OO_X\ar[r]& {\cal G}'_{\Pi}
}
$$
and the cokernel ${\cal S}_T$ of that morphism is supported on the subscheme $T$. So the above morphism can be completed to the following exact complex
$$
\xymatrix@R=12pt@C=12pt{
0\ar[r]& \OO_X (-\KX) \ar[r]& \Pi^{\ast} \otimes \OO_X\ar[r]& {\cal G}'_{\Pi}\ar[r]&{\cal S}_T \ar[r]&0.
}
$$
Dualizing and using the fact that $\Pi \otimes \OO_X \longrightarrow \OO_X (\KX)$ is surjective we deduce that the double dual $({\cal G}'_{\Pi})^{\ast\ast}$ fits  into the following exact sequence
$$
\xymatrix@R=12pt@C=12pt{
0\ar[r]& \OO_X (-\KX) \ar[r]& \Pi^{\ast} \otimes \OO_X\ar[r]& 
({\cal G}'_{\Pi})^{\ast\ast}\ar[r]&0.
}
$$ 
This implies the equality $({\cal G}'_{\Pi})^{\ast\ast} ={\cal G}'_{\Pi}$. Hence  
the cokernel ${\cal S}_T =0$ or, equivalently, $T=\emptyset$.

Once $\tau_{\Pi}$ has no zeros, the morphism in \eqref{tauPi-morph} is an epimorphism and ${\cal G}_{\Pi}$ is a part of the following exact sequence
\begin{equation}\label{GTeta}
\xymatrix@R=12pt@C=19pt{
0\ar[r]& {\cal G}_{\Pi}  \ar[r]& \bigwedge^2 \TT'_{\eta} (E'_0) \ar[r]^(0.45){ \wedge\tau_{\Pi}}& \OO_X (3\KX -2E'_1 )\ar[r]&0,}
\end{equation}
  This is a part of the
Koszul complex of $(\TT'_{\eta} (K_X -E'_1), \tau_{\Pi})$ tensored with $\OO_X (-2(K_X -E'_1)+E'_0))$. So we also have the exact sequence
\begin{equation}\label{GTeta1}
\xymatrix@R=12pt@C=19pt{
0\ar[r]& \OO_X (-2(\KX -E'_1)+E'_0)\ar[r]^(0.57){\tau_{\Pi}}& \TT'_{\eta} (-(K_X -E'))\ar[r]& {\cal G}_{\Pi}   \ar[r]&0.}
\end{equation}
The immediate consequence of this is that $E'\neq 0$, since otherwise
$H^0 (\TT_{\xi} (-K_X)) = H^0 (\TT'_{\eta} (-K_X)) \cong H^0({\cal G}_{\Pi})$ thus implying a nonzero global section of $\TT_{\xi} (-K_X)$ and hence a splitting of the extension sequence \eqref{ext}. So from now on we assume that the divisor $E'\neq 0$.

The exact sequence \eqref{GTeta1} will play an important role in the rest of the argument, since it will serve to impose serious restrictions on the divisors involved in it. We already know from Lemma \ref{l:l-modif} that the linear system $|\KX-E'_1|$ has at most $0$-dimensional base locus, so the divisor $(\KX-E'_1)$ is nef. Thus $(\KX-E'_1)^2 \geq 0$ and we proceed our considerations according to two possibilities:

-- $(\KX-E'_1)^2 > 0$,

-- $(\KX-E'_1)^2 = 0$. 

\vspace{0.2cm}
{\it Case 1: $(\KX-E'_1)^2 >0$.} This imposes
\begin{lem}\label{l:K-E'}
$h^0(\OO_X (\KX-E')) = 1$.
\end{lem}
\begin{pf}
By Lemma \ref{l:m(phi,.)} we know that $h^0(\OO_X (\KX-E'))\geq 1$. Assume $h^0(\OO_X (\KX-E')) \geq 2$. Then we write
$$
\KX-E'=M+N,
$$
where $M$ (resp. $N$) is the moving (resp. fixed) part of the linear system
$|\KX-E'_1|$. Tensoring the exact sequence \eqref{GTeta1} with $\OO_X (N)$ we obtain 
$$
\xymatrix@R=12pt@C=12pt{
0\ar[r]& \OO_X (-(\KX -E'_1) -M)\ar[r]^(0.68){\tau_{\Pi}}& \TT'_{\eta} (-M)\ar[r]& {\cal G}_{\Pi} (N)   \ar[r]&0.}
$$
Since the divisor $((\KX -E'_1) +M)$ is nef and big we deduce the vanishing
$H^1( \OO_X (-(\KX -E'_1) -M))=0$. Hence an isomorphism
$$
H^0( \TT'_{\eta} (-M))\cong H^0 ({\cal G}_{\Pi} (N)).
$$
This implies that $\TT'_{\eta} (-M)$ has a nonzero global section or, equivalently, a nonzero monomorpism
$$
\OO_X (M) \longrightarrow \TT'_{\eta}.
$$
which is in contradiction with Proposition \ref{pro:modif-eta}, 5), $(iii)$ saying that all line subsheaves of $\TT'_{\eta}$ have at most one dimensional space of global sections.
\end{pf}

The above result imposes further restrictions on the linear map
$$
\mbox{$m: \bigwedge^2 \HKX \longrightarrow H^0 (\OO_X (\KX -E'_1))$}
$$
found in Lemma \ref{l:l-modif}. Namely, recall the subspace $W_{\phi} :=ker(m(\phi,\bullet))$, see Lemma \ref{l:m(phi,.)}. We can be now  more precise about it.
\begin{lem}\label{l:K-E'}

1) For every $\phi \in \HKX$ with $C_{\phi} =(\phi=0)$ irreducible, the subspace
$W_{\phi} \subset \HKX$ is $3$-dimensional.

2) $h^0 (\OO_X (\KX -E'_1))=p_g-3$, i.e. $E'_1$ is contained in a plane.
\end{lem}
\begin{pf}
Combining the inequality \eqref{Wphi} with the equality $h^0(\OO_X (\KX-E')) = 1$ of Lemma \ref{l:K-E'} gives the estimate
$$
1=h^0 (\OO_X (\KX-E') \geq 2dim(W_{\phi}) -5
$$
from which we deduce $dim(W_{\phi}) \leq 3$. This together with the lower bound $dim (W_{\phi}) \geq 3$, see Lemma \ref{l:m(phi,.)}, 1), give
the equality $dim (W_{\phi}) = 3$.

From 1) of the lemma we deduce
$$
h^0 (\OO_X (\KX-E'_1)) \geq dim(Im(m(\phi,\bullet)) =p_g-dim (W_{\phi})=p_g-3.
$$
Since $E'_1$ can not be a line in view of \eqref{cond}, $(iii)$, the above inequality must be equality and hence $E'_1$ is contained in a plane.
\end{pf}

We are now in a position to show that the subspace $W_{\phi}$ intersects the subspace $e'_0 H^0 (\OO_X (\KX-E'_0))$ nontrivially.
\begin{lem}\label{l:E'}
The intersection $W_{\phi} \bigcap e'_0 H^0 (\OO_X (\KX-E'_0)) $ is nonzero.
\end{lem}
\begin{pf}
Let $\sigma_0$ be a generator of $H^0 (\OO_X (\KX-E'))$. The multiplication of the exact sequence \eqref{GTeta1} by $\sigma_0$ implies that all global sections
in $\GG_{\Pi} (\KX-E')$ are of the form $\tau_{\Pi} \wedge H^0 (\TT'_{\eta})$. Using this observation for the global section $(e'_0 \alpha'^{(2)}_{\eta} (\phi,\phi') -m(\phi,\phi')h)$ used in defining $\GG_{\Pi}$, we obtain the equation
$$
\sigma_0 (e'_0 \alpha'^{(2)}_{\eta} (\phi,\phi') -m(\phi,\phi')h)=\tau_{\Pi} \wedge \alpha'_{\eta} (\psi),
$$
for a unique nonzero $\psi \in \HKX$. Taking the exterior product (on the right) with $\alpha'_{\eta} (\psi)$, we obtain
$$
\sigma_0 (e'_0 \alpha'^{(3)}_{\eta} (\phi,\phi',\psi) -m(\phi,\phi')\psi)=0
$$
or, equivalently, 
\begin{equation}\label{psi-rel}
e'_0 \alpha'^{(3)}_{\eta} (\phi,\phi',\psi)=m(\phi,\phi')\psi.
\end{equation}
This means that $m(\phi,\phi')\psi$ vanishes on the divisor $E'_0 =(e'_0=0)$. The fact that the linear system $|Im(m)|$ has no fixed part, see Lemma \ref{l:l-modif}, insures that for a generally chosen $\phi,\phi' \in \HKX$, the section $m(\phi,\phi')$ does not vanish on any component of $E'_0$. Hence $\psi$ must be a multiple of $e'_0$, i.e. $\psi \in e'_0 H^0 (\OO_X (\KX-E'_0))$. Furthermore, we claim that the equation \eqref{psi-rel} implies that $m(\phi, \psi)=0$.
Indeed, substituting for the left hand side the Koszul coboundary expression in \eqref{cobound-formula} we deduce
$$
\phi m(\phi',\psi) - \phi' m(\phi, \psi) =0.
$$ 
Since $\phi$ and $\phi'$ have no components in common, the above can only occur for $m(\phi, \psi)=m(\phi',\psi)=0$. Hence $\psi$ is a nonzero element in the intersection $e'_0 H^0 (\OO_X (\KX-E'_0))\bigcap W_{\phi}$.
\end{pf}

Let $x$ be a nonzero element of $W_{\phi} \bigcap e'_0 H^0 (\OO_X (\KX-E'_0))$. Then the coboundary relation \eqref{cobound-formula} for $\phi,x$, and an arbitrary
$\phi'\in \HKX$ reads
$$
e'_0 \alpha'^{(3)}_{\eta} (\phi,\phi', x) =xm(\phi,\phi')+\phi m(\phi',x).
$$
Since $x$ is a multiple of $e'_0$ and $\phi$ is general and hence does not vanish on any component of $E'_0=(e'_0 =0)$, the above relation tells us that
$m(\phi',x)\in H^0 (\OO_X (\KX-E'_1))$ is a multiple of $e'_0$ as well. So, upon fixing a generator $\sigma_0$ of  the one dimensional space $H^0 (\OO_X (\KX-E'))$, we obtain a linear function $t_x: \HKX \longrightarrow \CC$ such that
$$
 m(\phi',x)=t_x(\phi')e'_0 \sigma_0, \,\forall \phi' \in \HKX.
$$ 
From this it follows 
\begin{equation}\label{dimWx}
dim (W_x) \geq p_g -1, 
\end{equation}
where $W_x:= ker(m(\bullet,x))$.

For any $\psi,\psi' \in W_x$ the coboundary formula \eqref{cobound-formula} gives
$$
e'_0 \alpha'^{(3)}_{\eta} (\psi,\psi', x) =xm(\psi,\psi') =e'_0 p_x m(\psi,\psi'),
$$
where $ p_x \in H^0 (\OO_X (\KX-E'_0))$ is such that $x=e'_0 p_x$. Hence the formula
$$
\alpha'^{(3)}_{\eta} (\psi,\psi', x) =p_x m(\psi,\psi'), \forall \psi,\psi' \in W_x.
$$
This tells us that the global sections $\alpha'_{\eta} (\psi)$, $\alpha'_{\eta} (\psi')$, $\alpha'_{\eta} (x)$ are algebraically dependent along the divisor
$D_x =(p_x =0) \in |\KX -E'_0|$, for any $\psi,\psi' \in W_x$. This in turn implies that $\alpha'_{\eta} (\psi)$, $\alpha'_{\eta} (\psi')$, $\alpha'_{\eta} (\psi'')$
are algebraically dependent along $D_x $, for all $\psi,\psi', \psi'' \in W_x$, or, equivalently, $\alpha'^{(3)}_{\eta} (\psi,\psi', \psi'')=p_x \sigma(\psi,\psi', \psi'')$, for a unique $\sigma(\psi,\psi', \psi'') \in H^0 (\OO_X (\KX-E'_1))$. This gives a linear map
$$
\mbox{$\sigma: \bigwedge^3 W_x \longrightarrow H^0 (\OO_X (\KX-E'_1))$}
$$
determined by the rule $ \psi \wedge \psi' \wedge \psi'' \mapsto \sigma(\psi,\psi', \psi'')$.

Fix $\psi \in W_x$ so that $C_{\psi} =(\psi=0)$ is irreducible\footnote{this is possible since the linear subsystem $|W_x| \subset |\KX|$ has at most one base point in view of the estimate \eqref{dimWx}.} and consider
the linear map
$$
\mbox{$\sigma(\psi,\bullet, \bullet): \bigwedge^2 W_x \longrightarrow H^0 (\OO_X (\KX-E'_1))$.}
$$
From the relation 
$\alpha'^{(3)}_{\eta} (\psi,\psi', \psi'')=p_x \sigma(\psi,\psi', \psi'')$ and Proposition \ref{pro:modif-eta}, 5), $(ii)$, the kernel of this map misses all decomposable tensors  $\psi' \wedge \psi''$ in $\bigwedge^2( W_x /\CC\psi)$. Hence the estimate
$$
p_g -3=  h^0 (\OO_X (\KX-E'_1))\geq dim (Im (\sigma(\psi,\bullet,\bullet))) \geq 2(dim( W_x /\CC\psi) -2) +1\geq 2(p_g-4)+1=2p_g -7,
$$
where the third inequality comes from \eqref{dimWx}.
The inequality implies that $p_g=4$. But then $X$ is a smooth quintic in $\PP^3$ and the condition \eqref{cond}, $(iii)$, implies that the canonical linear system $|\KX|$ has no reducible members\footnote{see a more general argument on p. \pageref{quintic} ruling out all quintic surfaces.}. This forces $E'=E'_1=\KX$ and contradicts $(\KX-E'_1)^2 >0$.

\vspace{0.2cm} 
{\it Case 2: $(\KX-E'_1)^2 =0$.} We claim the following.
\begin{lem}\label{l:quintic}
$\OO_X (\KX-E'_1)=\OO_X$ and $p_g=4$.
\end{lem}

\begin{pf}
Assume $\KX-E'_1 \neq 0$. Then the linear system
$|\KX-E'_1|$ is base point free and defines
a surjective morphism
$f:X \longrightarrow \PP^1$ with connected fibres such that $\KX-E'_1=aF$, where $F$ is the class of a fibre of $f$ and  $a$ is a positive integer.

 Our argument turns on the property of the sheaf $\TT'_{\eta}$ recorded in
Proposition \ref{pro:modif-eta}, 5), $(ii)$: \linebreak
$\alpha'^{(3)}_{\eta} (\phi,\phi',\phi'')\neq 0$ as long as the linear system defined by $\CC\{\phi,\phi',\phi''\}$ is $2$-dimensional and has at most $0$-dimensional base locus. This implies in particular that for a general $\phi \in \HKX$, the sheaf $\FF'_{[\phi]}$ can not have a subsheaf of rank $1$ with at least $2$-dimensional space of global sections. Indeed, otherwise
$\FF'_{[\phi]}$ has two linearly independent global sections $f_{\phi'}$ and $f_{\phi''}$ which are proportional, i.e. $f_{\phi'} \wedge f_{\phi''} =0$ as a global section of $\det(\FF'_{[\phi]})$, and then $0=f_{\phi'} \wedge f_{\phi'} =\alpha'^{(3)}_{\eta} (\phi,\phi',\phi'')$.

With the above remark in mind, we start with $\FF'_{[\phi]}$, for a general $\phi \in \HKX$, and recall its modification $\FF''_{[\phi]}$ along the curve $C_{\phi} =(\phi =0)$. This is given by the exact sequence 
$$
\xymatrix@R=12pt@C=12pt{
0\ar[r]&\FF''_{[\phi]} \ar[r]& \FF'_{[\phi]}  \ar[r]&\OO_{C_{\phi}} (\KX-E')\otimes \OO_{C_{\phi}} (B_{\phi})\ar[r]&0}, 
$$
see Proposition \ref{pro:F''phi}, \eqref {modifF'phi2}.
According to that proposition, $\FF''_{\phi}$ is a locally free sheaf of rank $2$ subject to:
$$
c_1 ( \FF''_{\phi})=\KX-E'=\KX-E'_1 -E'_0 =aF -E'_0.
$$
Furthermore, Proposition \ref{pro:F''phi}, 2), tells us that $H^0 (\FF''_{\phi})$ contains a subspace isomorphic to $W_{\phi} /\CC\phi$, where $W_{\phi}$ is the kernel of the homomorphism
$$
m(\phi,\bullet): \HKX \longrightarrow H^0 (\OO_X (\KX-E'_1)).
$$
 Hence the estimate
\begin{equation}\label{estimh0}
h^0(\FF''_{\phi}) \geq dim (W_{\phi} /\CC\phi) \geq p_g-1 - h^0 (\OO_X (\KX-E'_1))\geq p_g-1 - (p_g -3)=2,
\end{equation}
 where the last inequality uses the upper bound
\begin{equation}\label{estim}
 h^0 (\OO_X (\KX-E'_1))\leq p_g -3
\end{equation}
imposed by the condition that $E'_1$ can not be a line. We claim that the equality must hold. Indeed, if the inequality above is strict, then the estimate \eqref{estimh0} implies
\begin{equation}\label{h0F''}
h^0(\FF''_{\phi}) \geq 3.
\end{equation}
If a general global section of $\FF''_{\phi}$ has at most $0$-dimensional zero locus, then the Koszul sequence of such a section has the form
$$
\xymatrix@R=12pt@C=12pt{
0\ar[r]&\OO_X \ar[r]&\FF''_{\phi}\ar[r]& {\cal I}_Z (aF-E'_0) \ar[r]&0,
}
$$
where $Z$ is subscheme of dimension at most $0$ and ${\cal I}_Z$ is its ideal sheaf. 
Hence $(aF-E'_0)$ is effective and the divisor $E'_0$ must be contained in the fibers of $f$. Therefore, restricting the above sequence to a smooth fibre $F_t$ of $f$ disjoint from $Z$ gives
\begin{equation}\label{FF''phionFt}
\xymatrix@R=12pt@C=12pt{
0\ar[r]&\OO_{F_t} \ar[r]&\FF''_{\phi} \otimes \OO_{F_t}  \ar[r]&\OO_{F_t}   \ar[r]&0.
}
\end{equation}
Hence $h^0 (\FF''_{\phi} \otimes \OO_{F_t} )\leq 2$. This together with \eqref{h0F''} imply
$$
h^0 (\FF''_{\phi} (-F)) \geq h^0(\FF''_{\phi} )- h^0 (\FF''_{\phi}\otimes \OO_{F_t} ) \geq 3-2=1.
$$
From this and the inclusion $\FF''_{\phi} \hookrightarrow \FF'_{\phi}$, it follows that $\FF'_{\phi} \otimes \OO_X (-F)$ has a nonzero global section. Equivalently, we have a monomorphism $\OO_X (F) \longrightarrow \FF'_{\phi}$. This provides a a subsheaf of rank $1$ of $\FF'_{\phi}$ with the $2$-dimensional space of global sections and this, as we explained at the beginning of the proof, contradicts
Proposition \ref{pro:modif-eta}, 5), $(ii)$.
So we may assume that all nonzero global sections of $\FF''_{\phi}$  have  $1$-dimensional zero locus. Denote by $T$ the rational equivalence class of the divisorial part of a nonzero global section of $\FF_{\phi}$ . Then we have
a monomorphism
$\OO_X (T) \longrightarrow \FF''_{\phi}$ and $h^0 (\OO_X (T)) =1$. This means that all global sections of $\FF''_{\phi}$ vanish along an effective nonzero divisor $T$. Hence $h^0 (\FF''_{\phi} (-T))= h^0 (\FF''_{\phi}) \geq 3$. Furthermore, a general global section of $\FF''_{\phi} (-T)$ has at most $0$-dimensional zero locus. Hence the Koszul sequence of such a section has the form
\begin{equation}\label{KoszF'}
\xymatrix@R=12pt@C=12pt{
0\ar[r]&\OO_X \ar[r]&\FF''_{\phi} (-T) \ar[r]& {\cal I}_Z (aF-2T) \ar[r]&0,
}
\end{equation}
where $Z$ is a subscheme of codimension at least $2$ and ${\cal I}_Z$ is its ideal sheaf. From $h^0 ({\cal I}_Z (aF-2T)) =h^0 (\FF'_{\phi} (-T))-1 \geq 3-1=2$, it follows that  $aF-2T=N +M$, where $N$ (resp. $M$) is the fixed (resp. moving) part of the linear system $|aF-2T|$. Again, the divisors $T, N, M$ are contained in the fibres of the morphism $f:X \longrightarrow \PP^1$. 
Thus the restriction of \eqref{KoszF'} to a general smooth fibre $F_t$ of $f$ has the form
$$
\xymatrix@R=12pt@C=12pt{
0\ar[r]&\OO_{F_t} \ar[r]&\FF''_{\phi} \otimes \OO_{F_t}  \ar[r]&\OO_{F_t} \ar[r]&0
}
$$
and we deduce $h^0 (\FF''_{\phi}\otimes \OO_{F_t}) \leq 2$ with the conclusion  
$h^0 (\FF''_{\phi} (-F)) \geq 1$. Thus there are global sections of $\FF''_{\phi}$
with a moving divisorial part which contradicts our assumptions.

\vspace{0.2cm}
We now have the equality $h^0 (\OO_X (\KX-E'_1))= p_g -3$. This means that the linear span $P_{E'_1}$ of $E'_1$ is a plane. From this it follows that $E'_1$ has no multiple components, since for such a component the plane $P_{E'_1}$ would be the embedded tangent plane of $X$ along that component and this is impossible in view of the finiteness of the Gauss map of $X\subset \PP(\HKX^{\ast})$, see \cite{Z}.

Next we show the equality $E'_1 =E'$.
For this we use that $\KX-E'=K_X -E'_1 -E'_0 =aF-E'_0$. Hence all irreducible components of $E'_0$ must be contained in the fibres of $f$. In particular, $C^2 \leq 0$ for any irreducible component of $E'_0$. We recall that the support of $E'$ is $E'_1$, see Lemma \ref{l:l-modif}.
Therefore, the irreducible components of $E'_0$ are plane curves. But then the normal sheaf $\OO_C (C) =\OO_C ((d_C -4)\KX)$, where $d_C =\KX \cdot C$ is the degree of $C$. Hence the inequality
$$
0\geq C^2 =d_C (d_C -4)
$$
implying $d_C =3$ or $4$ (we are using here the condition \eqref{cond}, $(iii)$, excluding lines and conics). But the formula
$$
\KX =(\KX-E'_1) +E'_1=aF +E'_1
$$
shows that none of such curves can exist on $X$. Indeed, if there is $C$ in $E'_1$ with $d_C=3$ (i.e., $\KX \cdot C=-C^2=3$) and which is also part of a fibre of $f$, the formula for $\KX$ gives
$$
3=\KX\cdot C=E'_1 \cdot C=C^2 + (E'_1 -C) \cdot C=-3+(E'_1 -C) \cdot C.
$$
Hence $(E'_1 -C) \cdot C =6$ and this means that the component complementary to $C$ is a conic which is impossible in view of the  condition \eqref{cond}, $(iii)$. Running the same argument for $d_C=4$, i.e., $\KX \cdot C=4$, $C^2=0$, we obtain
$$
 4=\KX\cdot C=E'_1 \cdot C=C^2 + (E'_1 -C) \cdot C=(E'_1 -C) \cdot C.
$$
This means that the component complementary to $C$ is a line and this again is in contradiction with \eqref{cond}, $(iii)$.

With the understanding that $E'=E'_1$ and $\KX-E'=\KX-E'_1=aF$, we return to the exact sequence \eqref{GTeta1} to obtain
$$
\xymatrix@R=12pt@C=12pt{
0\ar[r]&\OO_{X} (-2aF) \ar[r]&\TT'_{\eta}(-aF)\ar[r]& {\cal G}_{\Pi}  \ar[r]&0.
}
$$
Tensoring with $\OO_X ((a-1)F)$ we deduce
$$
h^0 (\TT'_{\eta}(-F)) \geq h^0 ({\cal G}_{\Pi}((a-1)F))-h^1 (\OO_{X} (-(a+1)F)) \geq 3a - a=2a \geq 2.
$$
Hence a nonzero morphism $\OO_X (F) \longrightarrow \TT'_{\eta}$ which is in contradiction with Proposition \ref{pro:modif-eta}, 5), $(iii)$. This completes the proof of the assertion $\OO_X (\KX-E'_1) =\OO_X$.

We now turn to the assertion $p_g=4$. From $\KX-E'_1 =0$ and the effectiveness of
$\KX -E'=\KX-E'_1 -E'_0=-E'_0$ it follows that $E'_1=E'$. The cochain $m$ now is a linear map
$$
\mbox{$\bigwedge^2 \HKX \longrightarrow H^0 (\OO_X (\KX-E'_1))= H^0 (\OO_X).$}
$$
In particular, the kernel $W_{\phi}=ker(m(\phi,\bullet))$ is a subspace of $\HKX$ of codimension at most $1$. Hence the inequality \eqref{Wphi1} reads
$$
1=H^0 (\OO_X (\KX-E')) \geq 2 dim (W_{\phi}) -5 \geq 2(p_g-1)-5=2p_g -7
$$
or, equivalently, $p_g \leq 4$. The condition \eqref{cond}, $(i)$, saying that $\KX$ is very ample gives the equality $p_g = 4$. 
\end{pf}

From Lemma \ref{l:quintic} it follows that $X$ is a smooth quintic surface in $\PP^3$. Of course, for surfaces of degree $d\geq 4$ in $\PP^3 $ the injectivity of the cup-product 
\eqref{cp} is well known and is a consequence of Macaulay's theorem, see e.g \cite{V}, Theorem 18.19 and Remarque 18.30. Let us see, however, how the situation of a quintic can be ruled out by using the extension construction. \label{quintic}

The problem here is easy to detect: our parametrization $\alpha: \HKX \longrightarrow H^0 (\TT_{\xi})$ has the property that the global sections 
$(\phi' \alpha(\phi)- \phi \alpha(\phi'))$ of $H^0 (\TT_{\xi}(\KX))$, for all
$\phi, \phi' \in \HKX$, lie in the image of the injective homomorphism
$H^0 (\Omega_X (\KX)) \longrightarrow H^0 (\TT_{\xi}(\KX))$ induced by the monomorphism of our extension sequence \eqref{ext} tensored with $\OO_X (\KX)$.
But it is immediately seen from the normal sequence of $X \subset \PP^3$ that
$H^0 (\Omega_X (\KX))=0$. Hence  the relation
$$
\phi' \alpha(\phi)- \phi \alpha(\phi')=0,\, \forall \phi, \phi' \in \HKX.
$$
 But for $\phi, \phi'$ linearly independent and 
$\phi$ with $C_{\phi} =(\phi=0)$ irreducible the above relation tells us that
$\alpha(\phi)$ must vanish on $C_{\phi}$ or, equivalently, it has the form
$\alpha(\phi)=\phi \sigma$, where $\sigma$ is a nonzero global section of 
$\TT_{\xi}(-\KX)$. This in turn means that the extension sequence \eqref{ext}
splits, i.e. $\xi=0$.
 
\end{pf}

Denote by $Z^1_{\Pi}$ the $1$-dimensional part of the zero locus $Z_{\Pi}$ of the global section $\tau_{\Pi}$ defined in \eqref{tauPi-formula}. By Lemma \ref{l:zeroloc-tauPi} the divisor $Z^1_{\Pi} \neq 0$. We wish to understand
this divisor.
\begin{lem}\label{l:Z1Piproper}
Let $\Gamma_{\Pi}=(\alpha'^{(3)}_{\eta} (\phi,\phi',\phi'') =0)$, for some basis $\{\phi,\phi',\phi''\}$ of a $3$-dimensional subspace $\Pi$ of $\HKX$. Then, for a general $\Pi$, the divisor $Z^1_{\Pi} =\Gamma_{\Pi}$ and the global section $\tau_{\Pi}$ has the form
\begin{equation}\label{tauPi-formula1}
\tau_{\Pi} =\alpha'^{(3)}_{\eta} (\phi,\phi',\phi'') \tau,
\end{equation}
where $\tau$ is a nonzero global section of $\TT'_{\eta} (-\KX +E'_0)$.
\end{lem}
\begin{pf}
Recall that for a $3$-dimensional subspace $\Pi$ of $\HKX$ the section $\tau_{\Pi}$ is defined by the formula
\begin{equation}\label{tauPi-recallform}
\tau_{\Pi} =m(\phi, \phi')\alpha'_{\eta} (\phi'') -m(\phi, \phi'')\alpha'_{\eta} (\phi') +m(\phi', \phi'')\alpha'_{\eta} (\phi)
\end{equation}
for some basis $\{\phi,\phi',\phi''\}$ of $\Pi$. Taking the exterior product of 
$\tau_{\Pi}$ with $\alpha^{(2)} (\psi, \psi')$, for $\psi, \psi' \in \Pi$ we deduce that the zero locus $Z_{\Pi}=( \tau_{\Pi}=0)$ is contained in the base locus of the linear system generated by the divisors
$$
M(\psi, \psi')+\Gamma_{\Pi} :=\Big(m(\psi, \psi')\alpha'^{(3)}_{\eta} (\phi,\phi',\phi'')=0\Big), \forall \psi, \psi'\in \Pi,
$$
where $M(\psi, \psi')=(m(\psi, \psi')=0)$. In particular, $Z^1_{\Pi}$ is contained in the fixed part of the linear system
$$
\Gamma_{\Pi} +|M(\psi, \psi')|_{ \psi, \psi'\in \Pi }.
$$
By Corollary \ref{cor:m-3dim}, for a general choice of $\Pi$ the linear system $|M(\psi, \psi')|_{ \psi, \psi'\in \Pi }$ is fixed part free.\footnote{Corollary \ref{cor:m-3dim} requires the assumption $p_g \geq 5$. This is now superfluous, since in the end of the proof of Lemma \ref{l:zeroloc-tauPi} the case $p_g=4$ has been ruled out.} Hence $Z^1_{\Pi}$ is a component of $\Gamma_{\Pi}$. By Proposition \ref{pro:F'phi-gen-nopencil}, 3), the latter divisor is irreducible, for a general
$\Pi$. Therefore, $Z^1_{\Pi} =\Gamma_{\Pi}$, for a general $\Pi$, and we deduce 
$$
\tau_{\Pi}=\alpha'^{(3)}_{\eta} (\phi,\phi',\phi'') \tau,
$$
where $\tau$ is a nonzero global section  of $\TT'_{\eta} (K_X -E'_1) \otimes \OO_X (-2\KX +E')=\TT'_{\eta} (-\KX +E'_0)$.
\end{pf}

We are now ready to complete our argument. The first consequences of the lemma above are

$\bullet$ the global section $\tau$ in \eqref{tauPi-formula1} is unique,

$\bullet$ $h^0 (\OO_X (\KX -E'_0)) =1$.

The second statement is immediate from Proposition \ref{pro:modif-eta}, 5), $(iii)$. For the first one, we observe that under the homomorphism
\begin{equation}\label{p'0}
 p'_0 :H^0 (\TT'_{\eta} (-(\KX -E'_0))) \longrightarrow  H^0 (\OO_X (E'_0))
\end{equation}
induced by the epimorphism in the extension sequence 
\begin{equation}\label{ext-eta'-end}
\xymatrix@R=12pt@C=12pt{
0\ar[r]& {\cal P}\ar[r]&\TT'_{\eta} \ar[r]&\OO_X (\KX) \ar[r]&0,
}
\end{equation}
tensored with $\OO_X (-(\KX -E'_0))$ (see \eqref{ext-eta'} for that sequence),
the global section $\tau \in H^0 (\TT'_{\eta} (-(\KX -E'_0)))$ goes to $e'_0$, the section defining the divisor $E'_0$. This follows from applying the homomorphism
$$
p':H^0 (\TT'_{\eta} ((\KX -E'_1))) \longrightarrow  H^0 (\OO_X (2\KX -E'_1)),
$$
induced by the epimorphism in \eqref{ext-eta'-end} tensored with $\OO_X (\KX -E'_1)$, to the equation \eqref{tauPi-formula1}:
\begin{equation}\label{takingp'}
p'(\tau_{\Pi})=p'(\alpha'^{(3)}_{\eta} (\phi,\phi',\phi'') \tau) = \alpha'^{(3)}_{\eta} (\phi,\phi',\phi'') p'_0 (\tau),
\end{equation}
where $ p'_0$ is the homomorphism in \eqref{p'0}. On the other hand $p'(\tau_{\Pi})$ computed from the defining formula \eqref{tauPi-recallform} gives
$$
\begin{gathered}
p'(\tau_{\Pi})= p'\Big(m(\phi, \phi')\alpha'_{\eta} (\phi'') -m(\phi, \phi'')\alpha'_{\eta} (\phi') +m(\phi', \phi'')\alpha'_{\eta} (\phi)\Big)=\\
m(\phi, \phi') \phi'' -m(\phi, \phi'')\phi' +m(\phi', \phi'') \phi =e'_0 \alpha'^{(3)}_{\eta} (\phi,\phi',\phi''),
\end{gathered} 
$$
where the last equality is the Koszul coboundary relation \eqref{Kosz-coboune'0}. Putting this together with \eqref{takingp'} we deduce
$$
p'_0 (\tau)=e'_0.
$$
The uniqueness of $\tau$ now follows from the injectivity of $p'_0$, due to the fact that $ker (p'_0) \cong H^0 ({\cal P}(-(\KX-E'_0))) =0$.

The uniqueness of $\tau$ in the formula \eqref{tauPi-formula1} implies that
the formula holds for all  triples $\phi, \phi', \phi'' \in \HKX$.
Furthermore, taking the exterior product with $\alpha'^{(2)}_{\eta} (\phi,\phi')$ on the both sides of \eqref{tauPi-formula1} we obtain
$$
m(\phi,\phi') \alpha'^{(3)}_{\eta} (\phi,\phi',\phi'') =  \alpha'^{(3)}_{\eta} (\phi,\phi',\phi'') \big(\tau \wedge \alpha'^{(2)}_{\eta} (\phi,\phi') \big).
$$
Hence the formula for $m$:
$$
m(\phi,\phi')= \tau \wedge \alpha'^{(2)}_{\eta} (\phi,\phi') , \,\,\forall \phi,\phi' \in \HKX.
$$
Fix a nonzero global section $\sigma'_0$ of $\OO_X (\KX- E'_0)$ and multiply by it the above equation to obtain
\begin{equation}\label{m-a3}
\sigma'_0 m(\phi,\phi')= (\sigma'_0\tau) \wedge \alpha'^{(2)}_{\eta} (\phi,\phi') =\alpha'^{(3)}_{\eta} ( \phi_{0},\phi,\phi'),
\end{equation}
where $\phi_{0} =e'_0 \sigma'_0$ and the second equality uses the identity
$\sigma'_0\tau =\alpha'_{\eta} ( \sigma'_0 e'_0) =\alpha'_{\eta} (\phi_{0})$. But 
Lemma \ref{l:m(phi,.)} asserts that the kernel of $m(\phi,\bullet)$ is at least
$3$-dimensional. Thus we are sure to find $\phi'_0 \in \HKX$, linearly independent
of $\phi_0$ and $\phi$ such that $m(\phi,\phi'_0)=0$. Hence the equation \eqref{m-a3} implies
 $\alpha'^{(3)}_{\eta} ( \phi_{0},\phi,\phi'_0) =0$ and, for general $\phi \in \HKX$, this in contradiction with Proposition \ref{pro:modif-eta}, 5), $(ii)$.
The proof of Theorem \ref{th-tech} is now completed.

\vspace{1cm}
\begin{flushright}         
Universit\'e d'Angers\\
D\'epartement de Math\'ematiques
\\
2, boulevard Lavoisier\\
49045 ANGERS Cedex 01 \\
FRANCE\\
{\em{E-mail addres:}} reider@univ-angers.fr
\end{flushright} 


\begin{thebibliography}{9999}
\bibitem[Ga-Z]{Ga-Z}
U. Garra, F. Zucconi, Very ampleness and the Infinitesimal Torelli problem,
Math.Z. (2008),n°1, 31-46.
\bibitem[G]{G}
 Griffiths, Ph.(ed),  Topics in transcendental algebraic geometry,
 Ann.of Math.Stud., 106, Princeton Univ.Press, (1984).
\bibitem[G-H]{GH}
 Griffiths, Ph., Harris, J., Principles of algebraic geometry, New-York, 
 Wiley, 1978.  
\bibitem[G-S]{GS}
 Griffiths, Ph., Schmid, W., Recent development in Hodge theory:
 a discussion of techniques and results. In:
 Disctete subgroups of Lie groups and applications to moduli.
 Intern.Colloq., Bombay, 1973, Oxford Univ.Press, Bombay, 1975, 31-127. 
\bibitem[F]{Fl}
 Flenner, H., The Infinitesimal Torelli problem for zero sets of sections of vector
 bundles, Math.Z., 193, $n^{\circ}2$, (1986), 307-322.
\bibitem[Gr]{Gr}
 Green, M., Koszul cohomology and the geometry of projective varieties, II.
 J.Diff.Geom.20, $n^{\circ}1$ (1984) 279-289.
\bibitem[O-S-S]{[O-S-S]} Okonek, Ch., Schneider, M., Spindler H., Vector bundles on complex projective spaces, Birkh\"{a}user, 1988.
\bibitem[M]{M}
Miyaoka, Y., The Chern classes and Kodaira dimension of a minimal variety. In:
Algebraic geometry, Sendai, 1985, 449-476, Adv. Stud. Pure Math., North-Holland, Amsterdam, 1987.
\bibitem[Pa]{Pa}
 Pardini, R., Infinitesimal Torelli and Abelian covers of algebraic surfaces.
 In: Problems in the theory of surfaces and their classification (Cortona, 1988),
 247-257, Sympos.Math., XXXII, Academic Press, London,1991.  
\bibitem[Pe]{Pe}
 Peters, C., The local Torelli theorem II. Cyclic branched coverings.
 Ann.Scuola Norm.Sup.Pisa, Cl.Sci., (4)3, (1976), 321-340.
\bibitem[R]{R}
Reider, I., Vector bundles of rank 2 and linear systems on algebraic surfaces, Ann. of Math., 127 (1988), 309-316.
\bibitem[R1]{R1}
Reider, I., On the Infinitesimal Torelli Theorem for certain irregular surfaces of general type, Math.Ann. 280, 285-302 (1988).
\bibitem[R2]{R2}
Reider, I., Infinitesimal Torelli Theorem for double coverings of surfaces of general type, J. Algebraic Geometry, 14 (2005),691-704.
\bibitem[R3]{R3}
Reider, I., Nonabelian Jacobian of smooth projective surfaces, J. Differential Geometry, 74, 2006, 425-505.
\bibitem[R4]{R4}
Reider, I., Nonabelian Jacobian of smooth projective surfaces - a survey,
Science China, Math., Vol56, No 1, 1-42, 2013.
\bibitem[S]{S}
Siu, Y.-T., Curvature of the Weil-Petersson Metric in the Moduli space of compact Kahler-Einstein manifolds of negative first Chern class, in A. Howard et al. eds), Contributions to Several Complex Variables, 1986, 261-298.
\bibitem[T]{T}
Thomas, R.P., Derived categories for the working mathematician, 
arXiv:math/0001045v2[mathAG], 2001.
\bibitem[Ts]{Ts}
Tsuji, H., Stability of tangent bundles of minimal algebraic varieties, Topology
Vol 27, No. 4, 429-442, 1988.
\bibitem[V]{V}
Voisin, C., Th\'eorie de Hodge et g\'eom\'etrie alg\'ebrique complexe, Cours sp\'ecialis\'es v.10, SMF, 2002.
\bibitem[Z]{Z}
Zak, F.L., Projections of algebraic varieties, Math. Sb. (N.S.) 116 (1981), 593-602.

\end{thebibliography}
\end{document}